%% file: fundamentaalgroup.tex
\documentclass[11pt]{amsart}
\input{sdpStyle.inc}

\input{sdpDefinitions.inc}

\begin{document}
\title[Explicit relations and factors with uncountable fundamental
group]{Explicit examples of equivalence relations and II$_1$ factors
  with prescribed fundamental group and outer automorphism group}
\stepcounter{footnote}
\author{Steven Deprez}
\thanks{Research assistant of the Research Foundation -- Flanders
  (FWO) (till august 2011)}
\thanks{Partially supported by ERC grant VNALG-200749}
\thanks{Postdoc at the university of Copenhagen (from September 2011)}
\thanks{Partially supported by [Uffe's grant]}
\address{Department of mathematics, K.U.Leuven, Celestijnenlaan 200B, B-3001 Leuven (Belgium)}
\email{steven.deprez@wis.kuleuven.be}

\begin{abstract}
  In this paper we give a number of explicit constructions for II$_1$
  factors and II$_1$ equivalence relations that have prescribed
  fundamental group and outer automorphism group. We
  construct factors and relations that have uncountable fundamental
  group different from $\IRpos$. In fact, given any II$_1$ equivalence
  relation, we construct a II$_1$ factor with the same fundamental group.

  Given any locally compact unimodular second countable group
  $G$, our construction gives a II$_1$ equivalence relation $\RelR$ whose outer
  automorphism group is $G$. The same construction does not give a
  II$_1$ factor with $G$ as outer automorphism group, but when $G$ is
  a compact group or if $G=\SL^{\pm}_n\IR=\{g\in\GL_n\IR\mid \det(g)=\pm1\}$, then
  we still find a type II$_1$ factor whose outer automorphism group is
  $G$.
\end{abstract}

\maketitle

\tableofcontents
\section*{Introduction}
%
%
Type II$_1$ factors can be constructed in many different ways. For example
as group von Neumann algebras or using the group-measure space
construction. Two von Neumann algebras that are constructed in
different ways can nevertheless be isomorphic. For example, all the
group von Neumann algebras of ICC amenable groups are isomorphic
\cite{Connes:Injective}. So a number of invariants for type II$_1$
factors were introduced. Probably the two most natural among them are
the fundamental group and the outer automorphism group. Unfortunately,
both are very hard to compute. For type II$_1$ equivalence relations,
there are similar notions of fundamental group and outer automorphism
group.

Our knowledge about fundamental groups and outer automorphism groups of type
II$_1$ factors and equivalence relations advanced a lot in recent
years. The recent development of Popa's deformation/rigidity theory
allows us to study questions like ``which groups appear as the
fundamental group (or outer automorphism group) of a II$_1$ factor (or
equivalence relation).'' When we study this question, we will only be
interested in separable II$_1$ factors, i.e.\ II$_1$ factors that are
represented on a separable Hilbert space.

In this article I introduce a general method to create explicit
examples of type II$_1$ factors and equivalence relations with
prescribed fundamental group or prescribed outer automorphism
group. This general method provides explicit examples in
cases where we only knew pure existence results. For example, I give
explicit examples of type II$_1$ factors with uncountable fundamental
group different from $\IRpos$. Existence of such factors was already
proven in \cite{PopaVaes:ActionsOfFinfty,PopaVaes:FundamentalGroupGeneral}.
Moreover, I give an explicit
construction for type II$_1$ factors with prescribed second countable
compact outer automorphism group.

But I also give essentially new
examples. I show that every second countable unimodular locally
compact group is the outer automorphism group of a II$_1$ equivalence
relation. This result does not extend immediately to II$_1$ factors,
but I show that the groups $\SL^{\pm}_n\IR=\{g\in\GL_n\IR\mid
\det(g)=\pm1\}$ are outer automorphism groups of type II$_1$ factors.
Using a different style of example, I show that every subgroup of
$\IRpos$ that appears as the fundamental group of a type II$_1$
equivalence relation also appears as the fundamental group of a II$_1$
factor.

The fundamental group was introduced by Murray and von Neumann in
\cite{MvN:ROO4}. The fundamental group of a type II$_1$ factor is not
related to the fundamental group of a topological space. Instead it is
a group $\Fundg(M)\subset\IRpos$ of positive real numbers. Murray and
von Neumann could compute the fundamental group in only one case. They
showed that the fundamental group of the hyperfinite type II$_1$
factor $R$ is $\IRpos$ itself. They raised the question:
``What subgroups $\cF\subset\IRpos$ appear as the fundamental group of
a separable type II$_1$ factor?''

The outer automorphism group $\Out(M)$ is the quotient group of all
automorphisms of $M$ by the group of inner automorphisms, i.e.\ by the
group $\Inn(M)=\{\Ad_u\mid u\in \Unitary(M)\}$. Blattner shows that
the group of outer automorphisms of $R$ is huge: it contains an
automorphic copy of any second countable locally compact group
\cite{Blattner:AutomorphicGroupRepresentations}.
I study the question: ``which groups appear as the outer
automorphism group of a separable type II$_1$ factor?''

Similar notions of fundamental group and outer automorphism group are defined
for type II$_1$ equivalence relations.

The first significant progress on the questions raised above was
realized by Connes in
\cite{Connes:CountableFundamentalGroup}. He showed that the fundamental group
and the outer automorphism group of a property (T) factor are
countable. Examples of property (T) factors are the group von Neumann
algebras $\Lg(\Gamma)$ of property (T) ICC groups $\Gamma$. These were
the first
examples of II$_1$ factors with small fundamental group and outer
automorphism group. However, he could only show that the invariants
are countable, and even today we do not know the fundamental group nor
the outer automorphism group of any property (T) factor.

Popa's Deformation/rigidity theory allows, among many other things, to
compute outer automorphism groups and fundamental groups of type
II$_1$ factors and equivalence relations. As a first application to this theory, Popa
showed that the group von Neumann algebra $\Lg(\SL_2\IZ\ltimes\IZ^2)$
has trivial trivial fundamental group
\cite{Popa:Betti}, i.e.\ it contains only $1$. Later he showed that any
countable subgroup of $\IRpos$ appears as the fundamental group of
some type II$_1$ factor. Alternative constructions for the same result
were given in \cite{IoanaPetersonPopa:AFP, Houdayer:Fundg}. In
\cite{PopaVaes:ActionsOfFinfty,PopaVaes:FundamentalGroupGeneral}, Popa and Vaes
show that a large class of uncountable subgroups of $\IRpos$ appear as the
fundamental group of a type II$_1$ factor and of a II$_1$ equivalence
relation. However, they use a Baire category result to construct the
type II$_1$ factor. In this
paper I provide an explicit construction. This construction has
several additional advantages over the original construction. These
advantages will be detailed later on.

Moreover, I show that whenever
a subgroup of $\IRpos$ is the fundamental group of a type II$_1$
equivalence relation, then it is also the fundamental group of a type
II$_1$ factor, see theorem \ref{thm:SeqrelInSfactor}.

The first actual computation of an outer automorphism group was done
by Ioana, Peterson and Popa \cite{IoanaPetersonPopa:AFP}. They show that every
second countable compact abelian group appears as the outer
automorphism group of a type II$_1$ factor. Again, their construction
uses a Baire category argument. Later on Falgi\`eres and Vaes could
extend this result and show that every (not necessarily abelian)
second countable compact group appears as the outer automorphism group
of a type II$_1$ factor. Here I provide an explicit construction for
the latter.

For type II$_1$ equivalence relations, I even show
that every second countable unimodular locally compact group appears
as the outer automorphism group of a type II$_1$ equivalence
relation. I have not been able to extend this result to
outer automorphism groups of type II$_1$ factors. However, it is known
that several important subclasses of unimodular groups appear as the
outer automorphism group of a type II$_1$ factor. Vaes showed in
\cite{Vaes:Bimodules} that all countable discrete groups appear. As mentioned before, all
compact groups appear too. Moreover, I show that all
the groups of the form $\SL^{\pm}_n\IR=\{x\in\GL_n\IR\mid
\det(x)=\pm1\}$ with $n\in\IN$ appear as the outer automorphism group
of a type II$_1$ factor.

\section*{An overview of the construction}
\let\oldthetheorem=\thetheorem
\renewcommand\thetheorem{\Alph{theorem}}
Essentially all the results mentioned above use the same construction. 
This construction is very similar to the construction from
\cite{PopaVaes:ActionsOfFinfty,PopaVaes:FundamentalGroupGeneral}.

We start with two actions. Firstly, a ``general'' ergodic infinite
measure preserving action $\Lambda\actson (Y,\nu)$. Secondly, a
specific action $\Gamma\actson (X,\mu)$. For the moment we only assume
that this action is free, ergodic and probability measure
preserving. Suppose we also have a quotient morphism
$\pi:\Gamma\rightarrow\Lambda$. Then we can define a free infinite
measure preserving action $\Gamma\actson X\times Y$ by the formula
$g(x,y)=(g x,\pi(g)y)$. If we assume that $\ker\pi$ acts ergodically on
$X$, then this new action is also ergodic.

The orbit equivalence relation $\RelR=\Rel(\Gamma\actson X\times Y)$
is a type II$_\infty$ equivalence relation, and the von Neumann
algebra $M=\Lp^\infty(X\times Y)\rtimes\Gamma$ is a type II$_\infty$
factor. Let $U\subset X\times Y$ be a subset with finite measure. Then
the restricted equivalence relation $\RelR\restrict{U}$ is a type
II$_1$ equivalence relation and the following short exact sequence
describes the outer automorphism group and fundamental group of
$\RelR\restrict{U}$ in terms of the outer automorphism group of $\RelR$.
\[\begin{CD}
  1 @>>> \Out(\RelR\restrict U) @>>> \Out(\RelR) @>\Mod>>
  \Fundg(\RelR\restrict U) @>>> 1
\end{CD}\]
Every orbit equivalence $\Delta$ of $\Gamma\actson X\times Y$ with
itself scales the measure $\mu\times\nu$ by a constant
$\Mod(\Delta)$. Once we know $\Out(\RelR)$ and the group morphism
$\Mod:\Out(\RelR)\rightarrow \IRpos$, we also know the outer
automorphism group and the fundamental group of
$\RelR\restrict{U}$. The same relation holds between the outer
automorphism group of the II$_\infty$ factor $M$ and the outer
automorphism group and fundamental group of a finite corner $p M p$.

Let $\Delta\in\Autns(Y,\nu)$ be a non-singular automorphism that
commutes with the action of $\Lambda$. Then it is clear that
$\id\times\Delta$ is an orbit equivalence of $\Gamma\actson X\times
Y$. Because $\Gamma$ acts freely on $X$, such an orbit equivalence can
never be inner. We denote by $\Centr_{\Autns(Y,\nu)}(\Lambda)$ the group
of all non-singular automorphisms of $Y$ that commute with the action
of $\Lambda$. So the construction above embeds
$\Centr_{\Autns(Y,\nu)}(\Lambda)$ into $\Out(\RelR)$. 

Under strong conditions on the ``specific'' action $\Gamma\actson X$,
this one-to-one group morphism is actually an isomorphism (i.e.\ it is
also onto). In
\cite{PopaVaes:ActionsOfFinfty,PopaVaes:FundamentalGroupGeneral}, Popa
and Vaes give such a set of conditions on
$\Gamma\actson X$. However, they could not give any explicit example
of an action $\Gamma\actson X$ that satisfies their conditions.  They
used a Baire category argument to show the existence of such
actions. Theorem \ref{thm:intro:main:eqrel} below gives an alternative
set of conditions. This time, it is possible to give an explicit
example of an action $\Gamma\actson X$ that satisfies the conditions.
\begin{theorem}[see also theorem \ref{thm:proof:main:eqrel}]
  \label{thm:intro:main:eqrel}
  Let $\pi:\Gamma\rightarrow\Lambda$ be a quotient morphism and let
  $\Lambda\actson (Y,\nu)$ be an ergodic infinite measure preserving
  action. Suppose that the free, ergodic, probability measure
  preserving action $\Gamma\actson (X,\mu)$ satisfies the following
  conditions.
  \begin{itemize}
  \item The action $\ker\pi\actson X$ is $\Ufin$ cocycle superrigid
    and weakly mixing. See subsection \ref{subsect:prelim:CSR} for
    more details on cocycle superrigidity and weak mixing.
  \item There are no non-trivial group morphisms
    $\theta:\ker\pi\rightarrow \Lambda$.
  \item All the conjugations $\Delta$ of $\Gamma\actson X$ with itself are
    given by the formula $\Delta(x)=g x$ for some fixed element $g\in\Gamma$.
  \end{itemize}
  Consider the action $\Gamma\actson X\times Y$ given by
  $g(x,y)=(g x,\pi(g)y)$. Then the outer automorphism group of the type
  II$_\infty$ relation $\RelR=\Rel(\Gamma\actson X\times Y)$ is given
  by
  \[\Out(\RelR)=\Centr_{\Autns(Y,\nu)}(\Lambda).\]
  In particular, for every subset $U\subset X\times Y$ with finite
  measure, we know that
  \[\Fundg(\RelR\restrict{U})=\Mod(\Centr_{\Autns(Y,\nu)}(\Lambda))
  \text{ and }
  \Out(\RelR\restrict{U})=\Centr_{\Autmp(Y,\nu)}(\Lambda).\]
  Here we denoted the group of all measure preserving automorphisms of
  $Y$ by $\Autmp(Y,\nu)$.
\end{theorem}

In fact, there is a relatively easy example of an action
$\Gamma\actson X$ that satisfies the conditions of theorem
\ref{thm:intro:main:eqrel}. We give a short description here, and more
details can be found in proposition \ref{prop:ex:eqrel}. Consider the action of
$\Gamma=\GL_n\IQ\ltimes\IQ^n$ on $I=\IQ^n$ by affine
transformations. Define $\Gamma\actson (X,\mu)=(X_0,\mu_0)^I$ to be the
generalized Bernoulli action over a purely atomic base space $(X_0,\mu_0)$ with
unequal weights. Let $\pi:\Gamma\rightarrow\Lambda=\IQ^\times$ be the
determinant morphism. Then the action $\Gamma\actson (X,\mu)$ and the
quotient $\pi:\Gamma\rightarrow\Lambda$ satisfy the conditions of
theorem \ref{thm:intro:main:eqrel}.

In particular, for every infinite measure preserving ergodic action
$\IQ^\times\actson (Y,\nu)$, there is a type II$_1$ equivalence
relation $\RelR$ with
\[\Fundg(\RelR)=\Mod(\Centr_{\Autns(Y,\nu)}(\IQ^\times))
\text{ and }
\Out(\RelR)=\Centr_{\Autmp(Y,\nu)}(\IQ^\times).\]
Observe that we do not need the action $\IQ^\times\actson Y$ to be
free. Remark also that $\IQ^\times$ is isomorphic to $\{\pm1\}\times
\IZ^{(\IN)}$. So for every countable abelian group $\Lambda_0$ there is
a quotient morphism $\IQ^\times\rightarrow\Lambda_0$. So we can
compose any action of $\Lambda_0$ with this quotient and obtain an
action of $\IQ^\times$.

In \cite{Aaronson:IntrinsicNormalizingConstants,AaronsonNadkarni:EigenvaluesAndSpectra}, Aaronson and Nadkarni show that
$\Mod(\Centr_{\Autns(Y,\nu)}(\IZ))\subset \IRpos$  can have any
Hausdorff dimension between $0$ and $1$. This gives an explicit
construction of type II$_1$ equivalence relations whose fundamental
group can have any Hausdorff dimension. The existence of such
equivalence relations was already shown in \cite{PopaVaes:ActionsOfFinfty,PopaVaes:FundamentalGroupGeneral}.

Let $\cG\subset\Autmp(Y,\nu)$ be a closed abelian subgroup of measure
preserving automorphisms on $(Y,\nu)$. Assume that $\cG$ acts
ergodically on $Y$. Because $\Autmp(Y,\nu)$ is a
Polish group, we find a group morphism $\sigma:\IQ^\times\rightarrow
\cG$ with dense range. This defines an ergodic measure preserving
action of $\IQ^\times$ on $Y$. Observe that the centralizer of this action is
exactly $\Centr_{\Autns(Y,\nu)}(\IQ)=\Centr_{\Autns(Y,\nu)}(\cG)$.
From this it follows immediately that every second countable locally
compact abelian group appears as the outer automorphism group of some
type II$_1$ equivalence relation.

In the example above, the quotient group $\Lambda=\IQ^\times$ is
abelian. In section \ref{sect:ex:general}, we give a more complex example
of an action $\Gamma\actson (X,\mu)$. But this time the corresponding
quotient group $\Lambda$ is the free group $\FG_\infty$. This has the
following interesting consequence.
\begin{theorem}[see also theorems \ref{thm:fundg:cartan} and \ref{thm:ex:out:eqrel}]
  \label{thm:intro:ex:eqrel}
  Let $\cG\subset \Autmp(Y,\nu)$ be a closed group of measure
  preserving automorphisms of an infinite measure space, acting
  ergodically. Then there is a type II$_1$ equivalence relation
  $\RelR$ with
  \[\Fundg(\RelR)=\Mod(\Centr_{\Autns(Y,\nu)}(\cG))
  \text{ and }
  \Out(\RelR)=\Centr_{\Autmp(Y,\nu)}(\cG).\]
\end{theorem}

It follows that every second countable locally compact unimodular
group appears as the outer automorphism group of some type II$_1$
equivalence relation.

The situation for type II$_1$ factors is very similar, but it requires
an extra step. Remember that we want to control the outer automorphism
group of $M=\Lp^\infty(X\times Y)\rtimes\Gamma$. In general, let
$\Gamma$ act freely on an infinite measure space $(Y,\nu)$.
Denote $M=\Lp^\infty(Y)\rtimes\Gamma$ and $\RelR=\Rel(\Gamma\actson
Y)$. We have that
\begin{equation}
  \label{eqn:incl:obvious}
  \Out(\RelR)\ltimes\Cohom(\Gamma\actson Y)\subset \Out(M).
\end{equation}
This inclusion can be strict, see
\cite{ConnesJones:NonconjugateCartan,Popa:StrongRigidity2}. But if we know
that every automorphism of $M$ preserves the Cartan subalgebra
$\Lp^\infty(Y)$ up to a unitary, then the inequality in
(\ref{eqn:incl:obvious}) is actually an equality. It is a hard problem
to decide if every automorphism of $M$ preserves $\Lp^\infty(Y)$ but
by now we know many classes of actions $\Gamma\actson Y$ with this
property \cite{IoanaPetersonPopa:AFP,Popa:Betti,Popa:StrongRigidity1,Popa:StrongRigidity1,PopaVaes:StrongRigidity,PopaVaes:vNsuperrigidity}.

Combining theorem \ref{thm:intro:main:eqrel} with the results from
\cite{IoanaPetersonPopa:AFP} yields the following result.
\begin{theorem}[see theorem \ref{thm:proof:main:factor}]
  \label{thm:intro:main:factor}
  Let $\Gamma=\Gamma_1\free_\Sigma\Gamma_2$ be an amalgamated free
  product with $\Gamma_1\not=\Sigma\not=\Gamma_2$. Let
  $\pi:\Gamma\rightarrow\Lambda$ be a quotient morphism. Assume that
  \begin{itemize}
  \item $\Gamma_1\cap\ker\pi$ contains an infinite property (T) group
    $G$.
  \item the group $\Sigma$ is amenable.
  \item there are group elements $g_1,\ldots,g_n\in\Gamma$ such that
    the intersection $\bigcap_{i=1}^n g_i\Sigma g_i^{-1}$ is finite.
  \end{itemize}

  Let $\Gamma\actson (X,\mu)$ be a free, ergodic, p.m.p.\ action that
  satisfies the conditions of theorem \ref{thm:intro:main:eqrel}. For
  an arbitrary ergodic infinite measure preserving action
  $\Lambda\actson (Y,\nu)$, we define an action $\Gamma\actson X\times
  Y$ by the formula $g(x,y)=(g x,\pi(g)y)$.

  Denote $M=\Lp^\infty(X\times Y)\rtimes\Gamma$ and
  $\RelR=\Rel(\Gamma\actson X\times Y)$. Then we have that
  \[\Out(M)=\Out(\RelR)\ltimes\Cohom(\Gamma\actson X\times Y)
  =\Centr_{\Autns(Y,\nu)}(\Lambda)\ltimes\Cohom(\Lambda\actson Y).\]
  In particular, for every finite projection $p\in M$ we find that
  \[\Fundg(p M p)=\Mod(\Centr_{\Autns(Y,\nu)}(\Lambda))
  \text{ and }
  \Out(M)=\Centr_{\Autmp(Y,\nu)}(\Lambda)\ltimes
  \Cohom(\Lambda\actson Y).\]
\end{theorem}

Also in this case, we can give an explicit example of an action
$\Gamma\actson (Y,\nu)$ that satisfies the conditions of theorem
\ref{thm:intro:main:factor}. Just as in the equivalence relation case,
we can choose the quotient group $\Lambda$ to be $\FG_\infty$. For the
fundamental group, we obtain the following.
\begin{theorem}
  \label{thm:intro:fundg:cartan}
  Let $\cG\subset\Autmp(Y,\nu)$ be a closed ergodic subgroup of
  (possibly infinite) measure preserving automorphisms.
  Then there is a type II$_1$ factor $M_{\cG}$ with
  \[\Fundg(M_{\cG})=\Mod(\Centr_{\Autns(Y,\nu)}(\cG))\]
\end{theorem}

For the outer automorphism group, things are a little more subtle. The
cohomology group of an $\FG_\infty$-action is usually a non-locally
compact abelian group. Since we are most interested in locally compact
outer automorphism groups, we want to make this cohomology group
vanish, using a group different from $\FG_\infty$. We succeeded in
some specific cases.
\begin{theorem}[see theorems \ref{thm:ex:out:factor} and \ref{thm:ex:out:SLnR}]\ 
  \label{thm:intro:out:factor}
  \begin{itemize}
  \item If $\cG\subset \Autmp(Y,\nu)$ is a closed (not necessarily
    free or ergodic) subgroup of \emph{probability} measure preserving
    transformations, then we construct a type II$_1$ factor $M_{\cG}$ with outer
    automorphism group
    \[\Out(M_{\cG})=\Centr_{\Autmp(Y,\nu)}(\cG).\]
  \item For every $n\in \IN$, there is a type II$_1$ factor whose
    outer automorphism group is $\SL^{\pm}_n\IR$.
  \end{itemize}
\end{theorem}

The first of these two statements gives an explicit construction for
type II$_1$ factors with prescribed compact second countable outer
automorphism group. But many more groups can be realized this way. For
example $\Autmp(Y,\nu)$ for the non-atomic probability space $(Y,\nu)$
also appears as the outer automorphism group of a type II$_1$
factor.

In fact, we can further generalize theorem
\ref{thm:intro:main:factor}, using some kind of first quantisation
step. We replace the infinite measure space $(Y,\nu)$ by an arbitrary
semifinite von Neumann algebra $(B,\Tr)$ with a fixed semifinite
trace. The automorphism group of $(Y,\nu)$ is replaced by the
\emph{outer} automorphism group of $B$. As before, the most general
statement contains a 1-cohomology group of the action of $\Lambda$ on
$B$. This group adds more technical problems than in the case of a
measure space. I do not want to delve into these
technicalities here. Therefore we only state the result for the
fundamental group. More details are given in section \ref{sect:noncartan}.
We obtain the following result.
\begin{theorem}[see also theorem \ref{thm:fundg:noncartan}]
\label{thm:intro:fundg:noncartan}
Let $(B,\Tr)$ be any properly infinite but semifinite von Neumann
algebra with a specified semifinite trace. Let $\Lambda\subset
\Outtp(B,\Tr)$ be a group of $\Tr$-preserving outer automorphisms of
$B$. Assume that this group acts ergodically on the center of $B$. Then
there is a type II$_1$ factor $M_\Lambda$ with fundamental group
\[\Fundg(M_\Lambda)=\Mod(\Centr_{\Out(B)}(\Lambda)).\]
\end{theorem}

Unlike the other results in this paper, this is a pure existence
result. We use Ozawa's result \cite{Ozawa:NoUniversalII1Factor} that
for every $B$, there exists a property (T) group $G$ such that $G$
does not embed into the unitary group of $p B p$ for any finite
projection $p\in B$. In many interesting cases however, we know an explicit
example of such a group $G$. For example, with
$B=\Bounded(\ell^2(\IN))\otimes R$ or
$B=\Bounded(\ell^2(\IN))\otimes\Lg(\FG_\infty)$, we can take any property (T)
group $G$.

Every fundamental group $\cF(M)$ of a type II$_1$ factor appears in
the way of theorem \ref{thm:intro:fundg:noncartan}:
\[\Fundg(M)=\Mod(\Centr_{\Out(\Bounded(\ell^2(\IN)))\otimes M}(\{\id\})).\]
This gives an alternative characterisation of the set of all
fundamental groups of type II$_1$ factors. But this characterization
does not solve Murray and von Neumann's question, because it is not an
intrinsic description. We conjecture that every fundamental group
of a type II$_1$ factor appears as
$\Mod(\Centr_{\Out(\Bounded(\ell^2(\IN)))\otimes\Lg(\FG_\infty)}(\Lambda))$
for some subgroup $\Lambda$. If this conjecture would be true, then
this would provide an important step in the solution of Murray and von
Neumann's question.
\let\thetheorem=\oldthetheorem

\section{Preliminaries and notations}
\subsection{The fundamental group and the outer automorphism group of
  II\texorpdfstring{$_1$}{1} equivalence relations and II\texorpdfstring{$_1$}{1} factors}
Murray and von Neumann defined \cite{MvN:ROO4} the fundamental group
of a II$_1$ factor as
\[\Fundg(M)=\{\tau(p)/\tau(q)\mid p M p\cong q M q\},\]
where $p$ and $q$ run over the nonzero projections in $M$. Even though
this definition makes sense for all II$_1$ factors, we consider only
\emph{separable} von~Neumann algebras (i.e. von~Neumann algebras
acting on a separable Hilbert space). There is a similar notion of fundamental
group $\Fundg(\RelR)$ for a II$_1$ equivalence relation $\RelR$ . A II$_1$
equivalence relation is an ergodic, countable, probability measure preserving
equivalence relation on a standard
probability space $(X,\mu)$. The fundamental group of $\RelR$ is defined by
\[\Fundg(\RelR)=\{\mu(U)/\mu(V)\mid \RelR\restrict U\cong \RelR\restrict
V\}\]
where $\RelR\restrict U$ denotes the restriction of $\RelR$ to $U$,
i.e. for two points $x,y\in U$ we have that $x\RelR\restrict U y$ if
and only if $x\RelR y$.
Associated to a II$_1$ equivalence relation $\RelR$ is
a separable II$_1$ factor $\Lg(\RelR)$ and a Cartan subalgebra
$\Lp^\infty(X,\mu)\subset \Lg(\RelR)$ via the generalized
group-measure space construction of Feldman-Moore
\cite{FeldmanMoore:ErgodicEquivalenceRelations}. The
fundamental group of $\RelR$ is always a subgroup of
$\Fundg(\Lg(\RelR))$. This inclusion can be strict: for example in
\cite{Popa:SpectralGap}, based on results in
\cite{ConnesJones:NonconjugateCartan}, a II$_1$ equivalence relation
$\RelR$ is constructed such that $\Fundg(\RelR)=\{1\}$ while
$\Fundg(\Lg(\RelR))=\IR_+$. On the other hand, if for every $p,q\in A$
and every isomorphism $\alpha:p M p\rightarrow q M q$, the Cartan
subalgebra $p A$ is
mapped onto $uqAu^\ast$ for some unitary $u\in q M q$, then we have
$\Fundg(\Lg(\RelR))=\Fundg(\RelR)$.

To a non-singular action $\Gamma\actson
(X,\mu)$, one associates the orbit equivalence relation
$\Rel(\Gamma\actson X)$ given by $x\sim y$ iff.\ $x=g y$ for some
$g\in \Gamma$. If the action is ergodic and p.m.p.\ (probability measure
preserving) on a standard probability space, then
$\Rel(\Gamma\actson X)$ is a II$_1$ equivalence relation. If the
action $\Gamma\actson X$ is free, then we have that
$\Lg(\Rel(\Gamma\actson X))=\Lp^\infty(X)\rtimes\Gamma$.

We will encounter ergodic actions $\Gamma\actson (X,\mu)$ preserving
the infinite non-atomic standard measure $\mu$. The associated orbit equivalence relation
$\Rel(\Gamma\actson X)$ is a so-called II$_\infty$
equivalence relation. For any II$_\infty$ equivalence relation $\RelR$, the
restriction $\RelR\restrict U$ to a subset with $0<\mu(U)<\infty$ is a
II$_1$ equivalence relation. From now on, we only consider II$_1$ and
II$_\infty$ equivalence relations.

To an equivalence relation $\RelR$ on $(X,\mu)$, several groups are associated. The
full group $\fullg(\RelR)$ is the group of non-singular automorphisms
$\Delta$ of $(X,\mu)$ with $\Delta(x)\sim x$ for almost all $x\in
X$. Remark that all such $\Delta$ are measure preserving. An
automorphism of $\RelR$ is a non-singular automorphism of $(X,\mu)$
such
that $\Delta(x)\sim \Delta(y)$ iff $x\sim y$. We identify
automorphisms that are equal almost everywhere. The group of all
automorphisms of $\RelR$ will be denoted by $\Autr(\RelR)$.
Note that such an automorphism $\Delta$ preserves the measure in the
case of a II$_1$ relation, and scales the measure by a positive
constant $\Mod(\Delta)$ if $\RelR$ is a II$_\infty$ relation. The full
group $\fullg(\RelR)$ is a normal subgroup of $\Autr(\RelR)$. The quotient
$\Out(\RelR)=\Autr(\RelR)/\fullg(\RelR)$ is called the outer automorphism group.
If $\RelR$ is a II$_\infty$ relation and
$0<\mu(U)<\infty$, then we have the following short exact sequence
\[\begin{CD}
  1 @>>> \Out(\RelR\restrict U) @>>> \Out(\RelR) @>\Mod>>
  \Fundg(\RelR\restrict U) @>>> 1
\end{CD}\]
In particular, if we can identify $\Out(\RelR)$, we also know
$\Fundg(\RelR\restrict U)$ and $\Out(\RelR\restrict U)$.

The situation is similar for factors. Let $M$ be a factor. We denote
by $\Autf(M)$ the group of all automorphisms of $M$. If $M$ is a type II$_1$ factor,
then every automorphism of $M$ preserves the trace, and if $M$ is of
type II$_\infty$, then every automorphism scales the trace by a
positive constant $\Mod(\psi)\in \IRpos$. An automorphism
$\psi:M\rightarrow M$ is called inner if it is given by unitary
conjugation, i.e.\ there is a unitary $u\in M$ such that
$\psi(x)=uxu^\ast$ for all $x\in M$. The group of all inner
automorphisms is denoted by $\Inn(M)$.
This group is a normal subgroup of $\Autf(M)$, and we call the quotient
$\Out(M)=\Autf(M)/\Inn(M)$ the outer automorphism group of $M$.

If $M$ is a type II$_\infty$ factor and $p\in M$ is a projection with
finite trace, then we have the short exact sequence
\[\begin{CD}
  1 @>>> \Out(p M p) @>>> \Out(M) @>\Mod>>
  \Fundg(p M p) @>>> 1
\end{CD}\]

\subsection{Cocycle superrigidity and \texorpdfstring{$\Ufin$}{Ufin} groups}
\label{subsect:prelim:CSR}
Let $\Gamma\actson (X,\mu)$ be a probability measure
preserving action. We say that $\Gamma\actson (X,\mu)$ is
weakly mixing if the diagonal action of $\Gamma$ on $X\times X$ is
ergodic. Many equivalent definitions are known, see for example
\cite[Appendix D]{Vaes:Bourbaki}. For us, the most important
equivalent definition is the following.
\begin{lemma}[{see \cite[lemma 5.4]{PopaVaes:lattices}}]
  Let $\Gamma\actson (X,\mu)$ be a weakly mixing action.
  
  Let $(Z,d)$ be a Polish space with separable complete metric
  $d$. Suppose that a Polish group $\Lambda$ has a continuous action
  $\alpha$ by isometries on $Z$.
  
  Let $\Gamma\actson (Y,\nu)$ be another measure preserving
  action on a standard measure space and let
  $\omega:\Gamma\times Y\rightarrow \Lambda$ be a measurable
  map. If a measurable map $f:X\times Y\rightarrow Z$
  satisfies
  \[f(g x,g y)=\alpha_{\omega(g,y)}f(x,y),\]
  then there is a measurable map $f_0:Y\rightarrow Z$ such that
  $f(x,y)=f_0(y)$ almost everywhere.
\end{lemma}

For example, all (plain) Bernoulli actions of infinite groups are
weakly mixing. More generally, a generalized Bernoulli action
$\Gamma\actson (X_0,\mu_0)^I$ is
weakly mixing iff. it is ergodic iff. every orbit of $\Gamma\actson I$
is infinite.

Fix a Polish group 
$\Lambda$. A measurable map $\omega:\Gamma\times X\rightarrow\Lambda$ is called
a 1-cocycle if $\omega$ satisfies the cocycle relation
\[\omega(g,hx)\omega(h,x)=\omega(gh,x)\text{ for almost every }x\in
X\text{ and }g,h\in\Gamma.\]
Examples are given by Zimmer cocycles: Let $\Lambda\actson
(Y,\nu)$ be a free non-singular action and let $\Delta:X\rightarrow Y$ be an orbit
equivalence (i.e. $\Delta$ is an isomorphism between the orbit equivalence
relations $\Rel(\Gamma\actson X)$ and $\Rel(\Lambda\actson Y)$, or still
$\Delta$ is a non-singular isomorphism such that
$\Delta(\Gamma x)=\Lambda\Delta(x)$ for almost all $x\in X$). The
Zimmer cocycle is the unique (up to measure 0) map $\omega:\Gamma\times X\rightarrow
\Lambda$ such that $\Delta(gx)=\omega(g,x)\Delta(x)$ almost everywhere.
If $\Delta_1,\Delta_2$ are two orbit equivalences and there is an
$\eta\in[\Rel(\Lambda\actson Y)]$ with
$\Delta_1=\eta\circ\Delta_2$, then $\omega_1$ and $\omega_2$ are
homologous. This means that there is a function
$\varphi:X\rightarrow \Lambda$ such that
\[\omega_1(g,x)=\varphi(gx)^{-1}\omega_2(g,x)\varphi(x)\text{ for
  all }g\text{ and almost all }x\]

The group morphisms $\delta:\Gamma\rightarrow \Lambda$ can be
identified with the cocycles that do not depend on the
$X$-variable. Popa showed the following cocycle superrigidity
theorem. To state it, we need the notion of a $\Ufin$ group.
A Polish group $\Lambda$ is called a $\Ufin$ group if it can be realized as a
closed subgroup of the unitary group of a finite von~Neumann
algebra. Examples include all countable groups and all compact second
countable groups.
\begin{theorem}[{\cite[theorem
    0.1]{Popa:CocycleSuperrigidityMalleable}}] Suppose that 
  $\Gamma\actson (X,\mu)$ is
  an s-malleable action (for example a generalized Bernoulli action)
  of a group $\Gamma$ with an infinite almost normal subgroup $G$
  that has the relative property (T). Assume that $G$ acts weakly
  mixingly on $(X,\mu)$. Then every cocycle $\omega:\Gamma\times
  X\rightarrow \Lambda$, with a $\Ufin$ target group $\Lambda$, is
  cohomologous to a group morphism.
\end{theorem}
Groups $\Gamma$ with infinite almost normal subgroups $G$ that have the
relative property (T) are called weakly rigid. A subgroup
$G\subset\Gamma$ is called almost normal if $G\cap gGg^{-1}$ has
finite index in $G$, for all $g\in\Gamma$.

We say that an action $\Gamma\actson(X,\mu)$ is
$\Ufin$-cocycle superrigid if every cocycle $\omega:\Gamma\times
X\rightarrow 
\Lambda$, to a $\Ufin$ target group $\Lambda$, is cohomologous to a
group morphism.
For example, \cite[theorem 0.1]{Popa:CocycleSuperrigidityMalleable}
asserts that s-malleable actions of weakly rigid groups are
$\Ufin$-cocycle superrigid.

\subsection{Intertwining by bimodules in the semifinite case}
Given two subalgebras $A$ and $B$ in a finite von Neumann algebra $M$, we want to
decide if
$A$ and $B$ are unitarily conjugate or not. In \cite[section 2]{Popa:StrongRigidity1} and
\cite[theorem A.1]{Popa:Betti}, Popa introduced a
powerful technique to decide this type of question. In fact, we will need a
similar technique that works if
$M$ is a semifinite von Neumann algebra instead of a finite von
Neumann algebra. Such a generalization was developed
in \cite{ChifanHoudayer:BassSerreRigidity} and later an improved version was given in
\cite{HoudayerRicard:FreeArakiWoodsFactors}. We give an overview 
of the terminology and of the results that we will use.

Let $(B,\Tr_B)$ be a von Neumann algebra with a faithful normal
semifinite trace. From now on, a trace will always be assumed to be faithful and
normal. Any right 
Hilbert $B$-module $H$ is isomorphic to a submodule of $\ell^2(\IN)\otimes \Lp^2(B,\Tr_B)$,
so there is a projection $p\in\cB(\ell^2(\IN))\vnOtimes B$ such that
$H\cong p(\ell^2(\IN)\otimes \Lp^2(B,\Tr_B))$. The trace $\Tr\otimes
\Tr_B(p)$ of $p$ is an invariant 
of the module $H$, which we call the dimension
$\dim_{\Tr_B}(H_{B})$. Observe that the dimension 
of $H$ depends on the choice of trace $\Tr$, and even in the case
of a II$_\infty$ \emph{factor} $B$, the canonical trace $\Tr$ is only
determined up to a positive scalar multiple.

\begin{definition_theorem}[{\cite[lemma
    2.2]{HoudayerRicard:FreeArakiWoodsFactors},
    generalizing \cite[theorem 2.1]{Popa:StrongRigidity1}}]
  Let $(M,\Tr)$ be a von Neumann algebra with a
  semifinite trace. Let $p\in M$
  be a projection with finite trace and let $A\subset p M p$ be a von
  Neumann subalgebra. Suppose
  that $B\subset M$ is a unital von Neumann subalgebra such that the restricted trace
  $\Tr\restrict B$ is still semifinite. We say that $A$ embeds
  into $B$ over $M$ (denoted by $A\embeds_M B$) if one of the following
  equivalent conditions holds.
  \begin{itemize}
  \item There is an $A$-$B$ subbimodule $H\subset p\Lp^2(M)$ that has finite dimension over
    $(B,\Tr)$.
  \item There exist a natural number $n$, a projection $q\in \MatM_n(\IC)\otimes B$ with
    finite trace, a partial isometry $v\in p(\MatM_{1,n}(\IC)\otimes M)q$ and
    a unital $\ast$-homomorphism 
    $\theta:A\rightarrow q(\MatM_n(\IC)\otimes B)q$ such that
    $av=v\theta(a)$ for all $a\in A$. 
  \item There is \emph{no} sequence $(v_n)_n$ of unitaries in $A$ such that
    \[\norm{\E_B(xv_ny)}_{\Tr,2}\rightarrow 0\text{ for all }x,y\in M.\]
  \end{itemize}
  In terms of Popa's notion for finite von Neumann algebras, this is equivalent with
  \begin{itemize}
  \item There is a projection $q\in B$ with finite trace such that
    $A\embeds_{(p\vee q)M(p\vee q)} qBq$. 
  \end{itemize}
\end{definition_theorem}

Let $C\subset M$ be a unital subalgebra.
If there is a unitary $u\in M$ with $uCu^\ast\subset B$, then we have
$pCp\embeds_M B$ for every projection $p\in A$ with finite trace. The
converse implication does not hold in general, but if both $C$ and
$B$ are Cartan subalgebras, then the fact that $pC\embeds_M B$, for
some $p\in C$ with finite trace, is sufficient to conclude that $C$
is unitarily conjugate to $B$. In the case where $M$ is a II$_1$
factor, this fact was proven in \cite[theorem A.1]{Popa:Betti}. For a
proof of the generalization  to
the semifinite case, see \cite[proposition 2.3]{HoudayerRicard:FreeArakiWoodsFactors}. 
In proposition \ref{prop:inter:regular}, we give a different set of
conditions on two subalgebras $B,C\subset M$ such that $p B p\embeds_M
C$ and $qCq\embeds_M B$ implies that $B$ and $C$ are unitarily
conjugate in $M$.

Suppose that $M=B\rtimes\Gamma$ where $(B,\Tr)$ is a von Neumann
algebra with a semifinite trace and $\Gamma$ acts on $B$ preserving
$\Tr$. Let $A\subset p M p$ be a regular subalgebra of a finite corner
of $M$. Assume that we have $A\embeds_M B\rtimes\Sigma$ for some
(highly non-normal) subgroup
$\Sigma\subset \Gamma$. In \cite[theorem
6.16]{PopaVaes:StrongRigidity}, a criterion is given to show 
that then we also have $A\embeds_M B$. An improved version of this
criterion is given in \cite{HoudayerPopaVaes:AppvNsuperrigidity}: For any
subgroup $\Sigma\subset \Gamma$, we denote by $\zH_{\Sigma}$ the closed
linear span of all $A$-$B\rtimes\Sigma$ subbimodules of $p\Lp^2(M)$
that have finite dimension over $B\rtimes\Sigma$. Because $A$ is
regular, this space $\zH_{\Sigma}$ is invariant under left
multiplication by $p M p$, so $\zH_{\Sigma}$ is of the form
$\zH_{\Sigma}=p\Lp^2(M)\z_{\Sigma}$ for some projection $\z_{\Sigma}$ in the
relative commutant of $B\rtimes\Sigma$ inside $M$. This projection
$\z_{\Sigma}$ is uniquely determined if we also require that
$\z_{\Sigma}$ is smaller than the central support of $p$. By
definition, we have $\z_{\Sigma}\not=0$ if and only if $A\embeds_M
B\rtimes \Sigma$.

Given $g\in\Gamma$, we have $\z_{g\Sigma g^{-1}}=u_g\z_{\Sigma}u_g^\ast$.
For two subgroups $\Sigma_1,\Sigma_2\subset\Gamma$,
\cite{HoudayerPopaVaes:AppvNsuperrigidity}
shows that
$\z_{\Sigma_1\cap\Sigma_2}=\z_{\Sigma_1}\z_{\Sigma_2}$.
Remember that $\z_{\Sigma}$ is always contained in $M\cap
(B\rtimes\Sigma)^\prime$. Suppose that $u_g$ commutes with all of
$M\cap (B\rtimes\Sigma)^\prime$ for some $g\in\Gamma$. Then the
combination of both observations in the beginning of this paragraph
shows that we have $\z_{\Sigma\cap g\Sigma g^{-1}}=\z_{\Sigma}$.
By induction, they obtain the following.
\begin{corollary}[{see \cite[Corollary 7]{HoudayerPopaVaes:AppvNsuperrigidity}}]
  \label{cor:regularInNonnormal}
  Let $M=B\rtimes\Gamma$ be a the crossed product of an action of
  $\Gamma$ preserving a semifinite trace $\Tr$ on $B$.
  Suppose that a regular subalgebra $A\subset p M p$ of a finite
  corner of $M$ satisfies $A\embeds_M B\rtimes\Sigma$, for a subgroup
  $\Sigma\subset \Gamma$. If there are $g_1,\ldots,g_n\in\Gamma$
  such that
  \begin{itemize}
  \item the $u_{g_i}$ are contained in the relative bicommutant of
    $B\rtimes\Sigma$, i.e.\[u_{g_i}\in M\cap \left(M\cap
      (B\rtimes\Sigma)^\prime\right)^\prime,\]
  \item the intersection $\bigcap_{i=1}^ng_i\Sigma g_i^{-1}$ is finite,
  \end{itemize}
  then we also have $A\embeds_M B$.
\end{corollary}

Actually, \cite[Corollary 7]{HoudayerPopaVaes:AppvNsuperrigidity} is stated
with a finite trace $\Tr$, but its proof is still valid in the
semifinite case. A second result in
\cite{HoudayerPopaVaes:AppvNsuperrigidity} states that if $M$ is a
II$_1$ factor and if $\Gamma=\Gamma_1\free_\Sigma\Gamma_2$ is an
amalgamated free product with $\Gamma_1\not=\Sigma\not=\Gamma_2$,
then $\z_{\Sigma}$ is either $1$ or $0$. The proof uses the finiteness
of the trace. However, the following variant is still true. Although
the proof is almost identical to the proof of \cite[Proposition
8]{HoudayerPopaVaes:AppvNsuperrigidity}, we provide a proof for the
convenience of the reader.
\begin{proposition}[{see \cite[Proposition
    8]{HoudayerPopaVaes:AppvNsuperrigidity}}]
  \label{prop:inter:free}
Let $\Gamma=\Gamma_1\free_\Sigma\Gamma_2$ be an amalgamated free
product with $\Gamma_1\not=\Sigma\not=\Gamma_2$. Let
$\pi:\Gamma\rightarrow\Lambda$ be a quotient morphism such that
$\ker\pi$ is not contained in $\Sigma$. Let $\Gamma\actson (X,\mu)$ be
a \emph{free, p.m.p.}\ action, and let $\Lambda\actson (B,\Tr)$ be any
trace preserving action on a semifinite von Neumann algebra. Consider
the action of $\Gamma$ on $\Lp^\infty(X)\otimes B$ given by
$g(a\otimes b)=ga\otimes\pi(g)b$, and denote the crossed product by
$M=(\Lp^\infty(X)\otimes B)\rtimes\Gamma$.

Let $A\subset p M p$ be a regular subalgebra of a finite corner of
$M$. Then the projections $\z_\Sigma$, $\z_{\Gamma_1}$ and
$\z_{\Gamma_2}$, as defined above, coincide.
\end{proposition}
In particular, this
projection $\z_{\Sigma}=\z_{\Gamma_1}=\z_{\Gamma_2}$ is contained in
the center of $M$. So whenever $A\embeds_M (\Lp^\infty(X)\otimes
B)\rtimes\Gamma_i$ for $i=1$ or $i=2$, then it follows that
\[A\embeds_M (\Lp^\infty(X)\otimes B)
\rtimes (g_1\Sigma g_1^{-1}\cap\ldots\cap g_n\Sigma g_n^{-1})\]
for any finite number of elements $g_1,\ldots,g_n\in\Gamma$.
\begin{proof}
  We show that $\z_\Sigma=\z_{\Gamma_1}$. By symmetry it then follows
  that $\z_{\Sigma}=\z_{\Gamma_2}$. Observe that $\z_\Sigma$ and
  $\z_{\Gamma_1}$ are contained in the abelian von Neumann algebra
  $\Lp^\infty(X)\otimes \Centre{B}$. Moreover, $\z_\Sigma$ is smaller
  than $\z_{\Gamma_2}$.

  Denote by $S$ the set of all
  elements $g\in\Gamma$ that have a reduced expression whose first
  letter comes from $\Gamma_2\setminusG\Sigma$. For all $g\in S$, we
  see that $g\Gamma_1g^{-1}\cap\Gamma_1\subset \Sigma$, so it follows
  that $\z_{\Gamma_1}u_g\z_{\Gamma_1}u_g^{\ast}\leq \z_\Sigma$. We claim
  that
  \[q=\bigvee_{g\in S} \z_{\Gamma_1}
  u_g\z_{\Gamma_1}u_g^{\ast}=\z_{\Gamma_1}.\]
  Suppose this was not the case, then we denote
  $r=\z_{\Gamma_1}-q$. Whenever $g\in S$, it is clear that
  \[r\wedge u_g r u_g^\ast\leq u_g\z_{\Gamma_1}u_g^\ast\wedge
  z_{\Gamma_1}\wedge r=0\].

  On the other hand, there are nonzero projections in the abelian von
  Neumann algebras $r_1\in\Lp^\infty(X)$ and $r_2\in\Centre(B)$ such
  that $r_1\otimes r_2 \leq r$. In particular, $r_1$ is orthogonal to
  $u_g r_1 u_g^\ast$ for all $g\in\ker\pi\cap S$. If we can show that
  $\ker\pi\cap S$ is infinite, we would have shown that the
  probability space $X$ contains infinitely many mutually disjoint
  subsets of equal, nonzero measure. This contradiction proves that
  $r=\z_{\Gamma_1}$ and hence that $\z_{\Sigma}=\z_{\Gamma_1}$,
  finishing the proof of proposition \ref{prop:inter:free}.

  We still have to show that $\ker\pi\cap S$ is infinite. In fact, we
  show that there is an element $g\in\ker\pi$ of infinite order and
  such that $g^n\in S$ for all $n\in\IN$. Let
  $g\in\ker\pi\setminus\Sigma$ be any element. We consider four
  cases. If the first letter of
  $g$ comes from $\Gamma_2\setminusG\Sigma$ and its last letter comes from
  $\Gamma_1\setminusG\Sigma$, then $g^n$ is clearly contained in
  $S$ for all $n\in\IN$. If the first letter comes from
  $\Gamma_1\setminusG\Sigma$ and the last letter comes from
  $\Gamma_2\setminusG\Sigma$, then we have $g^{-n}\in S$. If both the
  first and the last letter of $g$ come from
  $\Gamma_1\setminusG\Sigma$, take any $h\in\Gamma_2\setminus
  \Sigma$. Then $hgh^{-1}g$ is an element of
  $\ker\pi$ whose first letter comes from $\Gamma_2\setminusG\Sigma$
  and whose last letter comes from
  $\Gamma_1\setminusG\Sigma$. Finally, if both the first and the last
  letter of $g$ come from $\Gamma_2\setminus\Sigma$, then we take any
  element $h\in\Gamma_1\setminusG\Sigma$, and we see that
  $(ghgh^{-1})^n\in S$ for all $n\in\IN$.
\end{proof}

\section{The equivalence relation case}
\label{sect:eqrel}
In this section we concentrate on fundamental groups and outer
automorphism groups of type II$_1$ equivalence relations. We mainly
prove theorem \ref{thm:intro:main:eqrel}, but in section
\ref{sect:noncartan}, we will meet cases where the action
$\Lambda\actson (Y,\nu)$ does not preserve any $\sigma$-finite measure
on $Y$. At the end of this section, in proposition
\ref{prop:ex:eqrel}, we give an explicit example of an action
$\Gamma\actson (X,\mu)$ that satisfies the conditions of theorem
\ref{thm:proof:main:eqrel}. This gives explicit examples of type II$_1$
equivalence relations with uncountable fundamental group different
from $\IRpos$.

\begin{theorem}[see also theorem \ref{thm:intro:main:eqrel}]
  \label{thm:proof:main:eqrel}
  Let $\Gamma\actson (X,\mu)$ be a free, ergodic, p.m.p.\ action and
  let $\pi:\Gamma\rightarrow\Lambda$ be a quotient morphism. Suppose
  that
  \begin{itemize}
  \item the action of $\ker\pi$ on $X$ is weakly mixing and
    $\Ufin$-cocycle superrigid.
  \item there are no non-trivial group morphisms
    $\theta:\ker\pi\rightarrow\Lambda$.
  \item $\Norm_{\Autmp(X,\mu)}(\ker\pi)=\Gamma$, i.e.\ all the conjugations
    of the action of $\ker\pi$ on $X$ are described by elements of
    $\Gamma$.
  \end{itemize}
  Let $\Lambda\actson (Y,\nu)$ be any ergodic non-singular action on a
  standard measure space, and consider the action $\Gamma\actson
  X\times Y$ given by $g(x,y)=(gx,\pi(g)y)$. Then the orbit
  equivalence relation $\RelR=\Rel(\Gamma\actson X\times Y)$ has outer
  automorphism group
  \[\Out(\RelR)=\Centr_{\Autns(Y)}(\Lambda).\]

  If the action $\Lambda\actson (Y,\nu)$ preserves the infinite
  measure $\nu$, we recover theorem \ref{thm:intro:main:eqrel}:
  \[\Fundg(\RelR\restrict U)=\Mod(\Centr_{\Autns(Y)}(\Lambda))
  \text{ and }
  \Out(\RelR\restrict U)=\Centr_{\Autmp(Y,\nu)}(\Lambda),\]
  for any subset $U\subset Y$ with finite measure.
\end{theorem}
In the proof of theorem
\ref{thm:proof:main:eqrel}, we encounter local automorphisms of a measure
space $(X,\mu)$. We say that a map $\Delta:X\rightarrow X$ is a local
automorphism if
there is a countable partition $X=\bigsqcup_i X_i$ into Borel
subsets, such that the restriction $\Delta\restrict{X_i}$ is a
non-singular isomorphism between the Borel sets $X_i$ and $\Delta(X_i)$.
\begin{proof}
  The formula $\Delta_0\mapsto\id_X\times\Delta_0$ defines an
  injective group
  morphism from the centralizer $\Centr_{\Autns(Y)}(\Lambda)$ to the outer automorphism
  group $\Out(\RelR)$.

  It remains to show that every orbit equivalence $\Delta:X\times
  Y\rightarrow X\times Y$ is of the form $\id_X\times \Delta_0$, up to
  an inner automorphism. Let $\Delta$ be such an orbit equivalence and
  consider its Zimmer cocycle
  $\omega:\Gamma\times X\times Y\rightarrow \Gamma$. We can consider
  the restriction
  $\omega\restrict{\ker\pi}$ as a cocycle $\tilde\omega:\ker\pi\times
  X\rightarrow \cG$ to the $\Ufin$ target group $\cG$ of measurable
  functions from $Y$ to $\Gamma$. Cocycle superrigidity yields a
  measurable function $\varphi:X\times Y\rightarrow \Gamma$ and group
  morphisms $\delta_y:\ker\pi\rightarrow\Gamma$ such that
  \[\delta_y(g)=\varphi(gx,y)\omega(g,x,y)\varphi(x,y)^{-1}\]
  for all $g\in\ker\pi$ and for almost all $(x,y)\in X\times Y$.
  Because there are no non-trivial group morphisms
  $\theta:\ker\pi\rightarrow\Lambda$, it follows that $\delta_y(\ker\pi)\subset\ker\pi$.
  Then \cite[proposition 3.6]{Popa:CocycleSuperrigidityMalleable} says that
  $\varphi(g(x,y))\omega(g,x,y)\varphi(x,y)$ is essentially
  independent of $x$, for all $g\in\Gamma$. We denote that unique
  essential value by $\delta(g,y)$.

  We define a local automorphism $\widetilde\Delta:X\times
  Y\rightarrow X\times Y$ by the formula
  $\widetilde\Delta(x,y)=\varphi(x,y)\Delta(x,y)$. Observe that
  $\widetilde\Delta(g(x,y))=\delta_y(g)\widetilde \Delta(x,y)$ for all
  $g\in\ker\pi$ and for almost all $(x,y)\in X\times Y$.
  Split $\widetilde\Delta$ into its components:
  $\widetilde\Delta=(\Delta_1(x,y),\Delta_2(x,y))$. Remember that
  $\delta_y(\ker\pi)$ is contained in $\ker\pi$.
  The function $\Delta_2$ satisfies
  \[\Delta_2(gx,y)=\Delta_2(x,y)\text{ for all }g\in \ker\pi\text{ and for
    almost all }(x,y)\in X\times Y.\]
  Since $\ker\pi$ acts ergodically on $X$, we see that $\Delta_2(x,y)$ does
  not depend on $x$. We write just $\Delta_2(y)$.

  \emph{Claim: we show that the formula $\Delta_{1,y}(x)=\Delta_1(x,y)$
    defines and automorphism $\Delta_{1,y}:X\rightarrow X$, for almost
    every $y\in Y$.}
  We can apply the same argument as before to $\Theta=\Delta^{-1}$. So
  we find a local automorphism $\widetilde \Theta:X\times
  Y\rightarrow X\times Y$ and group morphisms
  $\theta_y:\ker\pi\rightarrow\ker\pi$ such that
  $\widetilde\Theta(gx,y)=\theta_y(g)\widetilde\Theta(x,y)$ for all $g\in\ker\pi$ and
  almost all $(x,y)\in X\times Y$. Moreover, we can assume that the
  second component $\Theta_2$ of $\widetilde\Theta$ does not
  depend on the $X$-variable, and we consider it as
  $\Theta_2:Y\rightarrow Y$. We denote the first component of
  $\widetilde\Theta$ by $\Theta_1$, and we write
  $\Theta_{1,y}(x)=\Theta_1(x,y)$. So we have that
  $\widetilde\Theta(x,y)=(\Theta_{1,y}(x),\Theta_2(y))$ almost everywhere.

  Remark that
  $\widetilde\Delta\circ\widetilde\Theta$ and
  $\widetilde\Theta\circ\widetilde\Delta$ are inner local
  automorphisms in the sense that there are measurable maps
  $\varphi,\psi:X\times Y\rightarrow \Gamma$ such that
  \[
  \widetilde\Delta\left(\widetilde\Theta\left(x,y\right)\right)=\varphi(x,y)(x,y)\qquad
  \widetilde\Theta\left(\widetilde\Delta\left(x,y\right)\right)=\psi(x,y)(x,y)\qquad
  \text{almost everywhere}
  \]
  It follows that
  \[\varphi(gx,y)g(x,y)=\widetilde\Delta\left(\widetilde\Theta\left(gx,y\right)\right)=\delta_{\Theta_2(y)}(\theta_y(g))\varphi(x,y)(x,y)\]
  for all $g\in\ker\pi$ and almost all $(x,y)\in X\times Y$. By
  freeness, we see that
  $\varphi(gx,y)=\delta_{\Theta_2(y)}(\theta_y(g))\varphi(x,y)g^{-1}$
  almost everywhere. Weak mixing (see for example \cite[lemma
  5.4]{PopaVaes:lattices}) implies that $\varphi(x,y)$ does not depend
  on $x$, so we write $\varphi(y)=\varphi(x,y)$ a.e. The same argument
  works for $\psi$, and we also write just $\psi(y)$.

  Then it follows that $\Delta_2\circ\Theta_2$ and
  $\Theta_2\circ\Delta_2$ are inner and hence local
  automorphisms. Lemma \ref{lem:fundg:technical1} shows that
  $\Delta_2, \Theta_2$ and the $\Delta_{1,y}, \Theta_{1,y}$ are local automorphisms, for
  almost all $y\in Y$. Now we see that
  $\psi(y)^{-1}\circ\Theta_{1,\Delta_2(y)}$ is a left inverse for
  $\Delta_{1,y}$, for almost every $y\in Y$. Similarly,
  $\Theta_{1,y}\circ\varphi(y)^{-1}$ is a right inverse for
  $\Delta_{1,\Theta_2(y)}$. So for almost all $y$ in the
  non-negligible set $W=\Theta_2(Y)$, we know that $\Delta_{1,y}$ is a
  genuine automorphism of $X$. But for any $g\in\Gamma$, we see that
  $\delta(g,y)\circ\Delta_{1,y}=\Delta_{1,\pi(g)y}\circ g$. The
  ergodicity of $\Lambda\actson Y$ shows that in fact almost every
  $\Delta_{1,y}$ is a genuine automorphism. This finishes the proof of
  our claim.

  Observe that $\Delta_{1,y}$ is a conjugation for the action
  $\ker\pi\actson X$, for
  almost every $y\in Y$. Our third condition shows that
  $\Delta_{1,y}\in\Gamma$, so there is an element $\psi(y)\in\Gamma$
  such that $\Delta_{1,y}(x)=\psi(y)^{-1}x$ almost everywhere. Set
  $\Delta_0(y)=\psi(y)^{-1}\Delta_2(y)$ and observe that
  $\psi(y)^{-1}\widetilde\Delta(x,y)=(x,\Delta_0(y))$ almost
  everywhere.
  Hence, for every $g\in\Gamma$, we see that
  \begin{align*}
    &(gx,\Delta_0(\pi(g)y))\\
    &=\psi(\pi(g)y)^{-1}\varphi(g(x,y))\Delta(g(x,y))\\
    &=\underbrace{\psi(\pi(g)y)^{-1}\varphi(g(x,y))
      \omega(g,x,y)\varphi(x,y)^{-1}\psi(y)}_{\overset{\|}{\widetilde\omega(g,x,y)}}
    \,(x,\Delta_0(y))
  \end{align*}
  Since $\Gamma$ acts freely on $X$, we see that $\widetilde
  \omega(g,x,y)=g$ almost everywhere, and in particular, $\Delta_0$
  commutes with the action of $\Lambda$. We still have to show that
  $\Delta_0$ is an automorphism of $Y$. The same argument as before
  applies to $\Delta^{-1}$, so we find a local automorphism
  $\Theta_0:Y\rightarrow Y$ such that
  $(x,\Theta_0(\Delta_0(y)))\in\Gamma (x,y)$ almost everywhere. By
  freeness we see that $\Theta_0\circ \Delta_0=\id_Y$. Symmetrically,
  we find that $\Delta_0\circ\Theta_0=\id_Y$, so $\Delta_0$ is an
  automorphism of $Y$.

  We have shown that $\Delta$ is of the form $\id_X\times \Delta_0$ up
  to an inner automorphism.
\end{proof}

\begin{lemma}
  \label{lem:fundg:technical1}
  Suppose that $\Delta\!:\!X\!\times\! Y\rightarrow X\!\times\! Y$ and
  $\Theta\!:\!X\!\times\! Y\rightarrow X\!\times\! Y$ are local
  automorphisms. Suppose that the second component of $\Delta$ and
  $\Theta$ does not depend on the $X$-variable, so we can write
  \[\Delta(x,y)=(\Delta_{1,y}(x),\Delta_2(y))\text{ and }
  \Theta(x,y)=(\Theta_{1,y}(x),\Theta_2(y))\]
  almost everywhere.

  If $\Delta_2\circ\Theta_2$ and $\Theta_2\circ\Delta_2$ are local
  automorphisms of $Y$, then $\Delta_2$, $\Theta_2$, the
  $\Delta_{1,y}$ and the $\Theta_{1,y}$ for 
  almost all $y\in Y$ are local automorphisms.
\end{lemma}
\begin{proof}
  \emph{We first show that $\Delta_2$ and $\Theta_2$ are local
    automorphisms.}\\
  Remark that $\Delta^{-1}(X\times W)=X\times\Delta_2^{-1}(W)$ for all
  measurable subsets $W\subset Y$. Since $\Delta$ is non-singular, we
  see that $\Delta_2$ is non-singular. The same argument applies to
  $\Theta_2$.

  We show that $\Delta_2$ and $\Theta_2$ are partial automorphisms.
  We can partition $Y$
  into a countable disjoint union of subsets $Y=\bigsqcup_i Y_i$ such
  that every $\Psi_i=\Theta_2\circ\Delta_2\restrict{Y_i}$ is a non-singular
  automorphism. Similarly, we can partition $Y$ into a disjoint union
  $Y=\bigsqcup_j \widetilde Y_j$ such that every
  $\Phi_j=\Delta_2\circ\Theta_2\restrict{\widetilde Y_j}$
  is a non-singular automorphism. For every $i,j$, it follows that
  $\Delta_2\restrict{Y_i\cap \Delta_2^{-1}(\widetilde Y_j)}$ is a
  non-singular automorphism between $Y_i\cap\Delta_2^{-1}$ and
  $\widetilde
  Y_j\cap\Theta_2^{-1}(\Psi_i(Y_i\cap\Delta_2^{-1}(\widetilde Y_j)))$, so $\Delta_2$ is a local automorphism.
  The same argument applies to $\Theta_2$.

  \emph{Then we show that almost all the $\Delta_{1,y}$ are local
    automorphisms.}\\
  Let $W\subset Y$ be a subset such that $\Delta_2\restrict{W}$ is a
  non-singular isomorphism onto its image.
  Let $U$ be a measurable subset of $X\times W$ such
  that $\Delta\restrict{U}$ is a non-singular isomorphism onto its
  image $V=\Delta(U)$. We denote by $U_y$ the slice
  $U_y=\{x\in X\mid (x,y)\in U\}$.
  Consider the Radon-Nikodym derivative $F\!:\!V\rightarrow \IRpos$
  of $\Delta\restrict{U}$ and denote by $G\!:\!\Delta_2(Y_i)\rightarrow
  \IRpos$ the Radon-Nikodym derivative of
  $\Delta_2\restrict{Y_i}$. Then we compute that, for all measurable
  functions $f:X\rightarrow \IR$ and $g:Y\rightarrow \IR$, we have
  \begin{align*}
    &\int_{Y_i}g(\Delta_2(y))\int_{U_y}f(\Delta_{1,y}(x))\D\mu(x)\D\nu(y)\\
    &=\int_{U}(f\otimes g)(\Delta(x,y))\D(\mu\times\nu)(x,y)\\
    &=\int_{V}(f\otimes g)(x,y)F(x,y)\D(\mu\times\nu)(x,y)\\
    &=\int_{Y_i}g(\Delta_2(y))\int_{V_{\Delta_2(y)}}f(x)\frac{F(x,\Delta_2(y))}{G(\Delta_2(y))}\D\mu(x)\D\nu(y).
  \end{align*}
  This calculation shows that the Radon-Nikodym derivative of
  $\Delta_{1,y}\restrict{U_y}$ is given by the formula
  $H_y(x)=F(x,\Delta_2(y))/G(\Delta_2(y))$ and in particular,
  we see that $\Delta_{1,y}\restrict{U_y}$ is non-singular. Because
  $\Delta\restrict{U}$ is bijective, it follows that
  $\Delta_{1,y}\restrict{U_y}$ is a non-singular isomorphism onto its
  image, for almost all $y\in Y$ with $\mu(U_y)>0$.

  We can partition $X\times Y$ into a countable disjoint union of sets
  like $U$, so $\Delta_{1,y}$ is a local automorphism for almost all
  $y\in Y$.
\end{proof}

\begin{proposition}
  \label{prop:ex:eqrel}
  Consider the action of the affine group
  $\Gamma=\GL_n\IQ\ltimes\IQ^n$ on the countable affine space $I=\IQ^n$,
  with dimension $n\geq 2$. The determinant gives us a quotient morphism
  $\pi:\Gamma\rightarrow \Lambda=\IQ^\times$. Then the generalized
  Bernoulli action $\Gamma\actson (X,\mu)=(X_0,\mu_0)^I$, with atomic
  base space $(X_0,\mu_0)$ with unequal weights, satisfies the
  conditions of theorem \ref{thm:intro:main:eqrel}.
\end{proposition}
Observe that $\IQ^\times\cong\IZ/2\oplus\IZ^{\oplus\infty}$.
We do not need the action of $\IQ^\times$ on $(Y,\nu)$ to be
free. This way we see already that for any abelian group $\Lambda$
and any ergodic, measure preserving action of $\Lambda$ on $(Y,\nu)$,
we can realize $\Mod(\Centr_{\Autns(Y)}(\Lambda))$ as the fundamental group
of a type II$_1$ equivalence relation. In particular, we can do this
with $\Lambda=\IZ$. Using ergodic measures on $\IR$, Aaronson and Nadkarni
\cite{Aaronson:IntrinsicNormalizingConstants,AaronsonNadkarni:EigenvaluesAndSpectra}
show that for any number $0\leq\alpha\leq 1$, there is a measure
preserving action $\IZ\actson Y$ such that
$\Mod(\Centr_{\Autns(Y)}(\IZ))$ has Hausdorff dimension $\alpha$.
\begin{proof}
  We show that the action $\Gamma\curvearrowright
  (X,\mu)$ satisfies the conditions of theorem \ref{thm:intro:main:eqrel}.
  The second condition follows from the fact that $\SL_n\IQ$ has no
  non-trivial morphisms to abelian groups and that the conjugates of
  $\SL_n\IQ$ generate $\SL_n\IQ\ltimes\IQ^n$. For cocycle superrigidity,
  observe that $\IZ^n$ is an almost normal subgroup of $\SL_n\IQ\ltimes\IQ^n$ that
  has relative property (T) and that acts with infinite orbits on
  $I$. Popa's cocycle superrigidity theorem \cite[theorem
  0.1]{Popa:CocycleSuperrigidityMalleable} implies that
  $\ker\pi\curvearrowright X$ is $\Ufin$-cocycle superrigid.

  Let $\Delta:X\rightarrow X$ be a conjugation of
  $\ker\pi\curvearrowright X$. We show that $\Delta\in\Gamma$. By
  \cite[proposition 6.10]{Vaes:Bimodules}
  there is a measure preserving automorphism
  $\Delta_0:X_0\rightarrow X_0$ and a conjugation $\alpha:I\rightarrow
  I$ such that $\Delta(x)_{\alpha(i)}=\Delta_0(x_i)$. Since
  $(X_0,\mu_0)$ is atomic with unequal weights, $\Delta_0$ can only be
  the identity. Moreover, there is an automorphism
  $\delta:\SL_n\IQ\ltimes\IQ^n\rightarrow\SL_n\IQ\ltimes\IQ^n$ such that
  $\alpha(gi)=\delta(g)\alpha(i)$ for all $i\in I$ and
  $g\in\SL_n\IQ\ltimes\IQ^n$. As every affine subspace $V\subset I=\IQ^n$ is of
  the form $V=\Fix(g)$ for some $g\in\SL_n\IQ\ltimes\IQ^n$, the conjugation
  $\alpha$ maps affine subspaces of $\IQ^n$ to affine subspaces. Since
  $\alpha$ preserves inclusions, it also preserves the dimension and in
  particular, $\alpha$ maps lines to lines and planes to planes.
  Such a map is easily seen to be an affine transformation of $\IQ^n$,
  i.e. $\alpha\in \Gamma$. We obtained that $\Delta\in\Gamma$.
\end{proof}

\section{A preservation of Cartan result}
\label{sect:prescartan}
Let $\Gamma\actson (X,\mu)$ be a free, ergodic, measure preserving
action. For the orbit equivalence relation $\RelR=\Rel(\Gamma\actson
X)$, the previous section allows us to compute the group
$\Mod(\Autr(\RelR))$, in specific cases. For the group
measure space construction $M=\Lp^\infty(X)\rtimes\Gamma$, it is clear
that $\Mod(\Autf(M))$ contains at least $\Mod(\Autr(\RelR))$. Equality
holds if every automorphism $\psi:M\rightarrow M$ maps $\Lp^\infty(X)$
onto a unitary conjugate of $\Lp^\infty(X)$ inside $M$.

The following theorem provides a criterion to show, for some actions
$\Gamma\actson X$, that every automorphism $\psi$ of
$M=\Lp^\infty(X)\rtimes \Gamma$
preserves the Cartan subalgebra $\Lp^\infty(X)$ up to a unitary
conjugation in $M$. Its proof is a direct combination of results from
\cite{IoanaPetersonPopa:AFP}. The statement is very similar to
\cite[theorem 5.2]{PopaVaes:vNsuperrigidity}. Because the exact form
we need is not available in the literature, we provide a complete proof.

\begin{theorem}[see \cite{IoanaPetersonPopa:AFP}]
  \label{thm:prescartan}
  Let $\Gamma=\Gamma_1\free_\Sigma\Gamma_2$ be an amalgamated free
  product and consider a quotient morphism $\pi:\Gamma\rightarrow
  \Lambda$.
  Let $(B,\Tr)$ be a von Neumann algebra with a semifinite trace $\Tr$
  and assume that $B$ does not have minimal projections.
  Let $\beta;\Lambda\actson B$ be a $\Tr$-preserving action.
  Assume that
  \begin{itemize}
  \item $\Gamma_1\cap\ker\pi$ contains an infinite property (T) group $G$.
  \item the group $\Sigma$ is amenable and
    $\Gamma_1\not=\Sigma\not=\Gamma_2$.
  \item there is no $\ast$-homomorphism $\theta:\Lg(G)\rightarrow
    qBq$ for any projection $q\in B$
    with finite trace.
  \item there are
    group elements $g_1,\ldots,g_n\in\Gamma$ such that
    the intersection $\bigcap_{i=1}^n g_i\Sigma g_i^{-1}$ is finite.
  \end{itemize}

  Let $\alpha:\Gamma\actson (X,\mu)$ be any free, ergodic, p.m.p.\ action.
  Consider the action $\sigma:\Gamma\actson \Lp^\infty(X)\vnOtimes B$ given by
  $\sigma_g(a\otimes b)=\alpha_g(a)\otimes \beta_{\pi(g)}(b)$. Denote by $M$ the crossed product
  $M=(\Lp^\infty(X)\vnOtimes B)\rtimes\Gamma$.

  Then any $\ast$-automorphism $\psi$ of $M$ stably preserves
  $C=\Lp^\infty(X)\vnOtimes B$ up to a unitary, i.e.\ there is a unitary
  $u\in\Bounded(\ell^2(\IN))\vnOtimes\Bounded(\ell^2(\IN))\vnOtimes M$ such
  that
  \[
    u(\ell^\infty(\IN)\vnOtimes\Bounded(\ell^2(\IN))
    \vnOtimes\psi(C))u^\ast
    =\ell^\infty(\IN)\vnOtimes\Bounded(\ell^2(\IN))\vnOtimes C
  \]

  If $B$ is abelian, or if $B$ is properly infinite and the ergodic action
  $\Lambda\actson\Centre(B)$ does not preserve a probability measure
  (i.e. the orbit equivalence relation is of type I$_\infty$,
  II$_\infty$ or III),
  then there is a unitary $v\in M$ such that $v\psi(C)v^\ast=C$.
\end{theorem}
\begin{proof}
  Throughout the proof, we will use the intertwining-by-bimodules
  technique in the semifinite setting. Because
  \cite{IoanaPetersonPopa:AFP} is formulated a finite setting, need to
  be able to take finite corners. The following claim allows to do
  exactly that.

  Write $A=\Lp^\infty(X)$. Denote by $\cN=B\rtimes\Lambda$ and set
  $\cM=(\Lp^\infty(X)\rtimes\Gamma)\vnOtimes\cN$. There is a natural
  trace preserving embedding $\varphi$ of $M$ into $\cM$, given by the formula
  $\varphi((a\otimes b)u_g)=au_g\otimes bu_{\pi(g)}$.\\
  \emph{Claim: let $H\subset\Gamma$ be a subgroup, let $A_0$ be $\IC$
    or $A$, and let $Q\subset p M p$
    be a subalgebra of a finite corner of $M$, then we have that
    \[Q\embeds_{M}(A_0\otimes B)\rtimes H\text{ if and only
      if }\varphi(Q)\embeds_{\cM}(A_0\rtimes H)\otimes \cN.\]}
  If $Q$ embeds into $(A_0\otimes B)\rtimes H$ inside $M$, it is clear
  that $\varphi(Q)$ embeds into $\varphi((A_0\otimes B)\rtimes
  H)\subset(A_0\rtimes H)\otimes \cN$ inside
  $\varphi(M)\subset\cM$.

  Conversely, assume that $Q$ does not embed into $(A_0\otimes B)\rtimes
  H$ inside $M$. Then, by \cite[theorem 2.1]{HoudayerRicard:FreeArakiWoodsFactors},
  we find a sequence $(v_n)_n$ of unitaries in $Q$ such that
  \[\norm{\E_{(A_0\otimes B)\rtimes H}(xv_ny)}_2\rightarrow 0\text{ for
    all }x,y\in M.\]
  We have to show that
  \[\norm{\E_{(A_0\rtimes H)\vnOtimes \cN}(x\varphi(v_n)y)}_2\rightarrow
  0\text{ for all }x,y\in\cM.\]
  By Kaplanski's density theorem, it suffices to do this for $x=u_g$
  and $y=u_h$. Write the Fourier expansion of $v_n$ in $(A_0\otimes
  B)\rtimes \Gamma$ as $v_n=\sum_{k\in\Gamma}v_{n,k}u_k$. Then we
  compute that
  \begin{align*}
    \norm{\E_{(A_0\rtimes H)\vnOtimes \cN}(u_g\varphi(v_n)u_h}_2^2
    &=\sum_{k\in\Gamma, gkh\in H}\norm{\E_{A_0\otimes B}(v_{n,k})}_2^2\\
    &=\norm{\E_{(A_0\otimes B)\rtimes H}(u_gv_nu_h)}_2^2\rightarrow 0,
  \end{align*}
  proving that $\varphi(Q)$ does not embed into $(A_0\rtimes
  H)\otimes\cN$ inside $\cM$.

  Using this claim, we prove theorem \ref{thm:prescartan}. Let $\psi$
  be an automorphism of $M$ and let $p\in B$ be a projection with finite trace.
  Up to a unitary conjugation, we can assume that
  $q=\psi(p)$ is contained in $B$. We will use the following notations
  \begin{align*}
    \widetilde \cM_1&=(A\rtimes\Gamma_1)\vnOtimes q\cN q\\
    \widetilde \cM_2&=(A\rtimes\Gamma_2)\vnOtimes q\cN q\\
    \widetilde \cP&=(A\rtimes\Sigma)\vnOtimes q\cN q\\
    \widetilde \cM&=(A\rtimes\Gamma)\vnOtimes q\cN q = \widetilde
    \cM_1\free_{\widetilde\cP}\widetilde\cM_2\\
    \widetilde\psi&=\varphi\circ\psi:p M p\rightarrow \widetilde\cM
  \end{align*}
  \emph{Step 1: we show that $\widetilde\psi(A\otimes p B p)\embeds_{\widetilde\cM}\widetilde\cM_i$ for $i=1$ or $2$}
  As in \cite{IoanaPetersonPopa:AFP}, we consider the word-length
  deformation $(m_\rho)_\rho$ on $\widetilde
  \cM=\widetilde\cM_1\free_{\widetilde\cP}\widetilde\cM_2$, i.e.\ for
  every number $0<\rho<1$, we define a unital completely positive map
  $m_\rho:\widetilde \cM\rightarrow\widetilde\cM$ that is given by
  $m_\rho(au_g\otimes x)=\rho^{\len{g}}(au_g\otimes x)$ for all
  $a\in A, g\in\Gamma$ and $x\in q\cN q$. The notation $\len{g}$
  denotes the word-length of $g$ in the amalgamated free product
  $\Gamma=\Gamma_1\free_\Sigma\Gamma_2$. Observe that
  $\norm{m_\rho(x)-x}_2$ tends to $0$ as $\rho$ tends to 1, pointwise
  for all $x\in\widetilde \cM$.

  The subalgebra $Q=\widetilde\psi(p\Lg(G))$ has property (T), so
  the word-length deformation converges to $\id$ uniformly on the unit
  ball of $Q$. We also know that $\psi(p\Lg(G))$ does not embed into $B$
  inside $M$, so the claim above shows that $Q$ does not
  embed into $q\cN q$ inside $\widetilde \cM$. Because $A\rtimes\Sigma$ is amenable,
  it follows that $Q$ does not embed into $\widetilde \cP$ inside
  $\widetilde\cM$ (see \cite{Popa:Betti}).
  The unitaries $\widetilde\psi(pu_g)\in\widetilde\cM$ with $g\in G$,
  normalize the abelian
  von Neumann algebra $\widetilde\psi(pA)\subset \widetilde \cM$. By \cite[lemma
  5.7]{PopaVaes:vNsuperrigidity} (which is a version of
  \cite[proposition 1.4.1]{IoanaPetersonPopa:AFP}), we find a
  $0<\rho_0<1$ and a $\delta>0$ such that $\Tr(w^\ast
  m_{\rho_0}(w))\geq\delta$ for all unitaries $w\in
  \widetilde\psi(p(A\rtimes G))$.

  We also know that $\widetilde\psi(p(A\rtimes G))$ does not embed into
  $\cP$ inside $\cM$, so \cite[theorem 5.4]{PopaVaes:vNsuperrigidity}
  (which is a version of \cite[theorem 4.3]{IoanaPetersonPopa:AFP})
  implies that $\vnNorm_{\widetilde\cM}(\widetilde\psi(p(A\rtimes
  G)))^{\prime\prime}$ embeds into $\widetilde\cM_i$, for $i=1$ or $2$. This
  normalizer contains $\widetilde\psi(A\otimes p B p)$.

  \emph{We finish the proof of theorem \ref{thm:prescartan}.}\\
  We have just shown that $\widetilde \psi(A\otimes p B p)$ embeds into
  $\widetilde \cM_i$ inside $\widetilde \cM$. The claim in the
  beginning of this proof shows that $\psi(A\otimes p B p)$ embeds into
  $M_1$ inside $M$. Because $\ker\pi$ contains the infinite property
  (T) group $G$, it is certainly not contained in the amenable group
  $\Sigma$. By our last condition, we can apply proposition \ref{prop:inter:free}
  (which is a variant of \cite[proposition 8]{HoudayerPopaVaes:AppvNsuperrigidity})
  to obtain that $\psi(A\otimes p B p)\embeds_M A\otimes B$.

  If $B$ happens to be abelian, then $A\otimes B$ is a Cartan
  subalgebra of $M$, so \cite[proposition 2.3]{HoudayerPopaVaes:AppvNsuperrigidity}
  implies that $\psi(A\otimes B)$ is unitarily conjugate to $A\otimes
  B$.

  Otherwise, we have just shown that $\psi(A\otimes p B p)\embeds_M
  A\otimes B$, and by symmetry it follows that $A\otimes
  qBq\embeds_M\psi(A\otimes B)$. Proposition \ref{prop:inter:regular}
  below shows that
  $\psi(A\otimes B)$ is stably unitarily conjugate to $A\otimes B$
  inside $M$, i.e.\ there is a unitary
  $u\in\Bounded(\ell^2(\IN))\otimes\Bounded(\ell^2(\IN))\otimes M$
  such that
  \[u(\ell^\infty(\IN)\otimes\Bounded(\ell^2(\IN))\otimes
  \psi(A\otimes B))u^\ast=\ell^\infty(\IN)\otimes\Bounded(\ell^2(\IN))\otimes
  A\otimes B.\]

  If $B$ is properly infinite and there is no finite
  $\Lambda$-invariant measure on $\Centre(B)$, proposition
  \ref{prop:inter:regular} below shows that $\psi(A\otimes B)$ is actually unitarily
  conjugate to $A\otimes B$.
\end{proof}

Popa's intertwining-by-bimodules technique is used to decide if two
subalgebras $A,B\subset M$ of a II$_1$ factor are unitarily
conjugate. If $A$ and $B$ are unitarily conjugate, the it is clear
that $A\embeds_M B$ and $B\embeds_M A$. The converse implication is
false in general. However, there are several special cases in which
this converse implication holds.

For example, if $A$ and $B$ are Cartan subalgebras of $M$ and
$A\embeds_M B$, then \cite[theorem A.1]{Popa:Betti} (or
\cite[proposition 2.3]{HoudayerRicard:FreeArakiWoodsFactors} in the semifinite case) shows
that $A$ and $B$ are unitarily conjugate. We need a criterion that
allows non-abelian algebras $A$ and $B$. Such a criterion is given below.

\begin{proposition}
  \label{prop:inter:regular}
  Let $M$ be a type \II$_\infty$ factor and let $A,B\subset M$ be von
  Neumann subalgebras such that the restricted traces
  $\Tr\restrict{A}$ and $\Tr\restrict{B}$ are still
  semifinite. Suppose that
  \begin{itemize}
  \item $A$ and $B$ are regular subalgebras of $M$.
  \item $\Centre(A)^\prime\cap M=A$ and $\Centre(B)^\prime\cap M=B$.
  \end{itemize}
  If $pAp\embeds_M B$ and $qBq\embeds_M A$ for some projection $p\in
  A$ and $q\in B$ with finite trace, then $A$ and $B$ are stably
  unitarily conjugate in $M$, in the sense that
  \[A\otimes \Bounded(\ell^2(\IN))\otimes\ell^\infty(\IN)\quad\text{ and }\quad
  B\otimes \Bounded(\ell^2(\IN))\otimes\ell^\infty(\IN)\]
  are unitarily conjugate subalgebras of
  $M\otimes\Bounded(\ell^2(\IN))\otimes \Bounded(\ell^2(\IN))$
  
  If moreover $A$ and $B$ satisfy
  \begin{itemize}
  \item $A$ and $B$ are properly infinite, i.e.\ every central
    projection in $A$ respectively $B$ is an infinite projection in $A$ respectively
    $B$,
  \item there is no finite trace $\tau_{\Centre(A)}$ on $\Centre(A)$
    that is invariant under the action of $\vnNorm_{M}(A)$.
    Similarly there is no $\vnNorm_M(B)$-invariant finite trace
    $\tau_{\Centre(B)}$ on $\Centre(B)$,
  \end{itemize}
  then there is a unitary $u\in M$ with $u^\ast Au=B$
\end{proposition}
\begin{proof}
  \emph{Step 1: There are projections $p_1\in A$ and $q_1\in B$ such
    that $p_1\Lp^2(M)q_1$ contains a nonzero finite index
    $p_1Ap_1$-$q_1Bq_1$ subbimodule.}\\
  Consider the $A$-$B$ bimodules
  \begin{align*}
    \vnInterB(A\embeds_M\!B)&\!=\!\bigvee\!\!\left\{\!K\!\subset\! p_1\Lp^2(M)\,\left\vert
        \begin{aligned}
          &p_1\in A\text{ is a projection with }\Tr(p_1)<\infty\text{,
            and}\\
          &K\text{ is a }p_1Ap_1\text{-}B\text{ bimodule with }\dim_{\Tr}(\bimod{}{K}{B})\!<\!\infty
        \end{aligned}
        \!\right.\right\}\\
    \vnInterB(A\oembeds_M\!B)&\!=\!\bigvee\!\!\left\{\!K\!\subset\! \Lp^2(M)q_1\,\left\vert
        \begin{aligned}
          &q\in B\text{ is a projection with }\Tr(q_1)<\infty\text{,
            and}\\
          &K\text{ is an }\!A\text{-}q_1Bq_1\text{ bimodule with }\!\dim_{\Tr}(\bimod{A}{K}{})\!<\!\infty
        \end{aligned}
        \!\right.\right\}
  \end{align*}
  Because $A,B\subset M$ are regular subalgebras,
  both bimodules are actually non-zero $M$-$M$ bimodules, so it follows that
  $\vnInterB(A\embeds_MB)=\Lp^2(M)=\vnInterB(A\oembeds_MB)$. In particular,
  we find that, for any projections $p_1\in A$ and $q_1\in B$ with
  finite trace, the Hilbert space $p_1\Lp^2(M)q_1$ is spanned by $p_1Ap_1$-$q_1Bq_1$
  bimodules that have finite dimension over $q_1Bq_1$. Similarly, $p_1\Lp^2(M)q_1$
  is spanned by $p_1Ap_1$-$q_1Bq_1$ bimodules that have finite dimension
  over $p_1Ap_1$. Write $C$ for the von Neumann algebra of $p_1Ap_1$-$q_1Bq_1$
  bimodular operators on $p_1\Lp^2(M)q_1$. On $C$, we have two natural
  traces: a trace $\Tr_A$ that comes from the canonical trace on the
  commutant of $p_1Ap_1$, and a trace $\Tr_B$ that comes from the
  canonical trace on the commutant of $q_1Bq_1$. The fact that
  $\vnInterB(A\embeds_M B)=\Lp^2(M)=\vnInterB(A\oembeds_M B)$ implies that
  both traces $\Tr_A$ and $\Tr_B$ are semifinite traces on $C$. In
  particular, we find elements $x,y\in C$ such that $\Tr_A(x^\ast
  x)<\infty$, such that $\Tr_B(y^\ast y)<\infty$ and such that
  $xy\not=0$. Then we also
  have that $\Tr_A(x^\ast y^\ast yx)<\infty$ and $\Tr_B(x^\ast
  y^\ast yx)<\infty$, so any nonzero spectral projection
  $r=\chara_{[\delta,\infty[}(y^\ast x^\ast x y)$ is the projection onto a
  nonzero finite index $p_1Ap_1$-$q_1Bq_1$ subbimodule $K$ of $p_1\Lp^2(M)q_1$.

  \emph{Step 2: Write $A_1=p_1Ap_1$ and $B_1=q_1Bq_1$. We write the
    left action of $A_1$ on $K$ explicitly by $\lambda:A_1\rightarrow
    \Bounded(K)$ and we write the right action by
    $\rho:B_1^{\text{op}}\rightarrow\Bounded(K)$. Then there is a non-zero
    $A_1$-$B_1$ bimodular projection $r$ on $K$ such that
  \[r\lambda(\Centre(A_1))=r(\bimod{A_1}{\Bounded}{B_1}(K))r=r\rho(\Centre(B_1)).\]}\\
  Denote by
  $C=\bimod{A_1}{\Bounded}{B_1}(K)$ the von Neumann algebra of all
  $A_1$-$B_1$ bimodular operators on $K$, so
  $C=\lambda(A_1)^\prime\cap\rho(B_1^{\text{op}})^\prime$. Since $K$ is
  a finite index bimodule, we know that
  $\lambda(A_1)\subset\rho(B_1^{\text{op}})^\prime$ is a finite index
  inclusion, so $\lambda(\Centre(A_1))$ is a finite index subalgebra
  of $C$, see for example \cite[lemma A.3]{Vaes:Bimodules}. Because
  $\Centre(A_1)$ is an abelian von Neumann algebra, we find a
  nonzero projection $r_1\in C$ that satisfies
  $r_0Cr_0=\lambda(\Centre(A_1))p$. Similarly, we see that
  $\rho(\Centre(B_1))r_0$ is a finite index subalgebra of $r_0Cr_0$, so we
  find a nonzero projection $r\leq r_0$ in $C$ such that
  $rCr=\rho(\Centre(B))r$.

  \emph{Step 3: There is a partial isometry $v\in M$ such that
    $vv^\ast\in A$ and $v^\ast v\in B$, both having finite trace, and
    such that $v^\ast A v=v^\ast vBv^\ast v$.}
  Possibly passing to a subbimodule, we can assume that $\widetilde
  K=rK$ is finitely generated over $B_1$. So we find an $A_1$-$B_1$
  bimodular unitary
  \begin{equation*}
    \label{eqn:bimod}
    V:\bimod{\psi(A_1)}{q_2(\IC^n\otimes\Lp^2(B_1))}{B_1}\rightarrow \bimod{A_1}{\widetilde K}{B_1}
  \end{equation*}
  for some projection $q_2\in B_1^n$ and a finite index inclusion
  $\psi:A_1\rightarrow q_2 B^n_1q_2$. The conclusion of step 2 shows
  that $\psi(\Centre(A_1))=\psi(A_1)^\prime\cap q_2B_1^nq_2=\Centre(q_2B_1^nq_2)$.

  Consider the vector $\xi=[V(e_1\otimes \hat1),\ldots,V(e_n\otimes
  \hat1)]\in\IC^n\otimes\Lp^2(M)$. Observe that this vector satisfies
  the relation $x\xi=\xi\psi(x)$ for all $x\in A_1$. Polar
  decomposition gives a partial isometry $v\in \MatM_{1,n}\otimes M$
  such that $xv=v\psi(x)$ for all $x\in A_1$. Moreover, $vv^\ast\leq
  p_1$ has finite trace. Remark that $vv^\ast$
  commutes with $A_1$, so $vv^\ast$ is contained in $A_1\subset
  A$. Also, $v^\ast v$ commutes with $\psi(A_1)$. In particular, it
  commutes with $\psi(\Centre(A_1))=\Centre(q_2B_1^nq_2)$. Hence we
  have that $v^\ast v\in q_2B_1^nq_2\cap
  \psi(A_1)^\prime=\Centre(q_2B_1^nq_2)$. Possibly making $p_1$, $q_1$ and
  $q_2$ smaller, we can assume that $n=1$, $p_1=vv^\ast$ and $v^\ast
  v=q_2=q_1$. Now we have that
  $v^\ast\Centre(p_1Ap_1)v=\Centre(q_1Bq_1)$. Our first condition
  shows that $v^\ast Av=v^\ast v Bv^\ast v$.

  \emph{Step 4: The subalgebras $A$ and $B$ are stably unitarily conjugate.}
 We use the following
  notations:
  \begin{align*}
    \cA&=A\otimes\Bounded(\ell^2(\IN))
    & \widetilde\cA&=\cA\otimes\ell^\infty(\IN) = A\otimes\Bounded(\ell^2(\IN))\otimes\ell^\infty(\IN)\\
    \cB&=B\otimes\Bounded(\ell^2(\IN))
    & \widetilde\cB&=\cB\otimes\ell^\infty(\IN) = B\otimes\Bounded(\ell^2(\IN))\otimes\ell^\infty(\IN)\\
    \cM&=M\otimes\Bounded(\ell^2(\IN))
    & \widetilde\cM&=\cM\otimes\Bounded(\ell^2(\IN)) = M\otimes\Bounded(\ell^2(\IN))\otimes\Bounded(\ell^2(\IN))
  \end{align*}
  \emph{Step 4.1: there is a partial isometry $w\!\in\!
    \cM$ with $ww^\ast\!\in\!
    \Centre(\cA)$, $w^\ast w\!\in\!
    \Centre(\cB)$ and such that}
  \[w^\ast \cA w=w^\ast
  w\cB.\]
  We find an infinite sequence $(v_n)_n$ of partial isometries in
  $\cA$ such that $v_nv_n^\ast=p\otimes e_{1,1}$ for all
  $n$, and such that $\sum_n v_nv_n^\ast$ is the central support of
  $p\otimes e_{1,1}$ in $\cA$. Similarly,
  we find an infinite sequence $(\tilde v_n)_n$ of partial
  isometries in $\cB$ with $\tilde
  v_n\tilde v_n^\ast=q\otimes e_{1,1}$ for all $n$, and such that
  $\sum_n \tilde v_n^\ast\tilde v_n$ is a central projection in
  $\cB$.
  The partial isometry $w=\sum_n v_n v\tilde v_n^\ast$ satisfies the
  conditions of step 4.1.
  
  \emph{Step 4.2: there is a unitary $u\in\widetilde \cM$
    such that}
  \[u^\ast\widetilde \cA u= \widetilde \cB.\]
  Because $\widetilde \cA$ is a
  regular subalgebra of the factor $\widetilde \cM$, we
  see that the normalizer of $\widetilde\cA$ acts ergodically on
  $\Centre(\widetilde \cA)$. We
  find an infinite sequence $(w_n)_n$ of partial isometries in
  $\widetilde \cM$ with $w_nw_n^\ast=ww^\ast$ and $w_n^\ast
  w_n\in\Centre(\widetilde\cA)$, such that $w_n^\ast \widetilde \cA w_n=w_n^\ast
  w_n\widetilde \cA$ and such that $\sum_nw_n^\ast w_n=1$. Similarly,
  we find an infinite sequence $(\tilde w_n)_n$ of partial
  isometries in $\widetilde\cM$ with $\tilde w_n\tilde
  w_n^\ast=w^\ast w$ and $\tilde w_n^\ast\tilde w_n\in \Centre(\widetilde\cB)$,
  such that $\sum_n \tilde w_n^\ast \tilde w_n=1$ and such that
  $\tilde w_n^\ast \widetilde \cB\tilde w_n=\tilde w_n^\ast\tilde
  w_n\widetilde \cB$. The unitary $u=\sum_n w_n w \tilde w_n^\ast$
  satisfies $u\widetilde Au^\ast=\widetilde B$.

  \emph{Step 5: If $A$ and $B$ are properly infinite and $\Centre(A)$
    (respectively $\Centre(B)$) does not admit a finite
    $\Norm_M(A)$-invariant (resp.\ $\Norm_M(B)$-invariant) trace, then
    $A$ and $B$ are unitarily conjugate.}
  Take a unitary $u\in
  \Bounded(\ell^2(\IN))\otimes\Bounded(\ell^2(\IN))\otimes M$ that
  conjugates $\ell^\infty(\IN)\otimes\Bounded(\ell^2(\IN))\otimes A$
  onto $\ell^\infty(\IN)\otimes\Bounded(\ell^2(\IN))\otimes B$.
  The fact that $A$ is properly infinite means exactly that there is a
  partial isometry $w_1\in \MatM_{1,\infty}(\IC)\otimes A$ such that
  $w_1w_1^\ast=1$ and $w_1^\ast w_1=1\otimes 1\in
  \Bounded(\ell^2(\IN))\otimes A$. There is a similar isometry $\tilde
  w_2\in \MatM_{1,\infty}(\IC)\otimes B$. The unitary $u_1=(1\otimes \tilde
  w_1)u(1\otimes w_1^\ast)$ in $\Bounded(\ell^2(\IN))\otimes M$
  conjugates $\ell^\infty(\IN)\otimes A$ onto $\ell^\infty(\IN)\otimes B$.

  Because $\Centre(A)$ is abelian, we can identify it with the
  $\Lp^\infty(X)$ of some measure space $X$. We know that the ergodic
  action $\Norm_M(A)\actson X$ does not admit a finite invariant
  measure on $X$. So the associated orbit equivalence relation is of
  type I$_\infty$, II$_\infty$ or III. In any case there is a non-singular
  isomorphism $\Delta:\IN\times X\rightarrow X$ such that $\Delta(i,x)$
  is in the same orbit as $x$, for all $i\in \IN$ and almost all $x\in X$.
  This yields a unitary $w_2\in \MatM_{1,\infty}(\IC)\otimes M$ such
  that $w_2^\ast Aw_2=\ell^\infty(\IN)\otimes A$. We find a similar
  unitary $\tilde w_2$ for $B$. Now that unitary $u_2=\tilde
  w_2u_1w_2^\ast\in M$ conjugates $A$ onto $B$.
\end{proof}

In both steps 4.1 and 4.2 of the proof of proposition \ref{prop:inter:regular}, we have
to pass to an amplification. The following example explains why we can
not avoid this.
Let $A_0\subset M_0$ be any Cartan subalgebra of a II$_1$ factor. Set
\begin{align*}
  M&=M_0\otimes\MatM_2(\IC)\\
  A&=A_0\otimes\MatD_2(\IC)\\
  B&=A_0\otimes\MatM_2(\IC)
\end{align*}
With $v=1\otimes e_{1,1}$, we clearly have that $vAv^\ast = vv^\ast
B vv^\ast$. Of course $v$ does not extend to a unitary that
conjugates $A$ onto $B$ ($A$ and $B$ are not even isomorphic).

We really need the ``abelian'' amplification in step 4.1:
$A\!\otimes\!\Bounded(\ell^2(\IN))$ and $B\otimes\Bounded(\ell^2(\IN))$
are not unitarily conjugate in $M\otimes\Bounded(\ell^2(\IN))$
because the smallest projection $p\in A\otimes\Bounded(\ell^2(\IN))$
that has full central support in $A\otimes\Bounded(\ell^2(\IN))$,
has $\Tr(p)=2$ while $q=1\otimes
e_{1,1}\otimes e_{1,1}\in B\otimes\Bounded(\ell^2(\IN))$ has full central
support in $B\otimes\Bounded(\ell^2(\IN))$, and $\Tr(q)=1$.

Also the ``factorial'' amplification in step 4.2 is necessary:
the von Neumann algebras $A\otimes\ell^\infty(\IN)$ and
$B\otimes\ell^\infty(\IN)$ are not even isomorphic.

\section{Every fundamental group of an equivalence relation is the
  fundamental group of a factor}
\label{sect:SeqrelInSfactor}
Given a II$_1$ equivalence relation $\RelR$, we construct a type
II$_1$ factor $M$ whose fundamental group is
$\Fundg(M)=\Fundg(\RelR)$. For a type II$_1$ equivalence relation
$\RelR$ on a probability space $(X,\mu)$, denote by $\RelR^\infty$ the
infinite amplification on $X\times \IN$. Given
any type II$_1$ factor $Q$, we can construct a new type II$_\infty$
factor $\widetilde M=(Q\otimes \Lp^\infty(X\times
\IN))\free_{\Lp^\infty(X\times \IN)}\Lg(\RelR^\infty)$. For
convenience, we will denote $A=\Lp^\infty(X\times\IN)$ and
$P=\Lg(\RelR^\infty)$.

It is now clear that $\Fundg(\RelR)\subset \Mod(\Autf(\widetilde
M))=\Fundg(p\widetilde Mp)$ for any projection $p\in \widetilde M$ with finite
trace. We show that, for the right choice of $Q$, we actually have
$\Fundg(\RelR)=\Fundg(p\widetilde Mp)$.
\begin{theorem}
  \label{thm:SeqrelInSfactor}
  Let $\RelR$ be any type II$_1$ equivalence relation on
  $(X,\mu)$. Denote by $P=\Lg(\RelR^\infty)$ the generalized
  group--measure space construction of the infinite amplification of
  $\RelR$. Denote by $A=\Lp^\infty(X\times\IN)$ the corresponding
  Cartan subalgebra $A\subset P$.

  Let $Q$ be a type II$_1$ factor with trivial fundamental group and
  such that there exists a diffuse Cartan subalgebra $Q_0\subset Q$
  with relative property (T). Define a type II$_\infty$ factor
  $M=(Q\otimes A)\free_A P$. Then we have
  \[\Mod(\Autf(M))=\Fundg(\RelR).\]
\end{theorem}
Popa showed in \cite{Popa:Betti} that the type II$_1$ factor
$Q=\Lg(\SL_2\IZ\ltimes\IZ^2)$ satisfies the conditions of theorem
\ref{thm:SeqrelInSfactor}. In this case, the Cartan subalgebra $Q_0=
\Lg(\IZ^2)$ has the relative property (T) in $Q$.
\begin{proof}
  It is obvious that the fundamental group of $\RelR$ is included in
  $\Mod(\Autf(M))$. For the other inclusion, it is sufficient to show
  that every automorphism of $M$ preserves $A$ up to a unitary in $M$.

  The techniques to do this were developed  by Ioana, peterson and Popa
  in \cite{IoanaPetersonPopa:AFP}. In fact, we need a semifinite
  generalization of these results, but the original proofs generalize
  in a straightforward way. A complete proof for these generalizations
  is also given in my thesis \cite{Deprez:PhDthesis}. For each result
  we use, we give both relevant references.

  Let $\psi:M\rightarrow M$ be an automorphism of $M$, and fix a
  projection $p\in A$ with finite trace. Consider the word-length
  deformation $m_\rho:M\rightarrow M$ as in
  \cite{IoanaPetersonPopa:AFP}. Because $Q_0\subset Q$ has the
  relative property (T), we know that $m_\rho$ converges uniformly on
  the unit ball of $\psi(pQ_0)$.

  We show that $\psi(pQ_0)\nembeds_M
  A$. If this were not the case,  \cite[lemma 3.5]{Vaes:Bimodules}
  would yield a projection $q\in A$ with finite trace and such that
  \[qQ\subset q(M\cap A^\prime)\embeds_M \psi(p(M\cap Q_0^\prime)p).\]
  But by \cite[theorem 1.2.1]{IoanaPetersonPopa:AFP} (generalized to the
  semifinite setting in \cite[theorem 2.4]{ChifanHoudayer:BassSerreRigidity}),
  we know that the relative commutant of $pQ_0$ is in fact equal to
  $p(Q_0\otimes A)$. Then the type II$_1$ factor $qQ$ embeds into the
  abelian von Neumann algebra $\psi(p(Q_0\otimes A))$. This
  contradiction shows that $\psi(pQ_0)\nembeds_M A$.

  But all the unitaries of $\psi(pQ_0)$ commute with the abelian von Neumann
  algebra $\psi(pA)$, so \cite[proposition 1.4.4 and theorem
  4.3]{IoanaPetersonPopa:AFP} (or rather \cite[lemma 4.8 and theorem
  4.6]{Deprez:PhDthesis}) shows that $\psi(pA)$ embeds into $(Q\otimes A)$
  or into $P$. Since $\psi(pA)$ is regular in $\psi(p)M\psi(p)$, it follows that
  $\psi(pA)\embeds_M A$, see \cite[theorem 1.2.1]{IoanaPetersonPopa:AFP}
  (generalized to the
  semifinite setting in \cite[theorem
  2.4]{ChifanHoudayer:BassSerreRigidity}). Now, \cite[lemma
  3.5]{Vaes:Bimodules} shows that $q(M\cap A^\prime)\embeds \psi(M\cap
  A^\prime)$, for some projection $q\in A$.

  We can apply the same argument to $\psi^{-1}$, and we obtain that also
  $\psi(r(M\cap A^\prime))\embeds_M M\cap A^\prime$. Now proposition
  \ref{prop:inter:regular} shows that $\psi(M\cap A^\prime)$ is stably
  unitarily conjugate to $M\cap A^\prime$. So we find a unitary
  $u\in \Bounded(\ell^2(\IN))\otimes\Bounded(\ell^2(\IN))\otimes M$
  conjugating $\Bounded(\ell^2(\IN))\otimes\ell^\infty(\IN)\otimes
  \psi(M\cap A^\prime)$ onto $\Bounded(\ell^2(\IN))\otimes\ell^\infty(\IN)\otimes
  (M\cap A^\prime)$.

  Observe that $M\cap A^\prime$ is isomorphic to
  $Q^{\free\infty}\otimes A$. By \cite[theorem
  6.3]{IoanaPetersonPopa:AFP}, we know that $\free^\infty Q$ has
  trivial fundamental group. So we can assume that $u\in
  1\otimes\Bounded(\ell^2(\IN))\otimes M$.  Since $\RelR^\infty$ is a type II$_\infty$
  equivalence relation, we can assume that $u\in 1\otimes 1\otimes M$,
  so $\psi(M\cap A^\prime)$ is unitarily conjugate to $M\cap
  A^\prime$ inside $M$. The same is true for their respective centers, $\psi(A)$
  and $A$.
\end{proof}

\section{The II\texorpdfstring{$_1$}{1} factor case}
\label{sect:factor}
We want to prove a result similar to theorem
\ref{thm:intro:main:eqrel}, but for type II$_1$ factors instead of
equivalence relations. To do this, we combine theorem
\ref{thm:intro:main:eqrel} with theorem \ref{thm:prescartan}.

Let $\Gamma\actson (X,\mu)$ be a free, ergodic, p.m.p.\ action, let
$\pi:\Gamma\rightarrow\Lambda$ be a quotient morphism and let
$\Lambda\actson (Y,\nu)$ be an ergodic, infinite measure preserving
action. Define a new action $\Gamma\actson X\times Y$ by
$g(x,y)=(gx,\pi(g)y)$. Denote the corresponding group measure space
construction by $M=\Lp^\infty(X\times Y)\rtimes \Gamma$. Any
$\Delta\in\Centr_{\Autns(Y)}(\Lambda)$ defines an automorphism
$\psi_\Delta:M\rightarrow M$ by the formula $\psi_\Delta((a\otimes
b)u_g)=(a\otimes\Delta_\ast(b))u_g$. Two different automorphisms like
that are never unitarily conjugate. However, there is another natural
class of automorphisms of $M$. Let $\omega:\Lambda\times
Y\rightarrow\Circle^1$ be a 1-cocycle. Define unitaries
$v_s\in\Lp^\infty(Y)$ by $v_s(y)=\omega(s,s^{-1}y)$. Then
the formula $\varphi_\omega((a\otimes b)u_g)=(a\otimes
bv_{\pi(g)})u_g$ defines an automorphism of $M$. Two such
automorphisms $\varphi_{\omega}$ and $\varphi_{\tilde\omega}$ are
unitarily conjugate if and only if $\omega$ and $\tilde\omega$ are
cohomologous. Denote by $\Cohom(\Lambda\actson Y)$ the group of all
1-cocycles $\omega: \Lambda\times Y\rightarrow\Circle^1$, identifying
cohomologous cocycles.

The preceding paragraph shows that $\Out(M)$ contains at least the
group
\begin{equation}
  \label{eqn:centr-cohom}
  \Centr_{\Autns(Y,\nu)}(\Lambda)\ltimes\Cohom(\Lambda\actson Y)\subset \Out(M).
\end{equation}
Of course, the $\varphi_\omega$ are trace preserving, so they do not
contribute to $\Mod(\Autns(M))$.

Theorem \ref{thm:proof:main:factor} shows that under strong
conditions on $\Gamma\actson X$, the inclusion in
(\ref{eqn:centr-cohom}) becomes and equality. Proposition
\ref{prop:ex:factor} gives an explicit example of an action
$\Gamma\actson (X,\mu)$ that satisfies the conditions of theorem
\ref{thm:proof:main:factor}. As with the equivalence relation case,
the quotient group is $\Lambda=\IQ^\times$. By
\cite{Aaronson:IntrinsicNormalizingConstants,AaronsonNadkarni:EigenvaluesAndSpectra},
this gives explicit examples of type
II$_1$ factors whose fundamental group can have any Hausdorff
dimension $0\leq \alpha\leq 1$.

\begin{theorem}
  \label{thm:proof:main:factor}
  Let $\Gamma\actson (X,\mu)$ be a free, ergodic and probability
  measure preserving action, and let $\pi:\Gamma\rightarrow\Lambda$ be
  a quotient morphism. Suppose that
  \begin{itemize}
  \item the conditions of theorem \ref{thm:intro:main:eqrel} are
    satisfied.
  \item the conditions of theorem \ref{thm:intro:main:factor} are
    satisfied.
  \item $\ker\pi$ is a perfect group.
  \end{itemize}

  Let $\Lambda\actson (Y,\nu)$ be an ergodic measure preserving
  action. Consider the action of $\Gamma$ on $X\times Y$ given by
  $g(x,y)=(gx,\pi(g)y)$, and denote the crossed product by
  $M=\Lp^\infty(X\times Y)\rtimes\Gamma$. Then the outer automorphism
  group of $M$ is
  \[\Out(M)=\Centr_{\Autns(Y,\nu)}(\Lambda)\ltimes\Cohom(\Lambda\actson Y).\]
\end{theorem}
\begin{proof}
  By the discussion above, it is clear that the group in the right
  hand side is contained in the outer automorphism group of $M$.

  Let $\psi:M\rightarrow M$ be an automorphism of $M$. By theorem
  \ref{thm:prescartan}, we can assume that $\psi$ globally
  preserves $\Lp^\infty(X\times Y)$. Theorem 5.1 in
  \cite{Singer:AutomorphismsOfFiniteFactors} yields an orbit equivalence $\Delta:X\times
  Y\rightarrow X\times Y$ such that $\psi(a\otimes
  b)=\Delta_\ast(a\otimes b)$ for all
  $a\otimes b\in A\otimes B$. By theorem \ref{thm:intro:main:eqrel}, we can
  assume that $\Delta$ is of the form $\id\times\Delta_0$ for some
  $\Delta_0\in\Centr_{\Autns(Y)}(\Lambda)$.

  In other words, we have that $\tilde\psi=\psi_{\Delta_0}^{-1}\circ\psi$ acts as
  the identity on $A\otimes B$. Theorem 3.1 in \cite{Singer:AutomorphismsOfFiniteFactors}
  yields a cocycle $\omega:\Gamma\times X\times Y\rightarrow
  \Circle^1$ such that $\tilde\psi$ can be described in the following
  way. Denote by $v_g\in\Lp^\infty(X\times Y)$ the unitary that is
  defined by $v_g(x,y)=\omega(g,g^{-1}x,\pi(g^{-1}y))$, for every
  $g\in\Gamma$. Then $\tilde\psi$ is given by $\tilde\psi((a\otimes
  b)u_g)=(a\otimes b)v_gu_g$.

  Remember that $\ker\pi$ acts $\Ufin$-cocycle superrigidly on $X$, so
  we can assume that $\omega(g,x,y)$ does not depend on the
  $x$-variable, at least whenever $g\in\ker\pi$. Since $\ker\pi$ acts
  weakly mixingly on $X$, \cite[proposition
  3.6]{Popa:CocycleSuperrigidityMalleable} tells us that $\omega(g,x,y)$ is independent of $x$
  for all $g\in\Gamma$. From now on we consider $\omega:\Gamma\times
  Y\rightarrow \Circle^1$. For almost every $y\in Y$, the map
  $\ker\pi\ni g\mapsto\omega(g,y)\in\Circle^1$ is a group
  morphism. Since $\ker\pi$ is a assumed to be a perfect group, it
  follows that $\omega(g,y)=1$ for all $g\in\ker\pi$ and almost all
  $y\in Y$. This shows that $\omega$ splits over the quotient
  $\pi$. We consider $\omega$ as a cocycle for the action
  $\Lambda\actson Y$. We have shown that $\tilde\psi=\varphi_{\omega}$.

  For any automorphism $\psi:M\rightarrow M$, we have found a
  unitary $u\in M$, an automorphism
  $\Delta_0\in\Centr_{\Autns(Y)}(\Lambda)$ and a cocycle
  $\omega\in\Cohom(\Lambda\actson Y)$ such that $\psi$ is the composition
  $\Ad_u\circ\psi_{\Delta_0}\circ\varphi_\omega$. So we have shown
  that
  \[\Out(M)=\Centr_{\Autns(Y)}(\Lambda)\ltimes\Cohom(\Lambda\actson Y).\]
\end{proof}

\begin{proposition}
  \label{prop:ex:factor}
  Consider the group $\Gamma=\GL_{n+1}\IQ\free_{\T_{n}\IQ\ltimes \IQ^{n}}
  (\GL_{n}\IQ\ltimes\IQ^{n})$ ($n\geq 3$) where $\T_{n}\IQ$ is the
  group of upper triangular matrices in $\GL_{n}\IQ$. We view
  $\T_{n}\IQ\ltimes\IQ^{n}$ as the subgroup of $\GL_{n+1}\IQ$ consisting
  of matrices of the form \sMatrix{A&v\\0&1} with $A\in \T_{n}\IQ$
  and $v\in \IQ^{n}$. Define a quotient map $\pi:\Gamma\rightarrow
  \IQ^\times$ by $\pi(A)=\det(A)$ for all $A\in\GL_{n+1}\IQ$ and
  $\pi(A,v)=\det(A)$ when $(A,v)\in\GL_{n}\IQ\ltimes\IQ^{n}$. Set
  $H=\GL_{n}\IQ\subset \GL_{n}\IQ\ltimes\IQ^{n}$ and consider the
  generalized Bernoulli action $\Gamma\curvearrowright
  (X,\mu)=(X_0,\mu_0)^{\Gamma/H}$ over an atomic base space with
  unequal weights. For any ergodic action $\IQ^\times\curvearrowright
  (Y,\nu)$ preserving the infinite non-atomic measure $\nu$. Define an
  action $\Gamma\actson X\times Y$ by the formula
  $g(x,y)=(gx,\pi(g)y)$, and denote the crossed product by
  $M=\Lp^\infty(Y,\nu)\rtimes\Gamma$.
  Then we have that
  \[\Out(M)=\Centr_{\Autns(Y)}(\IQ^\times)\ltimes\Cohom(\IQ^\times\actson
  Y).\]
  It follows that
  \[\Mod(\Autf(M))=\Mod(\Centr_{\Autns(Y)}(\IQ^\times)).\]
\end{proposition}
\begin{proof}
  The group $\Gamma$ clearly satisfies the conditions of theorem
  \ref{thm:intro:main:factor} with $G=\SL_n\IZ\ltimes \IZ^n$,
  so it suffices to check the conditions
  of theorem \ref{thm:intro:main:eqrel}. There are no non-trivial group
  morphisms $\theta:\SL_{n+1}\IQ\rightarrow \IQ^\times$ nor
  $\theta:\SL_n\IQ\ltimes\IQ^n\rightarrow \IQ^\times$, so there are
  no non-trivial group morphisms to $\IQ^\times$ from
  $\ker\pi=\SL_{n+1}\IQ\free_{\ST_{n}\IQ\ltimes\IQ^n}(\SL_n\IQ\ltimes\IQ^n)$.
  We denoted $\ST_n\IQ=\T_n\IQ\cap\SL_n\IQ$.
  
  Let $\omega:\ker\pi\times X\rightarrow \cG$ be a cocycle to a
  $\Ufin$ target group.
  Observe that $\SL_{n+1}\IQ$ is a weakly rigid group, so by \cite[theorem
  0.1]{Popa:CocycleSuperrigidityMalleable} we
  can assume that $\omega\restrict{\SL_{n+1}\IQ}$ is a group
  morphism. In particular the restriction of $\omega$ to $\IQ^n$ does
  not depend on the $X$-variable. Proposition 3.6 in
  \cite{Popa:CocycleSuperrigidityMalleable} implies that
  $\omega\restrict{\SL_n\IQ\ltimes\IQ^n}$ is a group morphism. Because
  $\SL_{n+1}\IQ$ and $\SL_n\IQ\ltimes\IQ^n$ generate $\ker\pi$, we
  have that $\omega$ itself is a group morphism.
  
  It remains to prove the third condition. Let $\Delta$ be a
  self-conjugation of $\ker\pi\curvearrowright (X,\mu)$ with
  isomorphism $\delta:\ker\pi\rightarrow\ker\pi$. Observe that the
  rotations $A$ in $\SL_{n}\IQ$ by $90^\circ$, around any
  $n-2$-dimensional linear subspace of $\IQ^n$, can not be conjugated
  into $\T_n\IQ\ltimes\IQ^n$ inside $\GL_n\IQ\ltimes\IQ^n$. Hence
  they move every point $gH$ with
  $g\not\in \GL_n\IQ\ltimes\IQ^n$, while the points $gH$ with
  $g\in\GL_n\IQ\ltimes\IQ^n$ form a copy of the affine space
  $\IQ^n$. Therefore $H\cap\ker\pi$ acts with infinite orbits on
  $\Gamma/H\setminusb\{H\}$, so \cite[proposition
  6.10]{Vaes:Bimodules}
  implies that $\Delta$ is of
  the form $\Delta(x)_{\alpha(i)}=x_i$ where
  $\alpha:\Gamma/H\rightarrow \Gamma/H$ is a conjugation with
  isomorphism $\delta$.
  
  The isomorphism $\delta$ and the value $\alpha(H)$ completely
  determine $\alpha$. If $\alpha(H)=hH$, then $\delta(H\cap\ker\pi)=h
  (H\cap\ker\pi)h^{-1}$. Because $\Gamma_1=\SL_{n+1}\IQ$ and
  $\Gamma_2=\SL_n\IQ\ltimes\IQ^n$ are weakly rigid groups, $\delta$
  maps $\Gamma_1$, $\Gamma_2$ into conjugates of $\Gamma_i$ or
  $\Gamma_j$ with $i,j\in\{1,2\}$ (for an elementary proof, see for
  example point 4 of the proof of theorem 3.1 in
  \cite{DeprezVaes:Subfactors}). By symmetry, $\delta$ actually maps
  $\Gamma_1$ onto a conjugate $g_1\Gamma_1g_1^{-1}$ and $\Gamma_2$
  onto $g_2\Gamma_2 g_2^{-1}$. In particular, we find an automorphism
  $\delta_2=\Ad_{g_2^{-1}}\circ\delta\restrict{\Gamma_2}$ of $\Gamma_2$
  that maps $\SL_n\IQ$ onto a conjugate of itself. The
  same argument as in the proof of proposition \ref{prop:ex:eqrel} yields an
  $A\in\GL_n\IQ\ltimes\IQ^n$ such that $\delta_2=\Ad_A$. We can
  assume that $h=g_2A$. Now the
  isomorphism $\tilde\delta=\Ad_{A^{-1}g_2^{-1}}\circ \delta$ is
  identity on $\Gamma_2$, while it maps $\Gamma_1$
  onto a conjugate $\tilde g_1\Gamma_1\tilde g_1^{-1}$. But the only
  elements of $\Gamma$ that conjugate $\ST_n\IQ\ltimes\IQ^n$ into
  $\Gamma_1$ are the elements of $\Gamma_1$, so we may assume that
  $\tilde g_1=e$. Now $\tilde\delta\restrict{\Gamma_1}$ is an
  isomorphism of $\SL_{n+1}\IQ$ that is identity on
  $\ST_n\IQ\ltimes\IQ^n$. Such an isomorphism is easily seen to be
  the identity. We have proven that \[\alpha(gH)=g_2AgH\text{ for all
  }gH\in\ker\pi/(H\cap\ker\pi)=\Gamma/H,\]
  so we also have that $\Delta=g_2A\in\Gamma$.
\end{proof}

\section{A non-abelian generalization}
\label{sect:noncartan}
In this section we apply a kind of first quantization step to
theorem \ref{thm:intro:fundg:cartan}. We replace the measure space
$(Y,\nu)$ by type II$_\infty$ factor $B$, or more generally by a
properly infinite von Neumann algebra $B$ with a semifinite trace
$\Tr$. Then we consider a trace preserving action $\beta$ of a group
$\Lambda$ on $B$, such that $\Lambda$ acts ergodically on the center
of $B$. As above, we take a free, ergodic, p.m.p. action
$\alpha:\Gamma\actson (X,\mu)$ and a quotient
$\pi:\Gamma\rightarrow\Lambda$. We define a new action $\sigma$ of
$\Gamma$ on $\Lp^\infty(X,\mu)\vnOtimes B$ by the formula
$\sigma_g(a\otimes b)=\alpha_g(a)\otimes\beta_{\pi(g)}(b)$. Consider
the crossed product von Neumann algebra
$M=(\Lp^\infty(X)\vnOtimes B)\rtimes \Gamma$.

Every automorphism
$\psi_0$ of $B$ that commutes with the action of $\Lambda$ extends to an
automorphism $\psi$ of $M$ by $\psi((a\otimes
b)u_g)=(a\otimes\psi_0(b))u_g$. In fact, the requirement that $\psi_0$
really commutes with the action of $\Lambda$ is too strict. Let
$\psi_0$ be an automorphism of $B$. We say that $\psi_0$ commutes with
the action of $\Lambda$ up to a cocycle if there are unitaries
$(v_s)_s$ in $B$ such that
\begin{align*}
  \psi\circ\beta_s&=\Ad_{v_s}\circ\beta_s\circ\psi\quad\text{for all }s\in\Lambda\\
  \text{and}\quad v_{st}&=v_s\beta_s(v_t)\quad\text{for all }s,t\in\Lambda.
\end{align*}
If $\psi_0$ commutes with the action of $\Lambda$ up to the cocycle
$(v_s)_s$, then we can extend $\psi_0$ to an automorphism $\psi$ of
$M$ by the formula $\psi((a\otimes b)u_g)=(a\otimes
\psi_0(b)v_{\pi(g)})u_g$.

Denote by $\CCentr_B(\Lambda)$ the group of all pairs $(\psi_0,(v_s)_s)$
where $\psi_0$ is an automorphism of $B$ that commutes with the action
of $\Lambda$ up to the cocycle $(v_s)_s$. The group law in
$\CCentr_B(\Lambda)$ is defined by
\[(\psi_0,(v_s)_s)\cdot(\varphi_0,(w_s)_s)=(\psi_0\circ\varphi_0,(\psi_0(w_s)v_s)_s).\]
We call this group the cocycle centralizer of the action of $\Lambda$
on $B$.
Observe that the map $\theta:\CCentr_B(\Lambda)\rightarrow \Out(M)$
defined by
\[\theta(\psi_0,(v_s))((a\otimes b)u_g)=(a\otimes
\psi_0(b)v_{\pi(g)})u_g\]
is a group morphism. If there is a unitary $w\in B$ such that
$\psi_0=\Ad_w$ and such that $v_s=w\beta_s(w^\ast)$ for all
$s\in\Lambda$, then $\theta(\psi_0,(v_s))$ is trivial in
$\Out(M)$. Such trivial elements $(\Ad_w,(w\beta_s(w^\ast))_s)$ will
be called the inner elements of $\CCentr_B(\Lambda)$, and the group of
all inner elements of $\CCentr_B(\Lambda)$ is denoted by
$\ICCentr_B(\Lambda)$. This is a normal subgroup of
$\CCentr_B(\Lambda)$ and we denote the quotient by
$\OCCentr_B(\Lambda)=\CCentr_B(\Lambda)/\ICCentr_B(\Lambda)$.

Under strong conditions on the action $\Gamma\actson (X,\mu)$, theorem
\ref{thm:main:noncartan} shows that the group morphism
$\theta:\OCCentr_B(\Lambda)\rightarrow\Out(M)$ defined above is an
isomorphism.
\begin{theorem}
  \label{thm:main:noncartan}
  Let $\pi:\Gamma\rightarrow\Lambda$ be a quotient morphism and let
  $\alpha:\Gamma\actson (X,\mu)$ be a free, ergodic, p.m.p.\ action such that
  the restriction to $\ker\pi$ is still ergodic. Let $\beta:\Lambda\actson
  (B,\Tr)$ be a trace preserving action on a semifinite von Neumann
  algebra such that the restricted action of $\Lambda$ on $\Centre(B)$ is
  ergodic. Define an action $\sigma:\Gamma\actson
  \Lp^\infty(X)\vnOtimes B$ by $\sigma_g(a\otimes
  b)=\alpha_g(a)\otimes \beta_{\pi(g)}(b)$. Consider the crossed
  product von Neumann algebra $M=(\Lp^\infty(X)\vnOtimes B)\rtimes\Gamma$.

  Assume that
  \begin{itemize}
  \item the conditions of theorem \ref{thm:proof:main:eqrel} are satisfied.
  \item the conditions of theorem \ref{thm:proof:main:factor} are satisfied.
  \item $\ker\pi$ is a perfect group, i.e.\ there are no non-trivial
    group morphisms $\ker\pi\rightarrow\Circle^1$.
  \item $B$ is properly infinite and the action of $\Lambda\actson
    \Centre(B)$ does not preserve any finite measure equivalent to the
    spectral measure.
  \end{itemize}

  Then the outer automorphism group of $M$ is
  \[\Out(M)=\OCCentr_B(\Lambda).\]
  Every automorphism $\psi_0$ of $B$ that commutes with the action of
  $\Lambda$ up to a cocycle automatically scales the trace $\Tr$ by a
  constant $\Mod(\psi_0)$. For any projection $p\in M$ with finite
  trace, the fundamental group of the type II$_1$ factor $p M p$ is given by
  \[\Fundg(p M p)=\Mod(\CCentr_B(\Lambda)).\]
\end{theorem}
\begin{proof}

  We show that the group morphism
  $\theta:\OCCentr_B(\Lambda)\rightarrow\Out(M)$ defined above is in
  fact an isomorphism.
  We first show that $\theta$ is one-to-one. Suppose that
  $\theta(\psi_0,(v_s)_s)$ is inner in $M$. Then there is a unitary
  $u\in M$ such that
  \begin{align*}
    a&=uau^\ast\text{ for all }a\in \Lp^\infty(X)\\
    \psi_0(b)&=ubu^\ast\text{ for all }b\in B\\
    \text{ and }v_{\pi(g)}u_g&=uu_gu^\ast\text{ for all }g\in\Gamma.
  \end{align*}
  The first equation implies that $u$ commutes with $\Lp^\infty(X)$,
  and since $\Gamma$ acts freely on $X$, it follows that $u\in
  A\otimes B$. Because $\ker\pi$ acts ergodically on $X$, the third
  equation shows that in fact $u\in B$. Now it is clear that
  $\psi_0=\Ad_u$ and $v_s=u\beta_s(u^\ast)$ for all $s\in\Lambda$.

  To prove the surjectivity of $\theta$, let $\psi:M\rightarrow M$ be an
  automorphism of $M$. Write $A=\Lp^\infty(X)$.
  By theorem \ref{thm:proof:main:factor},
  we find a unitary $u\in M$ such that $u\psi(A\otimes
  B)u^{\ast}=A\otimes B$. From now on, we assume that $\psi(A\otimes
  B)=A\otimes B$.

  We can assume that $\psi(A\otimes B)=A\otimes B$.\\
  \emph{Step 1: there is a unitary $u\in M$ such that
    $u\psi(a)u^\ast=a$ for all $a\in A$.}\\
  Because $\psi(A\!\otimes\! B)=A\!\otimes\! B$, we also have that
  $\psi(A\!\otimes\!\Centre(B))=A\!\otimes\!\Centre(B)$. Identify $\Centre(B)$
  with $\Lp^\infty(Y,\nu)$ for some measure space $(Y,\nu)$. The
  action of $\Lambda$ on $B$ defines a non-singular action of
  $\Lambda$ on $Y$. The isomorphism $\psi:\Lp^\infty(X\times
  Y)\rightarrow\Lp^\infty(X\times Y)$ is of the form
  $\psi(f)=f\circ\Delta^{-1}$ for some orbit equivalence
  $\Delta:X\times Y\rightarrow X\times Y$ (\cite{Singer:AutomorphismsOfFiniteFactors}).

  By theorem \ref{thm:proof:main:eqrel}, we find an element $\varphi$
  in the full group of $\Rel(\Gamma\actson X\times Y)$ and an
  element $\Delta_0\in\Centr_{\Autns(Y)}(\Lambda)$ satisfying
  $\varphi(\Delta(x,y))=(x,\Delta_0(y))$ almost everywhere. The
  unitary $u=u_\varphi\in\Lp^\infty(X\times Y)\rtimes\Gamma\subset M$
  corresponding to $\varphi$ normalizes $A\otimes B$ and satisfies
  the relation $u\psi(a)u^\ast=a$ for all $a\in A$.

  From now on, we assume that $\psi(a)=a$ for all $a\in A$. The unitaries
  $\overline v_g=\psi(u_g)u_g^\ast$, for $g\in\Gamma$, commute with $A$. Since
  $\Gamma$ acts freely on $X$, it follows that $\overline v_g\in
  A\otimes B$ for all $g\in\Gamma$.

  \emph{Step 2: there is a unitary $u\in A\otimes B$ such that
    $u\psi(u_g)u^\ast u_g^\ast\in B$ for all $g\in\Gamma$.}\\
  We consider $\overline v_g$ as a function $\overline
  v_g:X\rightarrow\Unitary(B)$, and we define a cocycle
  $\omega:\ker\pi\times X\rightarrow \Unitary(B)$ by the formula
  $\omega(g,x)=\overline v_g(gx)$. We want to apply $\Ufin$-cocycle
  superrigidity to the cocycle $\omega$, but $\Unitary(B)$ is
  not necessarily a $\Ufin$ group.

  Let $p\in B$ be a projection with full central support in $B$ and
  such that $p B p$ is a finite von Neumann algebra. We find partial
  isometries $(w_n)_n$ in $B$ with $w_n^\ast w_n=p$ and such that
  $\sum_n w_nw_n^\ast=1$. Consider the projection $q=\psi(p)$ as a map
  $q:X\rightarrow B$. Consider the $\Tr$-preserving faithful normal
  semifinite extended
  center-valued trace $T:B^+\rightarrow \widehat{\Centre(B)^+}$.
  For every $g\in\ker\pi$, we know that $p$ commutes with $u_g$, so
  $q$ commutes with $\psi(u_g)$. It follows that
  \[q(gx)= (u_g^\ast qu_g)(x) = (u_g^\ast\psi(u_g)\, q\,
  \psi(u_g)^\ast u_g)(x) = \overline v_g(x)q(x)\overline v_g(x)^\ast\]
  for all $g\in\ker\pi$.

  Hence
  the map $x\mapsto T(q(x))$ is invariant under the action of
  $\ker\pi$ on $X$. So all the projections $q(x)$ are equivalent in
  $B$. We find a partial isometry $w\in A\otimes B$ with $ww^\ast=q$
  and such that $q_1=w^\ast w\in B$. Remember that $p$ is a finite
  projection with full central support, as a projection in $B$, or
  equivalently as a projection in $A\otimes B$. It follows that the
  same is true for $q_1$. Hence we find partial isometries
  $(\widetilde w_n)_n$ in $B$ with $\widetilde w_n^\ast \widetilde
  w_n=q_1$ and such that $\sum_n \widetilde w_n \widetilde
  w_n^\ast=1$.

  Consider the unitaries $\tilde v_g\!=\!w^\ast \psi(u_g)wu_g^\ast\!\in\! q_1
  (A\otimes B)q_1$, for all $g\in\ker\pi$. When we consider the
  $\tilde v_g$ as functions $\tilde v_g:X\rightarrow \Unitary(q_1Bq_1)$,
  they define a cocycle $\omega:\ker\pi\times X\rightarrow
  \Unitary(q_1Bq_1)$, by the formula $\omega(g,x)=\tilde
  v_g(gx)$. Because $\Unitary(q_1Bq_1)$ is a $\Ufin$ group, cocycle
  superrigidity yields a unitary $v\in \Unitary(A\otimes q_1Bq_1)$ such
  that $v\tilde v_g \sigma_g(v)^\ast\in B$ for all $g\in\ker\pi$.

  Define a unitary $u\in A\otimes B$ by
  \[u=\sum_n \widetilde w_n v w^\ast \psi(w_n^\ast)\in A\otimes B.\]
  For all $g\in\ker\pi$, we see that $v_g=u\psi(u_g)u^\ast u_g^\ast$
  is contained in $B$. We show that this is in fact true for all
  $g\in\Gamma$. Set $q_n=\widetilde w_n\widetilde w_n^\ast$ and
  observe that $q_n$ commutes with $v_g$ for all $g\in\ker\pi$.
  For any $g\in\Gamma$, we know that $v_g=u\psi(u_g)u^\ast u_g^\ast$
  is contained in $A\otimes B$, so we can consider $v_g$ as a function
  $v_g:X\rightarrow \Unitary(B)$. This function satisfies
  \[v_g(hx)=v_hv_g(x)v_{g^{-1}hg}^\ast\text{ for all }h\in\ker\pi.\]
  Since the $v_h$ with $h\in\ker\pi$ commute with the $q_n$, we see
  that the same relation holds for $q_nv_g$. This element $q_nv_g$ is
  contained in the polish space $\Lp^2(B,\widetilde \Tr)$ where
  $\widetilde \Tr$ is a faithful, normal, semifinite trace on $B$ such
  that $\widetilde \Tr(q_n)$ is finite. Weak mixing (see \cite[lemma 5.4]{PopaVaes:lattices})
  implies that $q_nv_g\in B$. This is true for all
  $n$, so $v_g\in B$ for every $g\in\Gamma$.
  We can assume that $v_g=\psi(u_g)u_g^\ast\in B$ for all
  $g\in\Gamma$.

  \emph{Step 3: $v_g$ only depends on $\pi(g)$ and there is an isomorphism
    $\psi_0\in\CCentr_{B}(\Lambda)$ such that, for all $a\in
    A$, $b\in B$ and $g\in\Gamma$, we have that
    $\psi((a\otimes b)u_g)\!=\!(a\otimes \psi_0(b))v_gu_g$.}\\
  For any $b\in B$ with $\Tr(b^\ast b)<\infty$, we see that
  $c=\psi(b)$ commutes with $A$ and so $c\in A\otimes B$. We can
  consider $c$ as a function $c:X\rightarrow \Lp^2(B,\Tr)$. Observe
  that, since $b$ commutes with all $u_g$ for $g\in\ker\pi$, we have
  \[c(gx)=v_gc(x)v_g^\ast\text{ almost everywhere and for all
  }g\in\ker\pi.\]
  Weak mixing (see \cite[lemma 5.4]{PopaVaes:lattices}) shows that $c$ is essentially
  constant, or still, that $c\in B$. So $\psi$ maps $B$ into $B$. By
  symmetry, it follows that $\psi(B)=B$. Write
  $\psi_0=\psi\restrict{B}$. Then we see that
  \[\psi((a\otimes b)u_g)=(a\otimes \psi_0(b))v_gu_g\text{ for all
  }a\in A, b\in B\text{ and }g\in\Gamma.\]

  We still have to show that $v_g$ only depends on $\pi(g)$. Since
  $(v_g)_g$ is a cocycle, it suffices to show that $v_g=1$ for all
  $g\in\ker\pi$. Observe that $v_g$ is in the center of $B$ whenever
  $g\in\ker\pi$, and the application $\ker\pi\ni g\mapsto
  v_g\in\Unitary(\Centre(B))$ is a group morphism. Since $\ker\pi$ is
  a perfect group, this morphism is trivial. We have shown that
  $v_g=1$ for every $g\in\ker\pi$.
\end{proof}

\section{A flexible class of examples}
\label{sect:ex:general}
We have given one explicit example of an action $\Gamma\actson
(X,\mu)$ and a quotient $\pi:\Gamma\rightarrow\Lambda$ that satisfy
the conditions of theorem \ref{thm:proof:main:factor}, see proposition
\ref{prop:ex:factor}. In this example, the quotient group $\Lambda$
was abelian. For the applications in section \ref{sect:applications},
we need more flexibility in the choice of $\Lambda$.

We construct a new class of examples as follows. Let $\Gamma_1$ be a
countable group and $\Sigma\subset \Gamma_1$ a subgroup. For any
countable group $\Lambda$, we can consider
$\Gamma=\Gamma_1\free_\Sigma (\Sigma \times\Lambda)$, together with
the obvious quotient morphism $\pi:\Gamma\rightarrow\Lambda$. Let
$H\subset \Gamma$ be a subgroup. We consider the generalized Bernoulli
action $\Gamma\actson (X,\mu)=(X_0,\mu_0)^{\Gamma/H}$, over an atomic
base space $(X_0,\mu_0)$ with unequal weights. Among the conditions we
have to check, the least standard one says that
$\Norm_{\Autmp(X,\mu)}(\ker\pi)=\Gamma$.

If we have that $\Stab_{\ker\pi}\{i\}j$ is infinite for all $i\not=j\in
\Gamma/H$, then we know that the group of conjugations of
$\ker\pi\actson X$ is given by
$\Norm_{\Autmp(X,\mu)}(\ker\pi)=\Norm_{\Perm(\Gamma/H)}(\ker\pi)$ (see
\cite[proposition 6.10]{Vaes:Bimodules}).

Observe that $\ker\pi$ is the infinite amalgamated free product of
copies of $\Gamma_1$, amalgamated over $\Sigma$. The copies of
$\Gamma_1$ correspond to the conjugates
$\lambda\Gamma_1\lambda^{-1}\subset\Gamma$, for $\lambda\in\Lambda$.

If we would choose $H$ to be a subgroup of $\Gamma_1$, then any
permutation $\eta:\Lambda\rightarrow\Lambda$ defines an isomorphism
$\delta_\eta:\ker\pi\rightarrow\ker\pi$ by the formula
$\delta_\eta(\lambda
g\lambda^{-1})=\eta(\lambda)g\eta(\lambda)^{-1}$. for all $g\in\Gamma$
and $\lambda\in\Lambda$. The formula $\alpha_\eta(g\lambda
H)=\delta_\eta(g)\eta(\lambda)H$ defines a conjugation
$\alpha_\eta\in\Norm_{\Perm(\Gamma/H)}(\ker\pi)$. The permutation
group of $\Lambda$ is uncountable, so it is certainly strictly larger
that $\Gamma$.

We choose $H$ in the following way. For every $\lambda\in\Lambda$, we
choose a ``sufficiently different'' subgroup
$H_\lambda\subset\Gamma_1$, and we define $H$ to be the subgroup of
$\Gamma$ that is generated by the $\lambda H_\lambda
\lambda^{-1}$. More precisely, we have the following theorem.
\begin{theorem}
  \label{thm:ex:general}
  Let $\Lambda$ be any countable group.
  Let $\Gamma_1$ be a property (T) group. Assume that $\Gamma_1$ is a
  perfect group and that there are no non-trivial group morphisms
  $\theta:\Gamma_1\rightarrow\Lambda$.

  Let $\Sigma\subset\Gamma_1$ be an amenable subgroup such that there
  is a finite set $g_1,\ldots,g_n\in\Gamma_1$ for which $\bigcap_i
  g_i\Sigma g_i^{-1}$ is finite. Assume that the elements of $\Sigma$
  are the only elements $g\in\Gamma_1$ for which $g\Sigma
  g^{-1}\cap\Sigma$ has finite index in $\Sigma$.

  For every $\lambda\in\Lambda$, choose
  an infinite subgroup $H_\lambda\subset\Gamma_1$, subject to the following
  conditions.
  \begin{itemize}
  \item If $\theta:\Gamma_1\rightarrow\Gamma_1$ is an automorphism
    that maps a finite index subgroup of $H_\lambda$ into $H_\mu$,
    then it follows that $\lambda=\mu$ and $\theta=\Ad_h$ for some
    $h\in H_\lambda$.
  \item If $\theta$ is an automorphism of $\Gamma_1$, then the
    intersection $\theta(H_\lambda)\cap\Sigma$ is trivial, for all
    $\lambda\in\Lambda$.
  \end{itemize}

  Consider the group
  $\Gamma=\Gamma_1\free_\Sigma(\Sigma\times\Lambda)$ with its natural
  quotient $\pi:\Gamma\rightarrow\Lambda$. Define $H$ to be the
  subgroup of $\Gamma$ that is generated by the $\lambda
  H_\lambda\lambda^{-1}$. Consider the generalized Bernoulli action
  $\Gamma\actson (X,\mu)=(X_0,\mu_0)^{\Gamma/H}$ over an atomic base
  space $(X_0,\mu_0)$ with unequal weights. Then the action
  $\Gamma\actson (X,\mu)$ satisfies the conditions of theorem
  \ref{thm:proof:main:factor}.
\end{theorem}
\begin{proof}
  The action $\Gamma\actson X$ was constructed so that it
  satisfies the conditions of theorem
  \ref{thm:intro:main:factor}. Because $\ker\pi$ is generated by
  copies of the perfect group $\Gamma_1$, it is itself a perfect
  group. By the same argument, there are no non-trivial group
  morphisms $\theta:\ker\pi\rightarrow\Lambda$. It remains to show
  that $\ker\pi$ acts cocycle superrigidly on $X$ and that
  $\Norm_{\Autmp(X)}(\ker\pi)=\Gamma$. To prove this second condition,
  it is sufficient to show two properties of the action
  $\ker\pi\actson \Gamma/H$. We show that $H$ acts with
  infinite orbits on $\Gamma/H\setminusG \{H\}$, and we show that
  $\Norm_{\Perm(\Gamma/H)}(\ker\pi)=\Gamma$.

  First we prove a general property of the subgroups $H_\lambda$ and
  $\Sigma$.\\
  \emph{Claim: Let $\delta:\Gamma_1\rightarrow\Gamma_1$ be a group
    automorphism and let $g$ be an element of $\Gamma$. Then we have
    the following.
    \begin{itemize}
    \item If $G\subset H_\lambda$ is a finite index subgroup such that
      $g\theta(G)g^{-1}$ is contained in $H$, then $\theta$ is of the
      form $\theta=\Ad_h$ for some $h\in \Gamma_1$ with
      $gh\lambda^{-1}\in H$.
    \item $\Sigma\cap g^{-1}Hg$ has infinite index in $\Sigma$.
    \item $\Sigma\cap H_\lambda=\{e\}$ for all $\lambda\in\Lambda$.
    \end{itemize}
  }
  Suppose that $g\subset H_\lambda$ is a finite index subgroup such
  that $g\theta(G)g^{-1}$ is contained in $H$. Possibly replacing
  $\theta$ by $\Ad_{g_0}\circ\theta$ for some $g_0\in\Gamma_1$, we can
  assume that $g$ is an
  element of minimal length in $Hg\Gamma_1$. So we can write
  $g=\lambda_0g_1\ldots g_n\lambda_n$ with $\lambda_0\in\Lambda$,
  $\lambda_1,\ldots,\lambda_n\in\Lambda\setminusG\{e\}$, with
  $g_1,\ldots,g_n\in\Gamma_1\setminusG\Sigma$ and such that
  $g_1\in\Gamma1\setminusG H_{\lambda_0}\Sigma$. The assumptions of
  the theorem yield an element $k\in G$ with
  $\theta(k)\not\in\Sigma$. It follows that the expression
  \[\lambda_0g_1\ldots
  g_n\lambda_n\theta(k)\lambda_n^{-1}g_n^{-1}\ldots g_1^{-1}\lambda_0^{-1}\]
  is a reduced expression for $g\theta(k)g^{-1}\in H$. Since $g_1$ was
  not contained in $H_{\lambda_0}\Sigma$, this is only possible if $n$
  was $0$. In that case, we have that $\theta(G)\subset
  H_{\lambda_0}$. It follows that $\lambda_0=\lambda$ and $\theta=\Ad_h$
  for some $h\in H_\lambda$. Hence $gh\lambda^{-1}=\lambda
  h\lambda^{-1}\in H$. In this process we replaced $g$ by an element
  $h_0gg_0$ in $Hg\Gamma_1$ and $\theta$ by $\Ad_{g_0^{-1}}$. This
  does not affect the conclusion.

  Suppose that $G$ was a finite index subgroup of $\Sigma$ such that
  $gGg^{-1}\subset H$. We can assume that $g$ has minimal length among
  the elements of $Hg\Lambda$. So we can write
  $g=\lambda_1g_1\ldots\lambda_ng_n$ with $\lambda_1\in\Lambda$,
  $\lambda_2,\ldots,\lambda_n\in\Lambda\setminusG\{e\}$, with
  $g_1,\ldots,g_n\in\Gamma_1\setminusG\Sigma$ and such that
  $g_1\not\in H_{\lambda_1}\Sigma$. By the conditions on $\Sigma$, we
  find $k\in G$ such that $g_nkg_n^{-1}\not\in\Sigma$. So we see that
  the expression
  \[\lambda_1g_1\ldots \lambda_n (g_nkg_n^{-1}) \lambda_n^{-1}\ldots
  g_1^{-1}\lambda_1\]
  is a reduced expression for an element of $H$. This is impossible:
  if $n\not=1$, then we know that $\lambda_1g_1$ can never be the
  first letters of an element of $H$. If $n=1$, then it would follow
  that $g_1kg_1^{-1}\in H_{\lambda_1}$, which is impossible because
  $e\not=k\in\Sigma$. This finishes the proof of the claim.

  \emph{Step 1: The action $\ker\pi\actson (X,\mu)$ is cocycle
    superrigid.}
  Let $\omega:\ker\pi\times X\rightarrow \cG$ be a cocycle with a
  $\Ufin$ target group. For every $\lambda\in\Lambda$, we know that
  the restriction $\omega\restrict{\lambda\Gamma_1\lambda^{-1}}$ is
  cohomologous to a group morphism, by Popa's cocycle superrigidity
  theorem \cite[theorem 0.1]{Popa:CocycleSuperrigidityMalleable}. So there are maps
  $\varphi_\lambda:X\rightarrow\cG$ and group morphisms
  $\theta_\lambda:\lambda\Gamma_1\lambda^{-1}\rightarrow\cG$ such that
  \[\omega(g,x)=\varphi_\lambda(gx)^{-1}\theta_\lambda(g)\varphi_\lambda(x)\]
  for all $\lambda\in\Lambda$ and $g\in\lambda\Gamma_1\lambda^{-1}$,
  and for almost all $x\in X$. These descriptions must match for
  $g\in\Sigma$, so we find that
  \[(\varphi_\lambda\varphi_\mu^{-1})(hx)=\theta_\lambda(h)(\varphi_\lambda\varphi_\mu)(x)\theta_\mu(h)^{-1}\]
  for all $h\in\Sigma$, $\lambda,\mu\in\Lambda$ and for almost all
  $x\in X$. The second assertion in our claim shows that $\Sigma$ acts
  weakly mixingly on $X$, so it follows that
  $\varphi_\lambda\varphi_\mu^{-1}$ is essentially constant (see for example
  \cite[lemma 5.4]{PopaVaes:lattices}). Hence we can assume that all
  the $\varphi_\lambda$ are the same. Because the
  $\lambda\Gamma_1\lambda^{-1}$ generate $\ker\pi$, we see that
  $\omega$ is cohomologous to a group morphism.

  \emph{Step 2: The group $H$ acts with infinite orbits on
    $\Gamma/H\setminusG\{H\}$.}
  We have to show that $gHg^{-1}\cap H$ has infinite index in $H$
  whenever $g\in\Gamma\setminusG H$. This follows immediately from the
  claim above.

  \emph{Step 3: We show that
    $\Norm_{\Perm(\Gamma/H)}(\ker\pi)=\Gamma$.}
  Let $\alpha:\Gamma/H\rightarrow\Gamma/H$ be a conjugation for the
  action of $\ker\pi$. Possibly replacing $\alpha$ by $i\mapsto
  g\alpha(i)$ for some $g\in\Gamma$, we can assume that
  $\alpha(H)=H$. Denote by $\delta:\ker\pi\rightarrow\ker\pi$ the
  group automorphism such that $\alpha(gi)=\delta(g)\alpha(i)$ for all
  $g\in\ker\pi$ and $i\in\Gamma/H$. Remark that $\Delta(H)=H$.

  Fix $\lambda\in\Lambda$. Because $\lambda\Gamma_1\lambda^{-1}$ has
  property (T), we find an element $g\in\Gamma$ such that
  $\delta(\lambda\Gamma_1\lambda^{-1})\subset g\Gamma_1 g^{-1}$ (see
  for example \cite[part 4 of the proof of theorem 3.1]{DeprezVaes:Subfactors} for an elementary proof). By
  symmetry we actually have that
  $\delta(\lambda\Gamma_1\lambda^{-1})=g\Gamma_1g^{-1}$. Define an
  automorphism $\theta:\Gamma_1\rightarrow\Gamma_1$ by
  $\theta(k)=g^{-1}\delta(\lambda k \lambda^{-1})g$. Then it follows
  that $g\theta(H_\lambda)g^{-1}\subset H$. Our claim above shows that
  $\theta$ is of the form $\theta=\Ad_h$ and that
  $h_\lambda=gh\lambda^{-1}\in H$.

  For every $\lambda\in\Lambda$, we have found an element
  $h_\lambda\in H$ such that $\delta(\lambda
  k\lambda^{-1})=h_\lambda\lambda k\lambda^{-1}h_{\lambda}^{-1}$ for
  all $k\in\Gamma_1$. For $k\in\Sigma$, these different descriptions
  have to match, so $h_\mu^{-1}h_\lambda$ commutes with $\Sigma$ for
  all $\lambda,\mu\in\Lambda$. The only element in $H$ that commutes
  with $\Sigma$ is $e$, so we have found one element $h\in H$ such
  that $\delta(k)=hkh^{-1}$ for all $k\in\ker\pi$. Replacing
  $\alpha(i)$ by $h^{-1}\alpha(i)$, we can assume that $\alpha$
  commutes with the action of $\ker\pi$. We show that then
  $\alpha=\id$. Let $kH$ be an element in $\Gamma/H$, then we know
  that $kHk^{-1}$ fixes $\alpha(kH)$. But step 2 shows that $kHk^{-1}$
  fixes only the point $kH$ itself. It follows that
  $\alpha(kH)=kH$. This works for any $k\in\Gamma$, so
  $\alpha=\id$. In the course of this argument, we changed $\alpha(i)$
  to $g^{-1}\alpha(i)$ for some $g\in\Gamma$, so we have shown that
  $\alpha(i)=gi$ for all $i\in\Gamma/H$.
\end{proof}

\begin{construction}
  \label{constr:ex:general}
  Let $R=\F_2[X]$ be the ring of polynomials over the field $\F_2$ of
  two elements and take a natural number $k\geq 1$. Consider the group
  $\Gamma_1=\SL_{3k}R\ltimes R^3$.

  Let $\Lambda$ be any countable group for which there are no non-trivial
  group morphisms $\theta:\Gamma_1\rightarrow\Lambda$. We will define
  the groups $\Sigma$ and $H_\lambda$ with
  $\lambda\in\Lambda$. Consider one symbol $\ast\not\in\Lambda$
  representing $\Sigma$. Take a
  one-to-one map $(\Lambda\sqcup\{\ast\})\times\{1,\ldots,k\}\ni(\lambda,i)\mapsto n_{\lambda,i}\in\IN$ and
  consider the matrices
  \[h_\lambda=\left(\begin{matrix}h_{\lambda,1}&&\\&\ddots&\\&&h_{\lambda,k}\end{matrix}\right)
  \quad\text{ where }\qquad
  h_{\lambda,i}=\left(\begin{matrix}0&X^{2n_{\lambda,i}+1}&1\\1&0&0\\0&1&0\end{matrix}\right).\]
  Define subgroups $H_\lambda\subset\Gamma_1$ by the relation
  $H_\lambda=\{g\in\Gamma_1\mid g(h_\lambda,0)=(h_\lambda,0) g\}$. Set $\Sigma=H_{\ast}$.

  As in theorem \ref{thm:ex:general}, consider
  $\Gamma=\Gamma_1\free_\Sigma(\Sigma\times\Lambda)$, with its natural
  quotient morphism $\pi:\Gamma\rightarrow\Lambda$. Let $H$ be the
  subgroup generated by the $\lambda H_\lambda\lambda^{-1}$. Consider
  the generalized Bernoulli action $\Gamma\actson
  (X,\mu)=(X_0,\mu)^{\Gamma/H}$ over an atomic base space
  $(X_0,\mu_0)$ with unequal weights.
\end{construction}
\begin{proposition}
  \label{prop:ex:general}
  The action $\Gamma\actson (X,\mu)$ and quotient $\pi:\Gamma\rightarrow \Lambda$
  as in construction \ref{constr:ex:general}, satisfy the conditions of
  theorem \ref{thm:intro:main:factor}.
\end{proposition}
\begin{proof}
  We check that $\Gamma_1$, $\Sigma$ and the $H_\lambda$ satisfy the
  conditions of theorem \ref{thm:ex:general}.

  The group $\Gamma_1$ is a property (T) group, see for example
  \cite[example 3.4.1]{BekkaDelaharpeValette}. This group $\Gamma_1$
  is a perfect group because $\Gamma_1$ is generated by elements of
  the form $x(1+a e_{i,j},0)x^{-1}$ with $x\in\Gamma_1$, $a\in R$ and
  $i\not=j$. The symbol $e_{i,j}$ denotes the matrix whose $i,j$-th
  component is $1$ while all other components are $0$. The element
  $x(1+a e_{i,j},0)x^{-1}$ is the commutator of $x(1+ae_{k,j})x^{-1}$
  with $x(1+e_{i,k})x^{-1}$ where $k\not=i,j$.

  Consider the field $E=\F_2(X)$ of rational functions over $\F_2$,
  and denote its algebraic closure by $\overline{E}$.
  Fix $\lambda\in\Lambda\sqcup\{\ast\}$ and $1\leq i\leq k$. Observe
  that the characteristic polynomial of $h_{\lambda,i}$ is equal to
  $Y^3+X^{2n_{\lambda,i}+1}Y+1$. Denote its three roots in
  $\overline{E}$ by
  $s_{\lambda,i,1}, s_{\lambda,i,2}$ and $s_{\lambda,i,3}$. Consider
  the extension fields $E_{\lambda,i,j}=E(s_{\lambda,i,j})$ and remark
  that lemma \ref{lem:ex:general} shows that $E_{\lambda,i,j}\cap
  E_{\lambda^\prime,i^\prime,j^\prime}=E$ whenever
  $(\lambda,i,j)\not=(\lambda^\prime,i^\prime,j^\prime)$.

  Consider the ring automorphism $\psi:R\rightarrow R$ that is defined
  by $\psi(X)=X+1$. Remark that the only outer automorphism of
  $\Gamma_1$ is induced by this ring automorphism $\psi$. We still
  denote the group automorphism by $\psi$, and we extend the ring
  automorphism $\psi:R\rightarrow R$ to $\overline{E}$. This
  extension is unique up to multiplication by an element of the Galois
  group of $\overline{E}$ over $E$. Again it follows from lemma
  \ref{lem:ex:general} that $E_{\lambda,i,j}\cap
  \psi(E_{\lambda^\prime,i^\prime,j^\prime})=E$ for any
  $\lambda,\lambda^\prime\in \Lambda\sqcup\{\ast\}$, $1\leq
  i,i^\prime\leq k$ and $1\leq j,j^\prime\leq 3$.

  The eigenvector $\xi_{\lambda,i,j}$ of $h_{\lambda,i}$ corresponding
  to the eigenvalue $s_{\lambda,i,j}$ is contained in
  $E_{\lambda,i,j}^3$ but not in $E^3$. Denote by
  $H_{\lambda,i}\subset\SL_3R$ the subgroup of all matrices that
  commute with $h_{\lambda,i}$. Observe that this is exactly the
  group of all matrices whose eigenvectors are $\xi_{\lambda,i,1},
  \xi_{\lambda,i,2}$ and $\xi_{\lambda,i,3}$.

  Let $g\in H_{\lambda,i}$
  be a nontrivial element. The $g$ has three distinct eigenvalues
  because if the eigenvalues corresponding to $\xi_{\lambda,i,1}$ and
  $\xi_{\lambda,i,2}$ were equal, then this eigenvalue must be
  contained in $E=E_{\lambda,i,1}\cap E_{\lambda,i,2}$. And hence so
  must the third eigenvalue. Now all three eigenvalues are invertible
  elements of $R$. The only such element is $1$, so $g$ would be $1$.
  Suppose that $k\in \MatM_3(R)$ is a possibly non-invertible matrix
  such that $gk\in kH_{\lambda^\prime,i^\prime}$ for some
  $\lambda^\prime$ and $i^\prime$. Then
  $k\xi_{\lambda^\prime,i^\prime,j}$ is an eigenvector of $g$ and is
  contained in $E_{\lambda^\prime,i^\prime,j}$. Hence we see that
  either $k\xi_{\lambda^\prime,i^\prime,j}=0$ or it follows that
  $\lambda^\prime=\lambda$, $i^\prime=i$ and $\xi_{\lambda,i,j}$ is an
  eigenvector of $k$. Because this holds for all $j=1,2,3$, we
  conclude that $k$ commutes with $h_{\lambda,i}$.
  Similarly, if we have that $gk\in
  k\psi(H_{\lambda^\prime,i^\prime})$, then it follows that $k=0$.

  We conclude that for every automorphism
  $\theta:\Gamma_1\rightarrow\Gamma_1$ and for every
  $\lambda\not=\lambda^\prime\in \Lambda\sqcup\{\ast\}$, the intersection
  $H_\lambda\cap \theta(H_{\lambda^\prime})$ is trivial. Moreover, if
  $\theta$ maps a finite index subgroup of $H_\lambda$ into itself,
  then we have that $\theta=\Ad_{h}$ for some $h\in H_\lambda$.

  Observe that $\Sigma$ is abelian, and hence amenable. Fix elements
  $g_i\in\Gamma_1\setminusG H_{\ast,i}$ and consider the diagonal
  matrix $g=\diag(g_1,\ldots,g_k)$. Then it is clear that $g\Sigma
  g^{-1}\cap\Sigma$ is trivial and therefore finite. We have checked
  all conditions of theorem \ref{thm:ex:general}.
\end{proof}

\begin{lemma}
  \label{lem:ex:general}
  Denote by $E=\F_2(X)$ the field of rational functions over
  $\F_2$. Let $a\in \F_2[X]\subset E$ be either $a=X^{2k+1}$ or
  $a=(X+1)^{2k+1}$ for some $k\in\IN$. Consider the polynomial $P_a(Y)=Y^3+aY+1$, and
  denote the splitting field of $P_a$ over $E$ by $E_a$. Then we have
  that
  \begin{itemize}
  \item $E_a$ is a degree 6 extension of $E$
  \item $E_a\cap E_b=E$ whenever $b\not=a$ is of the form $X^{2l+1}$
    or $(X+1)^{2l+1}$ for some $l\in\IN$.
  \end{itemize}
\end{lemma}
\begin{proof}
  Remark that $P_a(Y)$ is an irreducible polynomial for if $n,d\in
  \F_2[X]$ were coprime elements such that $P_a(n/d)=0$, then
  it follows that $n=d=1$, and hence that $a=0$. So the degree
  $[E_a:E]$ is at least $3$.

  Denote the roots of $P_a(Y)$ in the algebraic closure of $E$ by
  $s_1,s_2,s_3\in E_a$. Consider the element
  $t=s_1s_2^2+s_2s_3^2+s_3s_1^2\in E_a$. Observe that this element
  satisfies the second degree equation $Q_a(t)=t^2+t+1+a^3=0$. Suppose
  $t$ were an element of $E$. Since $1+a^3$ is a polynomial, the same
  is true for $t$. Then we know that $t^2+t$
  has even degree. This contradicts the fact that $1+a^3$ has odd
  degree. So we see that
  $E\subsetneqq E(t)\subsetneqq E_a$, where $E(t)$ is a degree 2
  extension. Hence $E_a$ is a degree 6 extension of $E$. Denote this
  intermediate field by $K_a=E(t)$.

  Let $a\not=b$ be either $b=X^{2l+1}$ or $(X+1)^{2l+1}$. The intersection
  of fields $E_a\cap E_b$ is a normal subfield of $E_a$. So it is
  either $E_a$ itself, $K_a$ or $E$. We have to show that only the
  last case can occur. In both of the former cases, we see that
  $K_a=K_b$, because $K_a$ is the unique subfield of $E_a$ that has
  degree 2 over $E$. So we find elements $x, y\in E$ such that $xt+y$
  satisfies $Q_b(xt+y)=0$, while $Q_a(t)=0$. This implies that
  $y^2+y=a^3+b^3$. Since the right hand side is a polynomial over
  $\F_2$, the same is true for $y$. Hence $y^2+y$ is a polynomial of
  even degree. If $k\not=l$, then
  we see that $a^3+b^3$ is a polynomial of odd degree. So we can assume that
  $a=X^{2k+1}$ and $b=(X+1)^{2k+1}$. But now $X$ divides $y^2+y$ but
  it does not divide $a^3+b^3$. We have shown that $E_a\cap E_b=E$.
\end{proof}

\def\niets{
\begin{construction}
  \label{constr:ex:general}
  Let $R=\F_2[X]$ be the ring of polynomials over the field $\F_2$ of
  two elements. Consider the group $\Gamma_1=\SL_3R\ltimes R^3$. Set
  $\Sigma=\T_3R\ltimes Re_1$ where $T_3R$ denotes the group of upper
  triangular matrices in $\SL_3R$ and where $e_1$ is the first
  standard basis vector in $R^3$. Let $\Lambda$ be any countable group
  such that there are no non-trivial group morphisms
  $\theta:\Gamma_1\rightarrow\Lambda$. Take a map
  $\lambda\mapsto P_\lambda(X)$ that associates to every element $\lambda
  \in\Lambda$ a non-trivial irreducible polynomial over $\F_2$. Choose
  this map such that $P_\lambda(X)\not=P_\mu(X)$ whenever $\lambda\not=\mu$
  and such that $P_\lambda(X)\not=P_\mu(X+1)$ for any
  $\lambda,\mu\in\Lambda$. Set
  \[h_\lambda=\left(
    \left(\begin{matrix}P_\lambda(X)+X&1&0\\XP_\lambda(X)+X^2+1&X&0\\0&0&1\end{matrix}\right),
    \left(\begin{matrix}0\\0\\0\end{matrix}\right)
  \right)\in\Gamma_1.\]
  Define the group $H_\lambda=\{h\in\Gamma_1\mid hh_\lambda
  h^{-1}=h_\lambda\text{ or }h_\lambda^{-1}\}$. Set
  $\Gamma=\Gamma_1\free_\Sigma(\Sigma\times\Lambda)$ and consider the
  canonical quotient $\pi:\Gamma\rightarrow\Lambda$. Let $H$ be the
  subgroup generated by the $\lambda H_\lambda\lambda^{-1}$. Consider
  the generalized Bernoulli action $\Gamma\actson
  (X,\mu)=(X_0,\mu)^{\Gamma/H}$ over an atomic base space
  $(X_0,\mu_0)$ with unequal weights.
\end{construction}
\begin{proposition}
  \label{prop:ex:general}
  The action $\Gamma\actson (X,\mu)$ and quotient $\pi:\Gamma\rightarrow \Lambda$
  as in construction \ref{constr:ex:general}, satisfy the conditions of
  theorem \ref{thm:intro:main:factor}.
\end{proposition}
\begin{proof}
  We only have to check the conditions of theorem
  \ref{thm:ex:general}. It is well-known that $\Gamma_1$ is a property
  (T) group and a perfect group. The subgroup $\Sigma$ is
  amenable and the only elements $g\in\Gamma_1$ for which $g\Sigma
  g^{-1}\cap\Sigma$ has finite index in $\Sigma$ are the elements of
  $\Sigma$. Moreover, we see that
  \[\left(\begin{matrix}0&0&1\\0&1&0\\1&0&0\end{matrix}\right)
  \Sigma
  \left(\begin{matrix}0&0&1\\0&1&0\\1&0&0\end{matrix}\right)\cap\Sigma=\{1\}.\]

  We still have to show that the subgroups $H_\lambda$ satisfy the
  three conditions of theorem \ref{thm:ex:general}. We prove this by
  studying the eigenvalues of $h_\lambda$. Denote by $F=\F_2(X)$ the
  field of rational functions over $\F_2$. Observe that the
  characteristic polynomial of $h_\lambda$ is $(Y-1)(Y^2+P_\lambda(X)
  Y+1)$. Take a root $Y=s_\lambda$ of $Y^2+P_\lambda(X)Y+1$ in the algebraic
  closure of $F$. The eigenvalues of $h_\lambda$ are $1$, $s_\lambda$
  and $s_\lambda^{-1}$, and the corresponding eigenvectors are
  $e_3=(0,0,1)$, $\xi_\lambda=(1,X+s_\lambda^{-1},0)$ and
  $\eta_\lambda=(1,X+s_\lambda,0)$. From this it is easy to see that
  all elements of $H_\lambda$ are of the form $(A^iB,Q(X)e_3)$ where
  $i=0,1$, $Q(X)\in R$, where
  \[A=\left(\begin{matrix}X&1&0\\X^1+1&X&0\\0&0&1\end{matrix}\right)\]
  and where $B$ is of the form
  \begin{equation}
    \label{eqn:formB}
    B=\left(\begin{matrix}S(X)+P(X)T(X)&T(X)&0\\(X\,P(X)+X^2+1)T(X)&S(X)&0\\0&0&1\end{matrix}\right)
  \end{equation}
  for some $S(X),T(X)\in R$. In any case, the square of every element
  of $H_\lambda$ is of the form $(B,0)$ with $B$ as in (\ref{eqn:formB}).
  Moreover, we see that $H_\lambda\cap\Sigma=\{(1,0)\}$ for all
  $\lambda\in\Lambda$.

  The automorphism group of $\Gamma_1$
  is just $\AutRing(R)\ltimes\Gamma_1$, where $\AutRing(R)$ is the group of
  ring automorphisms of $R$. The only non-trivial ring automorphism of
  $R=\F_2[X]$ is given by $\psi(X)=X+1$. For every
  $\lambda\in\Lambda$, fix a root $\psi(s_\lambda)$ in the
  algebraic closure of $F$, of the polynomial $Y^2+P_\lambda(X+1)Y+1$.

  Fix $\lambda\in\Lambda$ and let $\theta=\psi^i\circ\Ad_g\in\Autg(\Gamma_1)$ be an
  automorphism. Observe that $h_\lambda$ has infinite order in
  $\Gamma_1$ and that $h_\lambda^k$ ($k\not=0$) has three different
  eigenvalues in $F(s_\lambda)$, and two of them are not contained in
  $F$. Hence $\theta(h_\lambda)^k$ has two eigenvalues in
  $F(\psi^i(s_\lambda))\setminusG F$

  If $\theta$ would map a finite index subgroup of $H_\lambda$
  into $\Sigma$, then $\theta$ would in particular map some
  $h_\lambda^k$, for $k\not=0$, into $\Sigma$. But the elements of $\Sigma$ all have
  $1$ as their only eigenvalue. This contradicts our observation about
  the eigenvalues of $\theta(h_\lambda)^k$.

  If $\theta$ maps a finite index subgroup of $H_\lambda$ into some
  $H_\mu$ for $\mu\in\Lambda$, then we know that
  $\theta(h_\lambda)^{k}\in H_\mu$. Hence its square
  $\theta(h_{\lambda})^{2k}$ has eigenvectors $\xi_\mu,\eta_\mu,e_3\in
  F(s_\mu)$. So the eigenvalues of $\theta(h_\lambda)^{2k}$ are
  elements of $F(s_\mu)$. But we already know that two of them are in
  $F(\psi^i(s_\lambda))\setminusG F$, and lemma \ref{lem:F2:indep}
  shows that the intersection is $F(\psi^i(s_\lambda))\cap
  F(s_\mu)=F$, unless $\mu=\lambda$ and $i=0$. In this case, we see
  that $\theta(h_\lambda)^{2k}=h_{\lambda}^{\pm2k}$ and hence $g\in H_\lambda$.
\end{proof}

\begin{lemma}
  \label{lem:F2:indep}
  If $Q_1,\ldots,Q_n$ are distinct irreducible polynomials in
  $\F_2[X]$, and if $s_1,\ldots,s_n$ are elements of the algebraic
  closure of $F=\F_2(X)$ satisfying $s_i^2+Q_is_i+1=0$, then the set
  $\{1,s_1,\ldots,s_n\}$ is linearly independent over $F$.
\end{lemma}
\begin{proof}
  Assume that $a_0,\ldots,a_n\in F$ are such that
  $a_0+a_1s_1+\ldots+a_ns_n=0$. By induction on $n$ we can assume that
  $1,s_1,\ldots,s_{n-1}$ are linearly independent, and that the
  $a_1,\ldots,a_n$ are non-zero. Squaring the relation above, we find that
  \[0=a_0^2+\sum_{i=1}^na_i^2s_i^2
  =a_0^2+\sum_{i=1}^na_0^2+\sum_{i=1}^na_i^2Q_is_i
  =\sum_{i=0}^na_i^2+\sum_{i=1}^{n-1}a_i^2Q_is_i + a_0a_nQ_n + \sum_{i=1}^{n-1}a_ia_nQ_ns_i.
  \]
  Since we assumed that $1,s_1,\ldots,s_{n-1}$ is linearly
  independent, it follows that
  \begin{align}
    \label{eqn:F2:1}
    a_nQ_n&=a_iQ_i\text{ for all }i=1,\ldots,n\\
    \label{eqn:F2:2}
    a_0a_nQ_n&=\left(\sum_{i=0}^na_i\right)^2
  \end{align}.
  Because all these relations are homogeneous in the $a_i$, we can as
  well assume that all the $a_i$ are polynomials in $\F_2[X]$ and that they have
  no common factor. The first equation (\ref{eqn:F2:1}) shows that
  $Q_n$ divides all the $a_i$ with $1\leq i\leq n-1$. The second
  equation (\ref{eqn:F2:2}) shows that $Q_n$ also divides
  $a_0+a_n$. But then $Q_n^2$ divides the right hand side of
  (\ref{eqn:F2:2}) so $Q_n$ divides either $a_0$ or $a_n$. Since $Q_n$
  also divides $a_0+a_n$, we have shown that $Q_n$ divides all the
  $a_i$ with $i=0,\ldots,n$. This contradicts our assumption that the
  $a_i$ do not have a common factor, proving our claim.
\end{proof}}

\section{Applications}
\label{sect:applications}
\subsection{Fundamental groups}
In section \ref{sect:factor}, we have shown the following. Let
$\Lambda$ be a countable abelian group. Suppose $\Lambda$ acts
ergodically and measure preservingly on an infinite measure space
$(Y,\nu)$. Then there is a type II$_1$ factor $M_\Lambda$ with
fundamental group
$\Fundg(M_\Lambda)=\Mod(\Centr_{\Autns(Y)}(\Lambda))$. Using proposition
\ref{prop:ex:general}, we can generalize this result.

With $\Lambda=\FG_\infty$, proposition \ref{prop:ex:general} yields
and action $\Gamma\actson (X,\mu)$ and a quotient
$\pi:\Gamma\rightarrow\FG_\infty$ that satisfy the conditions of
theorem \ref{thm:intro:main:factor}. Let $\FG_\infty$ act
ergodically (but not necessarily freely) on an infinite measure space
$(Y,\nu)$. Then we find a type II$_1$ factor $M$ with fundamental
group
\[\Fundg(M)=\Mod(\Centr_{\Autns(Y)}(\FG_\infty)).\]

We give a slightly more appealing characterisation of the groups that
appear in the right hand side. Let $(Y,\nu)$ be an infinite measure
space and consider a closed subgroup $\cG$ of measure preserving
transformations on $Y$. Assume that $\cG$ acts ergodically on
$Y$. Since $\Autmp(Y,\nu)$ is a Polish group, so is
$\cG$ and hence we can take a group morphism
$\rho:\FG_\infty\rightarrow \cG$ that has dense range in $\cG$.
This defines an ergodic measure preserving action of $\FG_\infty$ with
\[\Centr_{\Autns(Y)}(\FG_\infty)=\Centr_{\Autns(Y)}(\cG).\]

We have proven the following theorem.
\begin{theorem}
  \label{thm:fundg:cartan}
  Let $(Y,\nu)$ be a standard infinite measure space. Let
  $\cG\subset\Autmp(Y,\nu)$ be a closed subgroup of measure preserving
  transformations of $Y$. Assume $\cG$ acts ergodically on $Y$. Then
  there is a type II$_1$ equivalence relation $\RelR_{\cG}$ and a type
  II$_1$ factor $M_{\cG}$ with fundamental group
  \[\Fundg(M_{\cG})=\Fundg(\RelR_{\cG})=\Mod(\Centr_{\Autns(Y)}(\cG)).\]
\end{theorem}

Using the results of section \ref{sect:noncartan}, we can further
generalize theorem \ref{thm:fundg:cartan} to theorem
\ref{thm:fundg:noncartan} below. However, we can not give an explicit
characterisation in full generality. The problem is the third
condition of theorem \ref{thm:prescartan}, which requires that $\ker\pi$
contains a property (T) subgroup $G$ such that there are no
$\ast$-homomorphisms $\theta:\Lg(G)\rightarrow p B p$ for any projection
$p\in B$ with finite trace. For many choices of $B$, any property (T)
group $G$ will suffice. For examples if $B$ is the hyperfinite II$_1$
factor, or if $B=\Bounded(\ell^2(\IN))\otimes \Lg(\FG_\infty)$.
For an arbitrary $B$, we use Ozawa's result from \cite{Ozawa:NoUniversalII1Factor}.

\begin{theorem}[see also theorem \ref{thm:intro:fundg:noncartan}]
  \label{thm:fundg:noncartan}
  Let $(B,\Tr)$ be any separable properly infinite but semifinite von
  Neumann algebra, with a given semifinite trace $\Tr$. Let
  $\Lambda\subset\Outtp(B,\Tr)$ be a countable group of trace
  preserving automorphisms of $B$. Suppose that $\Lambda$ acts
  ergodically on the center of $B$. Then there is a type II$_1$ factor
  $M_\Lambda$ with fundamental group
  \[\Fundg(M_\Lambda)=\Mod(\Centr_{\Out(B)}(\Lambda)),\]
  where $\Centr_{\Out(B)}(\Lambda)$ is the group of all outer
  automorphisms of $B$ that commute with $\Lambda$, as outer
  automorphisms.
\end{theorem}
Conversely, for every separable II$_1$ factor $M$, the fundamental
group is of this form:
\[\Fundg(M)=\Mod(\Autf(\Bounded(\ell^2(\IN))\otimes
M))=\Mod(\Centr_{\Out(\Bounded(\ell^2(\IN))\otimes M)}(\Lambda)).\]
\begin{proof}
Let $(B,\Tr)$ and $\Lambda$ be as in the statement.\\
\emph{Step 1: we can assume that there is no finite $\Lambda$-invariant
  measure on the center of $B$.}\\
Suppose there was such a finite measure. We replace $B$ by
$B_1=\ell^\infty(\IN)\vnOtimes B$ and $\Lambda_0$ by
$\Lambda_1=\Perm_\infty\times\Lambda_0$ where $\Perm_\infty$ is the group of
finite permutations of $\IN$. The group morphism
\[\theta:\Centr_{\Out(B)}(\Lambda)\ni\psi\mapsto
\id\otimes\psi\in\Centr_{\Out(B_1)}(\Lambda_1)\]
is an automorphism.

\emph{Step 2: there is genuine action $\beta:\FG_\infty\actson B$ such that
\[\Mod(\CCentr_B(\FG_\infty))=\Mod(\Centr_{\Out(B)}(\Lambda)).\]}
Take a quotient morphism
$\rho:\FG_\infty\rightarrow\Lambda\subset\Out(B)$. For every
elementary generator $a_n$ of $\FG_\infty$, choose a lift
$\beta(a_n)\in\Autf(B)$ for $\rho(a_n)\in\Out(B)$. This choice extends
to a group morphism $\beta:\FG_\infty\rightarrow \Autf(B)$, that is
necessarily a lift of the morphism $\rho:\FG_\infty\rightarrow\Out(B)$.
It is clear that
$\CCentr_{\Out(B)}(\FG_\infty)\subset\Centr_{\Out(B)}(\Lambda)$. Let $\psi$
be an automorphism of $B$ that commutes with $\Lambda$ up to
inner automorphisms. For every elementary generator
$a_n\in\FG_\infty$, choose a unitary $v_{a_n}\in B$ such that
$\psi\circ\beta_{a_n}=\Ad_{v_{a_n}}\circ\beta_{a_n}\circ\psi$.
This choice extends uniquely to a cocycle
$(v_s)_{s\in\FG_\infty}$ such that $\psi$ commutes with $\FG_\infty$
up to $(v_s)_s$. This shows that at least
\[\Mod(\CCentr_{B}(\FG_\infty))=\Mod(\Centr_{\Out(B)}(\Lambda)).\] 

\emph{Step 3: construction of $M_\Lambda$}
By \cite[Theorem
2]{Ozawa:NoUniversalII1Factor}, we find a property (T) group $G$ such
that there are no $\ast$-homomorphisms $\Lg(G)\rightarrow
qBq$ for any projection $q\in B$ with
finite trace. We can moreover assume that $G$ is ICC and quasifinite, meaning
that all the proper subgroups of $G$ are finite.
(see \cite[Theorem 1]{Ozawa:NoUniversalII1Factor} which is a variant
of \cite{Olshanskii:OnResidualingHomomorphisms} and
\cite{Gromov:HyperbolicGroups}) In particular, $G$ is a non-abelian
simple group and therefore a perfect group.

Set $\widetilde\Lambda=G\times\FG_\infty$ and construct a group
$\Gamma$, a quotient $\tilde\pi:\Gamma\rightarrow\widetilde\Lambda$
and an action $\alpha:\Gamma\actson (X,\mu)$ as in construction
\ref{constr:ex:general}. Consider the obvious quotient
$\pi:\Gamma\rightarrow\FG_\infty$. Define a new action $\sigma$ of
$\Gamma$ on $\Lp^\infty(X)\vnOtimes B$ by the formula
$\sigma_g(a\otimes b)=\alpha_g(a)\otimes \beta_{\pi(g)}(b)$. Denote
the crossed product by $N=(\Lp^\infty(X)\vnOtimes
B)\rtimes\Gamma$. Take any projection $p\in N$ with finite trace and
set $M_\Lambda=pNp$.

\emph{Step 3: the fundamental group of $M_\Lambda$ is \[\Fundg(M_\Lambda)=\Mod(\Centr_{\Out(B)}(\Lambda)).\]}
It is sufficient to show that the action $\Gamma\actson X$ and the
quotient $\pi:\Gamma\rightarrow \FG_\infty$ from step 2, satisfy the
conditions of theorem \ref{thm:main:noncartan}.

Proposition \ref{prop:ex:general} shows that $\ker\tilde\pi$
acts weakly mixingly and cocycle superrigidly on $X$. By
\cite[proposition 3.6]{Popa:CocycleSuperrigidityMalleable} the same is true for $\ker\pi$. Since $\FG_\infty$ has
the Haagerup property, there are no non-trivial group morphisms
$\theta:\ker\pi\rightarrow \FG_\infty$. Proposition
\ref{prop:ex:general} shows that
$\Norm_{\Autmp(X)}(\ker\tilde\pi)=\Gamma$, but we have to show that
$\Norm_{\Autmp(X)}(\ker\pi)=\Gamma$. It suffices to show that every
automorphism $\delta$ of $\ker\pi$ maps $\ker\tilde\pi$ onto
itself. Remember that $\ker\pi=\ker\tilde\pi\times G$. Since $G$ is an
infinite simple group, we know that there are
no non-trivial group morphisms from $G$ to any finite group. Because
$\Gamma_1=\SL_3{(\F_2[X])}\ltimes\F_2[X]^3$ is residually finite, there
are no non-trivial group morphism $\theta:G\rightarrow\Gamma_1$. There
can not be any non-trivial group morphisms
$\theta:G\rightarrow\ker\tilde\pi$ because $G$ has property (T) and
$\ker\tilde\pi$ is an (infinite) amalgamated free product of copies of
$\Gamma_1$. So our automorphism $\delta:\ker\pi\rightarrow\ker\pi$
maps $G$ into $G$. By symmetry, $\delta$ maps $G$ isomorphically onto
$G$. Since $G$ has trivial center, we see that $\delta$ maps
$\ker\tilde\pi=\Centr_{\ker\pi}(G)$ onto itself.

The conditions of theorem \ref{thm:prescartan} are satisfied by
construction. Since $\ker\pi$ is generated by copies of the perfect
groups $\Gamma_1$ and $G$, we see that $\ker\pi$ is a perfect group.
Theorem \ref{thm:main:noncartan} implies that the outer automorphism
group of $M_\Lambda$ is given by
\[\Out(M_\Lambda)=\OCCentr_{B}(\FG_\infty).\]
By step 2, we see that
\[\Fundg(M_\Lambda)=\Mod(\Centr_{\Out(B)}(\Lambda)).\]
\end{proof}

\subsection{Outer automorphism groups}
Theorem \ref{thm:intro:main:eqrel} gives us a way to compute the
outer automorphism group of type II$_\infty$ equivalence relations
$\RelR$ on an infinite measure space $(Y,\nu)$. For any set $U\subset
Y$ with finite measure, we have the short exact sequence
\[\begin{CD}
  1 @>>> \Out(\RelR\restrict U) @>>> \Out(\RelR) @>\Mod>>
  \Fundg(\RelR\restrict U) @>>> 1
\end{CD}.\]
Once we know $\Out(\RelR)$ and the group morphism $\Mod$, we can compute the
outer automorphism group and the fundamental group of the type II$_1$
equivalence relation $\RelR\restrict U$. A similar short exact
sequence exists for the fundamental group and the outer automorphism
group of type II$_1$ factors.

Up to now, we have only computed fundamental groups of type II$_1$
equivalence relations and factors. Theorems \ref{thm:intro:main:eqrel}
and \ref{thm:intro:main:factor}
allow us to also compute outer automorphism groups of type II$_1$
equivalence relations and factors.

More concretely, let $\Lambda\actson (Y,\nu)$, $\Gamma\actson (X,\mu)$
and $\pi:\Gamma\rightarrow\Lambda$ be as in theorem
\ref{thm:intro:main:eqrel}, and construct the II$_\infty$ relation
$\RelR$ as in this same theorem. Theorem \ref{thm:intro:main:eqrel}
shows that the outer
automorphism group of $\RelR$ is $\Out(\RelR)=\Centr_{\Autns(Y)}(\Lambda)$. But
$\RelR$ is a type II$_\infty$ equivalence relation. Take any subset
$U\subset X\times Y$ with finite measure. Then the outer automorphism
group of the type II$_1$ equivalence relation $\RelR\restrict{U}$ is
\[\Out(\RelR\restrict{U})=\Centr_{\Autmp(Y,\nu)}(\Lambda),\]
where $\Autmp(Y,\nu)$ is the set of all $\nu$-preserving automorphisms
of $Y$.

The example from construction \ref{constr:ex:general} satisfies the
conditions of theorem \ref{thm:intro:main:eqrel}. This proves the
following result.
\begin{theorem}
\label{thm:ex:out:eqrel}
For any closed subgroup $\cG\subset \Autmp(Y,\nu)$ of
$\nu$-preserving transformations that acts ergodically on a standard measure
space (finite or infinite), there is a type II$_1$ equivalence
relation $\RelR_{\cG}$ with outer automorphism group
\[\Out(\RelR_{\cG})=\Centr_{\Autmp(Y,\nu)}(\cG).\]
\end{theorem}
For any locally compact, second countable, unimodular group $G$, this
yields a type II$_1$ equivalence relation $\RelR$ with outer
automorphism group $\Out(\RelR)=G$: let $G$ be such a locally compact
second countable group unimodular group. Consider the measure space
$(G,h)$ where $h$ is the Haar measure on $G$. By definition, the left
translation action of $\cG=G$ on $(G,h)$ is measure preserving. The
centralizer of this action is precisely the right action of $G$ on
$(G,h)$. A similar result was obtained in \cite[theorem
1.1]{PopaVaes:ActionsOfFinfty}.
\begin{proof}[Proof of theorem \ref{thm:ex:out:eqrel}]
Let $\Lambda_0\subset\cG$ be a countable dense subgroup, and take a
quotient morphism $\rho:\FG_\infty\rightarrow \Lambda_0$. This gives
us an ergodic measure preserving action $\FG_\infty\actson (Y,\nu)$
with
$\Centr_{\Autmp(Y,\nu)}(\FG_\infty)=\Centr_{\Autmp(Y,\nu)}(\cG)$. Proposition
\ref{prop:ex:general} and theorem \ref{thm:intro:main:eqrel} show that
there is a type II$_1$ equivalence relation $\RelR_{\cG}$ with
\[\Out(\RelR_{\cG})=\Centr_{\Autmp(Y,\nu)}(\cG).\]
\end{proof}

Using exactly the same argument as above, we find a type II$_1$ factor
$M_{\cG}$ with outer automorphism group
\[\Out(M_{\cG})=\Centr_{\Autmp(Y,\nu)}(\cG)\ltimes\Cohom(\FG_\infty\actson
Y),\]
where the action of $\FG_\infty$ is given by a group morphism
$\theta:\FG_\infty\rightarrow \cG$ with dense range. This result has
not the same appeal as theorem \ref{thm:ex:out:eqrel} because the
group $\Cohom(\FG_\infty\actson Y)$ is huge. In two special cases we
can work around this problem. Theorem \ref{thm:ex:out:SLnR} below gives a
type II$_1$ factor $M$ with outer automorphism group
$\Out(M)=\SL^{\pm}_n\IR=\{g\in\GL_n\IR\mid \det(g)=\pm1\}$.

Besides that, theorem \ref{thm:ex:out:factor} shows that for every closed
subgroup $\cG$ of \emph{probability} measure preserving
transformations on $(Y,\nu)$, we find a II$_1$ factor $M_{\cG}$ with
outer automorphism group $\Out(M_{\cG})=\Centr_{\Autmp(Y,\nu)}(\cG)$.
For any compact group $G$, this yields an explicit construction of a
type II$_1$ factor with outer automorphism group equal to $G$. The
existence of such II$_1$ factors was already proven in \cite{FalguieresVaes:CompactOut}.
But the class of groups of the
form $\Centr_{\Autmp(Y,\nu)}(\cG)$ contains more than just the compact
groups. For example, we can take $\cG=\{\id\}$. Theorem
\ref{thm:ex:out:factor} shows that there is
a type II$_1$ factor whose outer automorphism group is precisely
$\Autmp(Y,\nu)$, the group of probability measure preserving
transformations on a standard probability space.

\begin{theorem}
  \label{thm:ex:out:SLnR}
  For every natural number $0\not=n\in\IN$ there is a type II$_1$ factor $M$
  whose outer automorphism group is
  \[\Out(M)=\SL^{\pm}_n\IR.\]
\end{theorem}
\begin{proof}
  Take $m=4n+1$, and consider the action of $\Lambda=\SL_m\IZ$ on
  $Y=\MatM_{m,n}(\IR)$ by left multiplication. This action clearly
  preserves the Lebesgue measure $\nu$ on $Y$, and its centralizer is
  \[\Centr_{\Autmp(Y,\nu)}(\Lambda)=\SL^{\pm}_n\IR.\]
  By \cite[theorem 1.3]{PopaVaes:lattices}, we know that the action
  $\Lambda\actson Y$ is $\Ufin$ cocycle superrigid. Since $\SL_m\IZ$
  is a perfect group, it follows that $\Cohom(\Lambda\actson
  Y)=\{1\}$.

  Choose any $k\geq m+1$ that is also a multiple of $3$. Using
  construction \ref{constr:ex:general}, theorem
  \ref{thm:intro:main:eqrel} and theorem \ref{thm:intro:main:factor},
  we construct a type II$_1$ factor $M$ with outer automorphism group
  \[\Out(M)=\SL^{\pm}_n\IR.\]
\end{proof}

\begin{theorem}
  \label{thm:ex:out:factor}
  Let $(Y,\nu)$ be a \emph{probability} space, and let $\cG\subset
  \Autmp(Y,\nu)$ be a closed subgroup of measure preserving
  transformations. Then there is a type
  II$_1$ factor $M_{\cG}$ with outer automorphism group
  $\Out(M_{\cG})=\Centr_{\Autmp(Y,\nu)}(\cG)$.
\end{theorem}
\begin{proof}
  This proof of theorem \ref{thm:ex:out:factor} uses a generalization
  of the co-induced action. The construction is explained in
  definition \ref{def:gen:co-ind} below and some properties are
  given in lemma \ref{lem:gen:co-ind} below.

  Let $\rho:\FG_\infty\rightarrow\cG$ be any group morphism with dense
  range. Denote by $a_n$ the $n$-th elementary generator of
  $\FG_\infty$.
  Consider the group
  \begin{equation}
    \label{eqn:Lambda}
    \Lambda=\underset{\overset{\Vert}{G}}{\SL_3\IZ}\ltimes(\underset{\overset{\Vert}{A}}{\IZ^3}\times((\underbrace{\IZ^3\times\IZ^3}_{B})\free
    \underset{\overset{\Vert}{F_1}}{\IZ^3}\free\underset{\overset{\Vert}{F_2}}{\IZ^3}\free\ldots)),
  \end{equation}
  where $\SL_3\IZ$ acts in the obvious way on each of the copies of
  $\IZ^3$. Because each copy of $\IZ^3$ plays a slightly different
  role in the following, we give them different names, as indicated in
  (\ref{eqn:Lambda}). Consider the obvious action of
  $\SL_3\IZ\ltimes(\IZ^3\times\IZ^3\times\IZ^3)=G\ltimes(A\times B)$
  on $I=\IZ^3\times\IZ^3\times\IZ^3$. This action extends to an action
  of $\Lambda$ on $I$, where the $F_n$ act trivially.

  We define a cocycle $\omega:\Lambda\times I\rightarrow \FG_\infty$
  by the following relations.
  \begin{align*}
    \omega(g,i)=e&\quad
    \text{ for all }g\in G\ltimes(A\times B)\text{ and all }i\in I\\
    \omega(f,(i_1,i_2,i_3))=a_n^{\det(f,i_2,i_3)}&\quad
    \text{ for all }f\in F_n\text{ and all }(i_1,i_2,i_3)\in I=\IZ^3\times\IZ^3\times\IZ^3.
  \end{align*}
  In the last formula, we denoted $\det(f,i_2,i_3)\in\IZ$ for the
  determinant of the matrix whose columns are $f,i_2$ and $i_3$.

  Consider the generalized co-induced action of $\FG_\infty\actson (Y,\nu)$,
  associated with $\omega:\Lambda\times I\rightarrow \FG_\infty$, as
  explained in definition \ref{def:gen:co-ind}. Denote the resulting
  action by $\Lambda\actson (Y_1,\nu_1)$. Lemma \ref{lem:gen:co-ind}
  shows that this action is ergodic and that
  \[\Centr_{\Autns(Y_1)}(\Lambda)=\Centr_{\Autns(Y)}(\cG).\]

  Following construction \ref{constr:ex:general}, we construct a group
  $\Gamma$, an action $\Gamma\actson (X,\mu)$ and
  a quotient $\pi:\Gamma\rightarrow\Lambda$. Consider the crossed
  product $M=\Lp^\infty(X\times Y_1)\rtimes\Gamma$.

  Observe that $\ker\pi$ is generated by copies of
  $\SL_3{(\F_2[X])}\ltimes\F_2[X]^3$. The group $\SL_3{(\F_2[X])}\ltimes\F_2[X]^3$ is generated by
  conjugations of the subgroup $\SL_3{(\F_2)}$, because
  $(1+P(X)e_{i,j},Q(X)e_i)=[(1+P(X)e_{k,j},Q(X)e_k),(1+e_{i,k},0)]$ whenever $i,j,k$ are
  different elements in $\{1,2,3\}$. Lemma \ref{lem:SL32} below shows
  that there are no non-trivial group morphisms
  $\theta:\SL_3{(\F_2)}\rightarrow\SL_3\IZ$. Hence there are none from
  $\ker\pi$ into $\SL_3\IZ$. Because $\ker\pi$ is a perfect group,
  there are no non-trivial group morphisms $\theta:\ker\pi\rightarrow \Lambda$.

  Proposition \ref{prop:ex:general} and theorem \ref{thm:proof:main:factor}
  show that
  \[\Out(M)=\Centr_{\Autmp(Y_1,\nu_1)}(\Lambda)\ltimes\Cohom(\Lambda\actson
  Y_1).\]

  We show that the cohomology group $\Cohom(\Lambda\actson Y_1)$ is
  trivial. Let $c:\Lambda\times Y_1\rightarrow\Circle^1$ be a
  cocycle. By Popa's cocycle superrigidity theorem \cite[theorem
  0.1]{Popa:CocycleSuperrigidityMalleable}, we can assume that $c(g,y)$ does not depend on $y$
  whenever $g\in A$. But $A\subset\Lambda$ is a normal subgroup that
  acts weakly mixingly on $Y_1$. Hence \cite[proposition 3.6]{Popa:CocycleSuperrigidityMalleable}
  shows that $c(g,y)$ is independent of $y$ for any $g\in\Lambda$, or
  still $c:\Lambda\rightarrow\Circle^1$ is a group morphism. This
  group morphism is trivial because $\Lambda$ is a perfect group.

  We have shown that $\Cohom(\Lambda\actson Y_1)$ is trivial and hence
  that
  \[\Out(M)=\Centr_{\Autmp(Y_1,\nu_1)}.\]
\end{proof}

\begin{lemma}
  \label{lem:SL32}
  There are no non-trivial group morphisms
  $\theta:\SL_3{(\F_2)}\rightarrow\SL_3\IZ$.
\end{lemma}
\begin{proof}
  Suppose that $\theta:\SL_3{(\F_2)}\rightarrow\SL_3\IZ$ is a group
  morphism. We can consider this morphism as a non-trivial
  representation $\theta:\SL_3{(\F_2)}\rightarrow\GL_3\IC$. The
  character of this representation is defined as
  $\chi_\theta(g)=\tr(\theta(g))$ and hence $\chi_\theta$ takes values
  in $\IZ$.
  According to the atlas of finite groups\cite{AtlasOfFiniteGroups},
  the character table of $\SL_3{(\F_2)}$ is the following:\\
  \centerline{
    \begin{tabular}{cccccc}
      \sMatrix{1&0&0\\0&1&0\\0&0&1}&\sMatrix{1&0&1\\0&1&0\\0&0&1}&
      \sMatrix{0&1&0\\0&0&1\\1&0&0}&\sMatrix{1&1&0\\0&1&1\\0&0&1}&
      \sMatrix{0&1&1\\1&1&0\\1&0&0}&\sMatrix{1&1&1\\1&1&0\\1&0&0}\\
      1&1&1&1&1&1\\
      3&-1&0&1&$r_7+r_7^2+r_7^4$&$r_7^3+r_7^5+r_7^6$\\
      3&-1&0&1&$r_7^3+r_7^5+r_7^6$&$r_7+r_7^2+r_7^4$\\
      6&2&0&0&-1&-1\\
      7&-1&1&-1&0&0\\
      8&0&-1&0&1&1\\
    \end{tabular}
  }\\
  where $r_7$ denotes a primitive 7-th root of unity. From this table,
  we see that all non-trivial characters of dimension 3 take at least
  one non-integer value. So $\theta:\SL_3{(\F_2)}\rightarrow\SL_3\IZ$
  is trivial.
\end{proof}

Let $\Lambda_0\subset\Lambda$ be an inclusion of groups, and suppose
that $\Lambda_0$ acts probability measure preservingly on a space $(Y_0,\nu_0)$.
One possible construction for the co-induced action of $\Lambda$
associated with $\Lambda_0\actson (Y_0,\nu_0)$ is the
following. Set $I=\Lambda/\Lambda_0$ and consider the space
$(Y,\nu)=(Y_0,\nu_0)^{I}$. The group $\Lambda$ acts in a natural way
on $I$. Choose representatives $g_i\in\Lambda$ for
the cosets of $\Lambda_0$. Then we define a cocycle
$\omega:\Lambda\times I\rightarrow\Lambda_0$ by the formula
$gg_i=g_{gi}\omega(g,i)$. The co-induced action of $\Lambda$ on
$(Y,\nu)$ is given by $(g y)_i=\omega(g,g^{-1}i)y_{g^{-1}i}$. Up to
conjugation, this action does not depend on the choice of
representatives $g_i\in\Lambda$.

We generalize this construction.
\begin{definition}
  \label{def:gen:co-ind}
  Let $\Lambda\actson I$ be an action of a countable group on a
  countable set. Let $\omega:\Lambda\times I\rightarrow \Lambda_0$ be
  a cocycle. Suppose that $\Lambda_0$ acts probability measure
  preservingly on $(Y_0,\nu_0)$. Define an action of $\Lambda$ on
  $(Y,\nu)=(Y_0,\nu_0)^I$ by the formula
  $(g y)_i=\omega(g,g^{-1}i)y_i$.
  This action is called the generalized co-induced action of
  $\Lambda_0\actson (Y_0,\nu_0)$, with respect to $\omega$.
\end{definition}
\begin{lemma}
  \label{lem:gen:co-ind}
  Let $\Lambda\actson I$ be an action of a countable group on a
  countable set, and let $\omega:\Lambda\times I\rightarrow\Lambda_0$
  be a cocycle. Suppose that $\Lambda_0\actson (Y_0,\nu_0)$ is a
  probability measure preserving action. Consider the generalized
  co-induced action $\Lambda\actson (Y,\nu)=(Y_0,\nu_0)^I$.
  \begin{itemize}
  \item If all orbits of $\Lambda\actson I$ are infinite, then the
    generalized co-induced action $\Lambda\actson(Y,\nu)$ is weakly mixing.
  \item Suppose that $\Lambda\actson I$ and $\omega$ satisfy the
    following three conditions.
    \begin{itemize}
    \item $\Lambda$ acts transitively on $I$.
    \item There exists an $i\in I$ such that (or equivalently, for all
      $i\in I$) $\omega$ maps the set
      $\Stab\{i\}\times\{i\}$ surjectively onto $\Lambda_0$
    \item There exists an $i\in I$ such that (or equivalently, for all
      $i\in I$) the subgroup $S_i=\{g\in\Lambda\mid gi=i\text{ and
      }\omega(g,i)=e\}$ acts with infinite orbits on $I\setminus\{i\}$.
    \end{itemize}
    Then we have that
    \[\Centr_{\Autmp(Y,\nu)}(\Lambda)=\Centr_{\Autmp(Y_0,\nu_0)}(\Lambda_0).\]
  \end{itemize}
\end{lemma}
\begin{proof}
  The first point can be proven exactly as the analogous result for
  generalized Bernoulli actions, see for example \cite[proposition
  2.3]{PopaVaes:StrongRigidity}.

  For every $\Delta_0\in\Centr_{\Autmp(Y_0,\nu_0)}(\Lambda_0)$, we
  define $\Delta\in\Centr_{\Autmp(Y,\nu)}(\Lambda)$ by the formula
  $\Delta(y)_i = \Delta_0(y_i)$. It remains to show that every
  automorphism $\Delta\in\Centr_{\Autmp(Y,\nu)}(\Lambda)$ is of this
  form. Let $\Delta:Y\rightarrow Y$ be a measure preserving
  automorphism that commutes with the action of $\Lambda$. Exactly the
  same argument as in the proof of \cite[proposition 6.10]{Vaes:Bimodules}
  yields p.m.p.\ automorphisms $\Delta_i:Y_0\rightarrow Y_0$ and a
  bijection $\alpha:I\rightarrow I$ such that that $Delta$ is given by
  \begin{equation}
    \label{eqn:Delta:form}
    (\Delta(y))_i=\Delta_i(y_{\alpha^{-1}(i)}).
  \end{equation}

  It follows that $\alpha(gi)=g\alpha(i)$ for all $g\in\Lambda$ and
  $i\in I$. For fixed $i\in I$, we see that $\alpha(i)$ is fixed under the
  action of $S_i$. But $S_i$ acts with infinite orbits on
  $I\setminus\{i\}$, so $\alpha(i)=i$. This works for all $i\in I$ so
  we see that $\alpha=\id$.

  From (\ref{eqn:Delta:form}) we also see that
  $\Delta_{gi}(\omega(g,i)y)=\omega(g,i)\Delta_{i}(y)$ for almost all
  $y\in Y_0$, and for all $g\in\Lambda$, $i\in I$. Fix $i\in I$. Since
  $\omega(\Stab\{i\}\times\{i\})=\Lambda_0$, we see that $\Delta_i$
  commutes with $\Lambda_0$, or still,
  $\Delta_i\in\Centr_{\Autmp(Y_0,\nu_0)}(\Lambda)$. Moreover, we see
  that $\Delta_{gi}=\Delta_i$ for all $g\in \Lambda$. Because
  $\Lambda$ acts transitively on $I$, it follows that
  $\Delta_j=\Delta_i$ for all $j\in I$. Writing $\Delta_0=\Delta_i$,
  we have shown that $\Delta$ is given by $\Delta(y)_i=\Delta_0(y_i)$.
\end{proof}

\newcommand{\etalchar}[1]{$^{#1}$}
\bibliography{references}{}
\bibliographystyle{sdpabbrv}
\end{document}

%% file: sdpStyle.inc
\usepackage{a4wide}
\usepackage{color}
\usepackage{pgf}
\usepackage[colorlinks, citecolor=green!50!black, linkcolor=red!50!black]{hyperref}
\usepackage{bookmark}
\usepackage{mathdots}
\parskip=\baselineskip

\let\oldlist=\list
\newlength\oldparskip
\def\list#1#2{\oldparskip=\parskip\parskip=0\baselineskip
\oldlist{#1}{#2}\parskip=0\baselineskip}
\let\oldendlist=\enditemize
\def\enditemize{\oldendlist\vspace*{-\oldparskip}}

\renewcommand{\address}[1]{\thanks{#1}}
\renewcommand{\email}[1]{\thanks{email: \texttt{#1}}}

\makeatletter
\def\@tocline#1#2#3#4#5#6#7{\relax
  \ifnum #1>\c@tocdepth 
  \else
    \begingroup \hyphenpenalty\@M
    \parskip=0cm
    \@ifempty{#4}{%
      \@tempdima\csname r@tocindent\number#1\endcsname\relax
    }{%
      \@tempdima#4\relax
    }%
    \parindent\z@ \leftskip#3\relax \advance\leftskip\@tempdima\relax
    \rightskip\@pnumwidth plus4em \parfillskip-\@pnumwidth
    #5\leavevmode\hskip-\@tempdima #6\nobreak\relax
    \hfil\hbox to\@pnumwidth{\@tocpagenum{#7}}\par
    \nobreak
    \endgroup
  \fi}
\makeatother

%% file: sdpDefinitions.inc
\usepackage{amsmath}
\usepackage{amssymb}
\usepackage{amsthm}
\usepackage{amscd}
\newtheorem{theorem}{Theorem}[section]
\newtheorem{lemma}[theorem]{Lemma}
\newtheorem{corollary}[theorem]{Corollary}
\newtheorem{proposition}[theorem]{Proposition}
\newtheorem{definition}[theorem]{Definition}

\newtheorem{definition_theorem}[theorem]{Definition/Theorem}

\newtheorem{construction}[theorem]{Construction}
\newcommand\len[1]{\left\lvert #1\right\rvert}
\newcommand\norm[1]{\left\lVert #1\right\rVert}

\DeclareMathOperator{\IC}{\mathbb{C}}
\DeclareMathOperator{\IR}{\mathbb{R}}
\DeclareMathOperator{\IRpos}{\mathbb{R}_{+}^{\ast}}
\DeclareMathOperator{\IQ}{\mathbb{Q}}
\DeclareMathOperator{\IZ}{\mathbb{Z}}
\DeclareMathOperator{\IN}{\mathbb{N}}
\DeclareMathOperator{\F}{F}
\DeclareMathOperator{\cG}{\mathcal{G}}
\DeclareMathOperator{\Circle}{S}
\DeclareMathOperator{\Perm}{Perm}
\def\MatrixGroup#1_#2#3{\mathchoice{#1_{#2}\hspace*{-0.3mm}#3}{#1_{#2}\hspace*{-0.4mm}#3}{#1_{#2}\hspace*{-0.4mm}#3\hspace*{0.3mm}}{#1_{#2}\hspace*{-0.4mm}#3\hspace*{0.3mm}}%
}
\newcommand\T{\MatrixGroup\OpT}
\newcommand\ST{\MatrixGroup\OpST}
\newcommand\SL{\MatrixGroup\OpSL}
\newcommand\GL{\MatrixGroup\OpGL}
\DeclareMathOperator{\OpGL}{GL}
\DeclareMathOperator{\OpT}{T}
\DeclareMathOperator{\OpST}{ST}
\DeclareMathOperator{\OpSL}{SL}
\newcommand{\FG}{\mathbb{F}}
\DeclareMathOperator{\free}{\ast}

\makeatletter
\def\quot{\@ifnextchar[{\sdp@quotleftA}{\sdp@quotright}}
\def\sdp@quotleftA[#1]#2{\@ifnextchar[{\sdp@quotboth[{#1}]{#2}}{\sdp@quotleft[{#1}]{#2}}}
\def\sdp@quotright#1[#2]{#1/#2}
\def\sdp@quotboth[#1]#2[#3]{#1\backslash #2 / #3}
\def\sdp@quotleft[#1]#2{#1\backslash #2}
\makeatother
\DeclareMathOperator{\Cohom}{H_1}
\DeclareMathOperator{\Centr}{Centr}
\DeclareMathOperator{\Norm}{Norm}

\DeclareMathOperator{\CCentr}{CCentr}
\DeclareMathOperator{\ICCentr}{InnCCentr}
\DeclareMathOperator{\OCCentr}{OutCCentr}
\DeclareMathOperator{\Stab}{Stab}
\DeclareMathOperator{\Fix}{Fix}
\DeclareMathOperator{\OpAut}{Aut}
\newcommand\Autmp{\OpAut_{\text{mp}}}
\newcommand\Autns{\OpAut_{\text{ns}}}

\newcommand\Outtp{\Out_{\text{tp}}}
\newcommand\Autf{\OpAut}
\newcommand\Autr{\OpAut}
\newcommand\Autg{\OpAut}
\newcommand\AutRing{\OpAut}
\DeclareMathOperator{\Mod}{mod}
\DeclareMathOperator{\Out}{Out}
\DeclareMathOperator{\Unitary}{\mathcal{U}}
\DeclareMathOperator{\Fundg}{\mathcal{F}}
\DeclareMathOperator{\cF}{\mathcal{F}}
\DeclareMathOperator{\Inn}{Inn}
\newcommand{\Ufin}{\mathcal{U}_{\text{fin}}}

\def\Rel(#1\actson #2){\mathcal{R}(#1\actson #2)}
\def\RelR{\mathcal{R}}

\def\fullg(#1){\left[#1\right]}
\def\fullpg(#1){\left[\left[#1\right]\right]}

\DeclareMathOperator{\Centre}{\mathcal{Z}}
\DeclareMathOperator{\vnNorm}{\mathcal{N}}
\DeclareMathOperator{\Bounded}{\mathcal{B}}
\DeclareMathOperator{\MatM}{M}
\DeclareMathOperator{\MatD}{D}
\DeclareMathOperator{\Lp}{L}
\DeclareMathOperator{\Lg}{L}
\DeclareMathOperator{\cM}{\mathcal{M}}
\DeclareMathOperator{\cN}{\mathcal{N}}
\DeclareMathOperator{\cP}{\mathcal{P}}

\DeclareMathOperator{\cA}{\mathcal{A}}
\DeclareMathOperator{\cB}{\mathcal{B}}
\DeclareMathOperator{\E}{E}
\DeclareMathOperator{\Tr}{Tr}
\DeclareMathOperator{\tr}{tr}
\DeclareMathOperator{\vnOtimes}{\overline{\otimes}}
\newcommand{\bimod}[3]{\strut_{#1}#2\strut_{#3}}
\DeclareMathOperator*{\embeds}{\prec}
\DeclareMathOperator*{\nembeds}{\not\prec}

\DeclareMathOperator*{\oembeds}{\succ}

\def\vnInterB(#1){\zH_{#1}}
\DeclareMathOperator{\diag}{diag}
\DeclareMathOperator{\Ad}{Ad}
\DeclareMathOperator{\actson}{\curvearrowright}
\DeclareMathOperator{\id}{id}

\DeclareMathOperator{\setminusb}{--\,}
\DeclareMathOperator{\setminusG}{--\,}

\DeclareMathOperator{\D}{d}
\newcommand\II{II}
\newcommand\restrict[1]{\vert_{#1}}
\makeatletter
\DeclareMathOperator{\zH}{H}
\def\smallscripts#1{\@ifnextchar_{\smallscripts@{#1}}{\@undefined}}
\def\smallscripts@#1{\@ifnextchar_{\smallscripts@subP{#1}}{\@ifnextchar^{\smallscripts@superP{#1}}{\smallscripts@none{#1}}}}
\def\smallscripts@subP#1_#2{\@ifnextchar^{\smallscripts@subsuperP{#1}{#2}}{\smallscripts@sub{#1}{#2}}}
\def\smallscripts@subsuperP#1#2^#3{\smallscripts@both{#1}{#2}{#3}}
\def\smallscripts@superP#1^#2{\@ifnextchar_{\smallscripts@supersubP{#1}{#2}}{\smallscripts@super{#1}{#2}}}
\def\smallscripts@supersubP#1#2_#3{\smallscripts@both{#1}{#3}{#2}}
\def\smallscripts@none#1{{\displaystyle #1} }
\def\smallscripts@sub#1#2{{\displaystyle #1}_{\scriptscriptstyle #2}}
\def\smallscripts@super#1#2{{\displaystyle #1}^{\scriptscriptstyle #2}}
\def\smallscripts@both#1#2#3{{\displaystyle #1}_{\scriptscriptstyle #2}^{\scriptscriptstyle #3}}
\makeatother
\newcommand\chara{\smallscripts{\chi}}
\newcommand\z{\smallscripts{z}}
\newcommand{\sMatrix}[1]{$\left(\!\begin{smallmatrix}#1\end{smallmatrix}\!\right)$}